\documentstyle[12pt]{article}
\begin{document}
\newtheorem{Def}{Definition}[section]
\newtheorem{theorem}{Theorem}[section]
\newtheorem{lemma}{Lemma}[section]
\newtheorem{remark}{Remark}[section]
\newtheorem{prop}{Proposition}[section]
\newtheorem{corollary}{Corollary}[section]
\newtheorem{claim}{Claim}[section]
\newtheorem{step}{Step}[section]
\newtheorem{example}{Example}[section]
\newtheorem{subsn}{Subsection}[section]
\title
{Some fully nonlinear problems on manifolds with boundary of
negative admissible curvature }
\author{\ Aobing Li\thanks{Partially supported by GRF Grant
6041477} \ \ \ \ \& \ \ \ Huan Zhu
\\ Department of Mathematics\\ City University of Hong Kong\\
83 Tat Chee Avenue\\
Kowloon, Hong Kong }
\date{}
\maketitle
\newcommand {\Bbb}{}

\setcounter{section}{0}

\begin{section}
{Introduction}
\end{section}

Let $(M^n,g)$ denote a compact smooth Riemannian manifold with no boundary of dimension $n\ge 3$. The Yamabe problem is to search a metric $\tilde {g}$ in the conformal class $[g]$ of $g$ such that $\tilde {g}$ has a constant scalar curvature $R_{\tilde {g}}=c$. Write $\tilde{g}=u^{\frac{4}{n-2}}g$. The Yamabe problem is equivalent to solve
\begin{equation}
-L_gu=cu^{\frac{n+2}{n-2}},\quad u>0\quad\mbox{in}~~M,
\label{yamabe}
\end{equation}
where $L_g=\Delta_g-\frac{n-2}{4(n-1)}R_g$ is the conformal Laplacian of $g$,  and $c=-1$, $0$, or $1$.

Let $\phi_1$ be a positive eigenfunction of the first eigenvalue $\lambda_1$ of $-L_g$, i.e.
\[
\lambda_1=\inf\limits_{\phi\in H^1(M)\setminus\{0\}}\frac{\int_M |\nabla\phi|_g^2+\frac{n-2}{4(n-1)}R_g\phi^2~dvol_g}{\int_M\phi^2~dvol_g},
\]
and $-L_g\phi_1=\lambda_1\phi_1$. A direct calculation yields that
\[
R_{\phi_1^{\frac{4}{n-2}}g}=-\frac{n-2}{4(n-1)}\phi_1^{\frac{-n-2}{n-2}}L_g\phi_1=\frac{n-2}{4(n-1)}\phi_1^{\frac{-4}{n-2}}\lambda_1.
\]

After replacing $g$ by $\phi_1^{\frac{4}{n-2}}g$, we assume the scalar curvature of the background metric $g$ has a definite sign, that is, either

\[
R_g>0,\quad\mbox{or}\quad R_g\equiv 0,\quad\mbox{or}\quad R_g<0.
\]

Consider the functional
\[
{\cal Q}(\phi)=\frac{\int_M |\nabla\phi|_g^2+\frac{n-2}{4(n-1)}R_g\phi^2~dvol_g}{(\int_M\phi^{\frac{2n}{n-2}}~dvol_g)^{\frac{n-2}{n}}}.
\]
$u$ is a solution of the equation (\ref{yamabe}), then $u$ is a critical point of the above functional ${\cal Q}$. It is a simple consequence of the H$\ddot{o}$lder inequality that
\[
\lambda(M^n,g):=\inf\limits_{\phi\in H^1(M)\setminus\{0\}}{\cal Q}(\phi)>-\infty.
\]

In \cite{Ya}, Yamabe approached the problem by attempting to prove that a minimizing sequence of ${\cal Q}$ will converge to a minimizer. Trudinger (\cite{Tr}) pointed out that the convergence failed on the standard sphere $(S^n, g_{round})$, and Trudinger was able to fix Yamabe's arguments when $\lambda(M^n,g)\le 0$. In general, we know
\[
\lambda(M^n,g)\le \lambda(S^n,g_{round}).
\]
In(\cite{Aub}), Aubin proved the convergence of the minimizing sequence if
\[
\lambda(M^n,g)< \lambda(S^n,g_{round}).
\]
When the manifold $M^n$ is not locally conformally flat, it was proved by Aubin, for $n\ge 6$, and that by Schoen, for $n=3,4,5$, that $\lambda(M^n,g)< \lambda(S^n,g_{round})$. When the manifold is locally conformally flat and not conformally diffeomorphic to the standard sphere, Schoen established the compactness result of the solutions to the equation (\ref{yamabe}) using a deep result of his joint work with Yau in \cite{SY}, therefore confirmed the existence of the solutions.

For $(M^n,g)$, an $n-$dimensional($n\ge 3$) smooth Riemannian compact manifold with boundary, a similar problem is to look for a metric $\tilde {g}\in [g]$ having constant scalar curvature on $M^n$ and constant mean curvature on the boundary $\partial M$. Let $\tilde {g}=u^{\frac{4}{n-2}}g$. The problem is equivalent to searching a solution of the following equation
\begin{equation}
\{
\begin{array}{lcl}
-L_gu&=&c_1u^{\frac{n+2}{n-2}}\quad\mbox{on}~~M^n\\
B_g u&=&c_2 u^{\frac{n}{n-2}}\quad\mbox{on}~~\partial M,
\end{array}
\label{yamabebdry}
\end{equation}
where the boundary operator $B_g=\frac{2}{n-2}\frac{\partial}{\partial\nu}+h_g$, $h_g$ is the mean curvature of $g$ w.r.t. the unit outer normal $\frac{\partial}{\partial\nu}$, and $c_1$, $c_2$ denote two constants. When $c_2=0$, the problem is variational. In fact, the equation (\ref{yamabebdry}) is the Euler-Lagrange equation of the functional

\[
{\cal F}(\phi)=\frac{\int_M |\nabla\phi|_g^2+\frac{n-2}{4(n-1)}R_g\phi^2~dvol_g+\frac{n-2}{2}\int_{\partial M}h_g\phi^2~dS_g}{(\int_M\phi^{\frac{2n}{n-2}}~dvol_g)^{\frac{n-2}{n}}},
\]
and we have
\[
\lambda(M,g):=\inf\limits_{\phi\in H^1(M)\setminus\{0\}}{\cal F}(\phi)>-\infty.
\]

Cherrier (\cite{Cher}) proved that the $\inf\cal{F}$ is achieved by a smooth positive function if

\begin{equation}
\lambda(M,g)<\lambda(S^n_+,g_{round}),
\label{bdrycher}
\end{equation}
where $(S^n_+,g_{round})$ is the standard half sphere. When $c_2=0$ in the equation (\ref{yamabebdry}), Escobar (\cite{Esc}) obtained the existence of the solution for a large class of manifolds by achieving (\ref{bdrycher}). For the general constant $c_2$, let $\phi_1$ be a smooth positive function of the eigenvalue problem
\[
\{
\begin{array}{lcl}
-L_g\phi_1&=&\lambda_1 \phi_1\quad\mbox{on}~~M^n\\
B_g \phi_1&=&0 \quad\mbox{on}~~\partial M,
\end{array}
\]
where
\begin{equation}
\lambda_1:=\inf\limits_{\phi\in H^1(M)\setminus\{0\}}\frac{\int_M |\nabla\phi|_g^2+\frac{n-2}{4(n-1)}R_g\phi^2~dvol_g+\frac{n-2}{2}\int_{\partial M}h_g\phi^2~dS_g}{\int_M\phi^2~dvol_g}.
\label{type}
\end{equation}
Then
\[
\{
\begin{array}{lcl}
R_{\phi_1^{\frac{4}{n-2}}g}&=&\frac{4(n-1)}{n-2}\lambda_1 \phi_1^{\frac{-4}{n-2}}\quad\mbox{on}~~M^n\\
h_{ \phi_1^{\frac{4}{n-2}}g}&=&0 \quad\mbox{on}~~\partial M.
\end{array}
\]
Replacing $g$ by $ \phi_1^{\frac{4}{n-2}}g$, we may assume one of the following three cases holds, i.e.,
\[
\begin{array}{lclcl}
R_g>0,&&R_g<0,&&R_g=0.\\
&\mbox{or}&&\mbox{or}&\\
h_g=0&&h_g=0&&h_g=0\\
\end{array}
\]

We say the equation (\ref{yamabebdry}) is of positive/negative/zero type if $\lambda_1$ as defined in (\ref{type}) is positive/negative/zero respectively (see \cite{HL1} for more discussion). When $c_2=0$, by the Hopf lemma, the positive/negative/zero type implies that $c_1>0/c_1<0/c_1=0$. In \cite{Esc1}, Escobar proved that the equation (\ref{yamabebdry}) is solvable for some $c_2>0$ and some $c_2<0$ under certain hypothesis. In \cite{HL1}, and \cite{HL2}, Han and Li confirmed the existence of the solutions to the equation (\ref{yamabebdry}) when the manifold is of positive type and is locally conformally flat with umbilic boundary or with non totally umbilic boundary of dimension $n\ge 5$. In this paper, we will study the equation (\ref{yamabebdry}) of negative type. More generally, we will study a fully nonlinear version of the negative type being stated as follows.

Let $Ric_g$ denote the Ricci tensor of $g$. Consider the modified Schouten tensor of $g$ as introduced in \cite{LL2}
\[
A^t_g:=\frac{1}{n-2}\Big (Ric_g-\frac{t R_g}{2(n-1)}g \Big),\qquad t\le 1.
\]

Note that $A^0_g=Ric_g$ and $A^1_g=A_g$ is the Schouten tensor (see \cite{Eis}). Schouten tensor as a $(0,2)$ tensor appears in the decomposition of the Riemann tensor, i.e., the Riemann tensor can be decomposed as the direct sum of the Weyl tensor and the Kulkarni-Numizu product of $A_g$ with $g$. In \cite{LL2}, we introduced $A_g^t$ up to a constant multiple. In fact, we introduced the 
tensor $sA_g+\frac{(1-s)R_g}{2(n-1)}g=sA_g^t$ with $t=n-1-\frac{n-2}{s}$.

Assume that
\begin{equation}
\Gamma\subset\Bbb{R}^n\quad\mbox{is an open convex symmetric cone with vertex at the origin}
\label{cone}
\end{equation}
satisfying
\begin{equation}
\Gamma_n:=\{\lambda=(\lambda_1,\cdots\lambda_n)\in\Bbb{R}^n|~\lambda_1>0,\cdots\lambda_n>0\}\subset\Gamma\subset\Gamma_1:=\{\lambda\in\Bbb{R}^n|~\sum\limits_{i=1}^{n}\lambda_i>0\},
\label{gammacone}
\end{equation}
where $\Gamma$ being symmetric means that
\[
(\lambda_1,\cdots\lambda_n)\in\Gamma\Longleftrightarrow (\lambda_{i_1},\cdots\lambda_{i_n})\in\Gamma
\]
for any permutation $(i_1,\cdots,i_n)$ of $(1,\cdots,n)$.

For $\alpha_0\in (0,1)$, assume that
\begin{equation}
f\in C^{2,\alpha_0}(\Gamma)\cap C^0(\bar{\Gamma})\quad\mbox{is concave, homogeneous of degree 1 and symmetric in }\lambda_i,
\label{f0}
\end{equation}
satisfying
\begin{equation}
f|_{\partial\Gamma}=0,\quad\nabla f\in\Gamma_n\quad\mbox{on}~~\Gamma,
\label{f1}
\end{equation}
\begin{equation}
\lim\limits_{s\to\infty}f(s\lambda)=\infty,\quad\lambda\in\Gamma,
\label{f2}
\end{equation}
and

\begin{equation}
f(\lambda)\le\frac{1}{\bar{\epsilon}}\sum\limits_{i=1}^n\lambda_i,\quad\sum\limits_{i=1}^n f_{\lambda_i}(\lambda)\ge\bar{\epsilon}\quad\mbox{on the level set}~ \{f=1\}
\label{f3}
\end{equation}
for some constant $\bar{\epsilon}>0$.

Notice that $f$ is homogeneous of degree $1$. Therefore $f_{\lambda_i}$ is homogeneous of degree $0$ and the above assumption (\ref{f3}) also holds in $\Gamma$.

Let $\lambda_g(A^t_g)$ denote the eigenvalues of $A_g^t$ w.r.t. the metric $g$. A fully nonlinear problem of negative admissible curvature is to look for a metric $\tilde{g}\in [g]$ solving
\begin{equation}
\{
\begin{array}{lcl}
f(-\lambda_{\tilde{g}}(A^t_{\tilde{g}}))&=&1,\quad -\lambda_{\tilde{g}}(A^t_{\tilde{g}})\in\Gamma\quad\mbox{on}~~M\\
h_{\tilde{g}}&=&c\quad\mbox{on}~~\partial M,
\end{array}
\label{yamabegeneral}
\end{equation}
if $-\lambda_g(A^t_g)\in\Gamma$ on $M$ and $h_g\le 0$ on $\partial M$, where $c$ is a constant.

When $(f,\Gamma)=(\sigma_k^{\frac 1k},\Gamma_k)$, the problem is the $k-$th Yamabe problem of negative admissible type, where
\[
\sigma_k(\lambda)=\sum\limits_{1\le i_1<\cdots<i_k\le n}\lambda_{i_1}\cdots\lambda_{i_k},\quad
\Gamma_k:=\{\lambda\in\Bbb{R}^n|\sigma_1(\lambda)>0,\cdots,\sigma_k(\lambda)>0\}.
\]

It is well-known that $(\sigma_k^{\frac 1 k},\Gamma_k)$ satisfies assumptions (\ref{cone})-(\ref{f3}). In particular, when $k=1$, the problem (\ref{yamabegeneral}) is equivalent to solving the equation (\ref{yamabebdry}) of negative type. This is
because $\sigma_1(-\lambda_{\tilde{g}}(A^t_{\tilde{g}})))=-\frac{1}{n-2}(1-\frac{nt}{2(n-1)})R_{\tilde{g}}$, and the assumption $-\lambda_{g}(A^t_g)\in\Gamma_1$, $h_g\le 0$ is to say that $R_g<0$ and $h_g\le 0$, which implies that $\lambda_1<0$ by taking $\phi\equiv 1$ in (\ref{type}). Conversely, if the equation is of negative type, we can assume $R_g<0$ and $h_g=0$. Clearly the solution $u$ of the equation (\ref{yamabebdry}) also gives a solution $\tilde{g}=u^{\frac{4}{n-2}}g$ to the problem (\ref{yamabegeneral}).

In \cite{GV}, Gursky and Viaclovsky proved that, for $t<1$, there exists a unique solution $\tilde{g}\in [g]$ solving
\[
\sigma_k(-\lambda_{\tilde{g}}(A^t_{\tilde{g}}))=1,\quad -\lambda_{\tilde{g}}(A^t_{\tilde{g}})\in\Gamma_k
\]
if the compact manifold of dimension $n\ge 3$ has no boundary and $-\lambda_g(A^t_g)\in\Gamma_k$.
\begin{theorem}
Let $(M^n,g)$ be an $n-$dimensional ($n\ge 3$) compact smooth Riemannian manifold with
$\partial M\neq\emptyset$, and let $(f,\Gamma)$ be a pair satisfying (\ref{cone})-(\ref{f3}). Assume that $-\lambda_g(A_g^t)\in\Gamma$ on $M$ and $h_g\le 0$ on $\partial M$. Then, for $c\le 0$ and for $t<1$, there exists a unique solution $\tilde {g}=e^{2v}g$ solving the problem (\ref{yamabegeneral}). Moreover,
\[
\|v\|_{C^{4,\alpha_0}(M^n,g)}\le C,
\]
where $C>0$ is a universal constant depending only on $(M^n,g)$, $(f,\Gamma)$, $\alpha$, $t$, and $|c|$.
\label{constantresult}
\end{theorem}
The next theorem is a more general result.
\begin{theorem}
Let $(M^n,g)$ be an $n-$dimensional ($n\ge 3$) compact smooth Riemannian manifold with
$\partial M\neq\emptyset$, and let $(f,\Gamma)$ be a pair satisfying (\ref{cone})-(\ref{f3}). Assume
that $-\lambda_g(A^t_g)\in\Gamma$ on $M$ and $h_g\le 0$ on $\partial M$. Given any $0<\phi\in C^{2,\alpha_0}(M^n)$ and any $0\ge \psi\in C^{3,\alpha_0}(\partial M)$, then, for $t<1$, there exists a unique solution $\tilde {g}=e^{2v}g$ solving
\begin{equation}
\{
\begin{array}{lcl}
f(-\lambda_{\tilde{g}}(A^t_{\tilde{g}}))&=&\phi,\quad -\lambda_{\tilde{g}}(A^t_{\tilde{g}})\in\Gamma\quad\mbox{on}~~M\\
h_{\tilde{g}}&=&\psi\quad\mbox{on}~~\partial M.
\end{array}
\label{generalproblem}
\end{equation}
Moreover
\[
\|v\|_{C^{4,\alpha_0}(M^n,g)}\le C,
\]
where $C>0$ is a universal constant depending only on $(M^n,g)$, $(f,\Gamma)$, $\phi$, $\psi$, $\alpha$, and $t$.
\label{generalresult}
\end{theorem}

In the above theorems, we do not assume the boundary $\partial M$ is umbilic or the manifold is locally conformally flat near $\partial M$, so when we establish the a-priori estimates on the boundary, we can not assume $\partial M$ is totally geodesic, which offers a very useful geodesic normal coordinates, i.e., locally, one direction of the geodesic normal coordinates is the normal direction and all the other directions of coordinates are the tangent directions of $\partial M$. On the general manifolds, the lack of such coordinates causes the a-priori estimates much more difficult to
obtain. The Yamabe problem of the negative case can avoid this particular assumption on the boundary of the manifold since the problem is variational and the minimizing sequence is convergent. However, our problem (\ref{generalproblem}) may not even be variational. To overcome this difficulty, we introduce a very useful coordinates near $\partial M$ in Section 4, called the tubular neighborhood normal coordinates. Such coordinates allow us get rid of the assumption of the umbilic boundary, which is very important in the following theorem. As an application of the above theorems, we affirm the existence of certain Riemannian metrics on a general compact smooth differential manifold with some boundary.
\begin{theorem}
Let $(f,\Gamma)$, $\phi,\psi$ be as in Theorem~\ref{generalresult}. Any compact $n-$dimensional ($n\ge 3$) smooth differential manifold with some boundary always admits a
smooth Riemannian metric $\tilde {g}$ with the negative Ricci tensor satisfying
\[
\{
\begin{array}{lcl}
&&\det(-Ric_{\tilde{g}})=1 \quad\mbox{on}~~M,\\
&&h_{\tilde{g}}=0\quad\mbox{on}~~\partial M.\\
\end{array}
\]
More generally, for $t<1$, any compact $n-$dimensional ($n\ge 3$) smooth differential manifold with some boundary always admits a $C^{4,\alpha_0}$ Riemannian metric $\tilde {g}$ satisfying

\[
\{
\begin{array}{lcl}
f(-\lambda_{\tilde{g}}(A^t_{\tilde{g}}))&=&\phi,\quad -\lambda_{\tilde{g}}(A^t_{\tilde{g}})\in\Gamma\quad\mbox{on}~~M\\
h_{\tilde{g}}&=&\psi\quad\mbox{on}~~\partial M.
\end{array}
\]
\label{existence}
\end{theorem}

We want to point out that a similar problem of positive admissible curvature has been studied by quite a few people and many important results have been obtained such as \cite{CGY1}, \cite{GV1}, \cite{GW}, \cite{JLL}, \cite{LL1}, \cite{LL2}, \cite{STW} and the references therein. If we write the equation (\ref{existence}) in $v$ with $\tilde{g}=e^{2v}g$. Then the equation becomes a fully nonlinear elliptic equation in $v$ with the exact form being given in section 2. In general, fully nonlinear elliptic equations involving $f(\lambda(D^2 v))$ have been studied by Caffarelli, Nirenberg and Spruck (\cite{CNS}) and many others. Fully nonlinear equations involving $f(\lambda(\nabla_g^2 v+g))$ have been investigated by Li (\cite{L}), Urbas (\cite{U}) and others.

We organize our paper as follows. In section 2, we present some prerequisites and prove the uniqueness of the solution. We establish the $C^0$ estimates in section 3. In section 4, we introduce the tubular neighborhood normal coordinates and discuss some of its properties. In the next two sections, we use such coordinates to derive the gradient estimates and the Hessian estimates. In section 7, we establish the existence of the solution to the equation (\ref{generalproblem}). In the last section, we prove the Theorem~\ref{existence}.

{\bf Acknowledgment}: The first author would like to express her appreciation to Professor Yanyan Li for his valuable suggestions on possible topics to work on.
\begin{section}
{Uniqueness}
\end{section}

For $\tilde{g}\in [g]$, write $\tilde{g}=e^{2v}g$. We have the conformal transformation
\[
\{
\begin{array}{l cl}
A_{\tilde{g}}^t&=&-W_g^v+e^{2 v}A_g^t\\
h_{\tilde{g}}&=& (h_g+v_{\nu})e^{-v},\\
\end{array}
\]
where $\frac{\partial}{\partial\nu}$ is the unit outer normal of $g$ on $\partial M$ and
\[
W_g^v:=\nabla_g^2 v+\frac{1-t}{n-2}(\Delta_g v)g +\frac{2-t}{2}|\nabla v|_g^2 g-d v\otimes d v.
\]
The equation (\ref{generalproblem}) is equivalent to solving
\begin{equation}
\{
\begin{array}{l c l}
f(\lambda_g(W_g^v- A_g^t))&=&\phi (x) e^{2 v},\quad \lambda_g(W_g^v- A_g^t)\in\Gamma \quad\mbox{on}~~M\\
h_g+v_{\nu}&=&e^{v}\psi (x)\quad\mbox{on}~~\partial M.\\
\end{array}
\label{realgeneral}
\end{equation}
{\bf Proof of the Uniqueness.} In this section, we give an independent proof of the uniqueness of the solution for $t\le 1$ even though we can see this later from the method of continuity and a suitable homotopy for $t<1$. Let $v_1,v_2$ be two solutions of the equation (\ref{realgeneral}), and let $g_i=e^{2v_i}g$ for $i=1,2$. Write $g_2=e^{2w}g_1$ with $w=v_2-v_1$. Then $v_2$ is a solution of the equation (\ref{realgeneral}) is to say
\[
\{
\begin{array}{l c l}
f(\lambda_{g_1}(W_{g_1}^w- A_{g_1}^t))&=&\phi (x) e^{2 w},\quad \lambda_{g_1}(W_{g_1}^w- A_{g_1}^t)\in\Gamma\quad\mbox{on}~~M\\
h_{g_1}+w_{\nu_1}&=&e^w\psi (x)\quad\mbox{on}~~\partial M,\\
\end{array}
\]
where $\frac{\partial}{\partial\nu_1}$ is the unit outer normal w.r.t. $g_1$ on $\partial M$.

Note that $v_1$ is also a solution of (\ref{realgeneral}), so $h_{g_1}=\psi$ and the above equation becomes
\begin{equation}
\{
\begin{array}{l c l}
f(\lambda_{g_1}(W_{g_1}^w- A_{g_1}^t))&=&\phi (x) e^{2 w},\quad \lambda_{g_1}(W_{g_1}^w- A_{g_1}^t)\in\Gamma\quad\mbox{on}~~M\\
w_{\nu_1}&=&(e^w-1)\psi (x)\quad\mbox{on}~~\partial M.\\
\end{array}
\label{uniquew2}
\end{equation}

Let $w(x_0)=\max\limits_{M}w$.
\begin{lemma}
$w(x_0)\le 0$.
\label{uniquel1}
\end{lemma}
{\bf Proof of the Lemma~\ref{uniquel1}}

{\bf Case 1}. If $x_0$ is an interior point of $M$, then $\nabla_{g_1}w(x_0)=0$, $\nabla_{g_1}^2 w(x_0)\le 0$, and
\[
W_{g_1}^w(x_0)=\nabla_{g_1}^2 w(x_0)+\frac{1-t}{n-2}(\Delta_{g_1}w)(x_0)g_1(x_0)\le 0,
\]
which, together with (\ref{f1}), implies that
\[
\phi(x_0)e^{2 w(x_0)}= f(\lambda_{g_1}(W_{g_1}^w-A_{g_1}^t)(x_0))\le f(\lambda_{g_1}(-A_{g_1}^t)(x_0))=\phi (x_0),
\]
therefore $e^{2w(x_0)}\le 1$, i.e., $w(x_0)\le 0$.

{\bf Case 2}. If $x_0\in \partial M$, then $w_{\nu_1}(x_0)\ge 0$. By the second equation in (\ref{uniquew2}), we know that
\[
0\le w_{\nu_1}(x_0)= (e^w-1)(x_0)\psi (x_0),
\]
so either $(e^w-1)(x_0)\le 0$ when $\psi (x_0)<0$, or $w_{\nu_1}(x_0)=0$ when $\psi (x_0)=0$, that is, when $\psi (x_0)<0$, we have $w(x_0)\le 0$, and when $\psi (x_0)=0$, $w_{\nu_1}(x_0)=0$ gives $\nabla_{g_1}w(x_0)=0$, therefore $\nabla_{g_1}^2w(x_0)\le 0$. We can proceed as in case 1 to conclude that $w(x_0)\le 0$. Lemma~\ref{uniquel1} has been established. $\clubsuit$
\vskip 5pt

Let $w(y_0)=\min\limits_M w$.
\begin{lemma}
$w(y_0)\ge 0$.
\label{uniquel2}
\end{lemma}
{\bf Proof of the Lemma~\ref{uniquel2}}

{\bf Case 1}. If $y_0$ is an interior point of $M$, then $\nabla_{g_1}w(y_0)=0$, $\nabla_{g_1}^2 w(y_0)\ge 0$, and
\[
W_{g_1}^w(y_0)=\nabla_{g_1}^2 w(y_0)+\frac{1-t}{n-2}(\Delta_{g_1}w)(y_0)g_1(y_0)\ge 0,
\]
which implies that
\[
\begin{array}{l cl}
\phi(y_0)e^{2 w(y_0)}&=& f(\lambda_{g_1}(W_{g_1}^w-A_{g_1}^t)(y_0))\ge f(\lambda_{g_1}(-A_{g_1}^t)(y_0))=\phi (y_0),
\end{array}
\]
therefore $e^{2w(y_0)}\ge 1$, i.e., $w(y_0)\ge 0$.

{\bf Case 2}. If $y_0\in \partial M$, then $w_{\nu_1}(y_0)\le 0$ and
\[
0\ge w_{\nu_1}(y_0)= (e^w-1)(y_0)\psi (y_0),
\]
so either $\psi (y_0)<0$, we have $w(y_0)\ge 0$, or $\psi (y_0)=0$,
then $w_{\nu_1}(y_0)=0$, which implies that $\nabla_{g_1}w(y_0)=0$, therefore $\nabla_{g_1}^2w(y_0)\ge 0$. We can proceed as in case 1 to conclude that $w(y_0)\ge 0$. Lemma~\ref{uniquel2} has been established. $\clubsuit$
\vskip 5 pt
Combining Lemma~\ref{uniquel1} and Lemma~\ref{uniquel2}, we have $w\equiv 0$, that is, $v_1\equiv v_2$. The uniqueness of the solution of the equation (\ref{generalproblem}) has been proved. $\clubsuit$
\begin{remark}
When $k=1$, $c=0$, the uniqueness of the solution has been obtained by Cheerier in \cite{Cher} and it implies that the solution must be the unique minimum point of ${\cal F}$.
\end{remark}
\begin{section}
{$C^0$ estimates}
\end{section}
When the manifold has some boundary, the $C^0$ estimate is not a trivial consequence of the maximum
principle anymore. In this section, we obtain $C^0$ estimates by establishing the upper bounds and the lower bounds individually.
\begin{lemma}
Let $(M^n,g)$ and $(f,\Gamma)$ be as in Theorem~\ref{generalresult}. For $t\le 1$, let $v$ be a $C^2$ solution of the equation (\ref{generalproblem}). Then there exists a universal constant $C>0$ depending only on $(M^n, g, t)$, $(f,\Gamma)$, $\phi$ and $\psi$ such that
\[
v\le C.
\]
\label{c0l1}
\end{lemma}
{\bf Proof of the Lemma~\ref{c0l1}.} In this paper, if not specified, we will use $C>0$ to denote a universal constant with the dependence as being described in the statement of the Lemma~\ref{c0l1}, but may change from line to line. Since $\lambda_g({A_g^t})\in\Gamma\subset\Gamma_1$ and $h_g\le 0$, we have $R_g<0$ and $h_g\le 0$, from which, we know that $(M^n,g)$ is of negative type. Hence we can find $g_0=e^{2v_0}g$ such that

\[
\{
\begin{array}{l cl}
R_{g_0}<0\quad\mbox{on}~~M&&\\
h_{g_0}=0\quad\mbox{on}~~\partial M.&&\\
\end{array}
\]

Write $e^{2v}g= e^{2\tilde {v}}g_0$ with $\tilde{v}=v-v_0$. Then $\tilde {v}$ satisfies
\begin{equation}
\{
\begin{array}{l c l}
f(\lambda_{g_0}(W_{g_0}^{\tilde{v}}- A_{g_0}^t))&=&\phi (x) e^{2 \tilde{v}},\quad \lambda_{g_0}(W_{g_0}^{\tilde{v}}- A_{g_0}^t)\in\Gamma\quad\mbox{on}~~M\\
{\tilde{v}}_{\nu_0}&=&e^{\tilde{v}}\psi (x)\quad\mbox{on}~~\partial M,\\
\end{array}
\label{c0e1}
\end{equation}
where $\frac{\partial}{\partial\nu_0}$ is the unit outer normal of $g_0$ on $\partial M$.

Let $\tilde {v}(x_0)=\max\limits_M(\tilde {v})$.

{\bf Case 1.} If $x_0$ is an interior point of $M$, then $\nabla_{g_0}\tilde{v}(x_0)=0$, $\nabla^2_{g_0}\tilde {v}(x_0)\le 0$ and therefore $W_{g_0}^{\tilde{v}}(x_0)\le 0$. Hence
\[
\lambda_{g_0}(W_{g_0}^{\tilde{v}}- A_{g_0}^t)(x_0)\in\Gamma
\]
implies that $\lambda_{g_0}(- A_{g_0}^t)\in\Gamma$. Thus by (\ref{f2}) and (\ref{f3}),
\[
\begin{array}{l cl}
e^{2\tilde{v}(x_0)}\phi (x_0)&=& f(\lambda_{g_0}(W_{g_0}^{\tilde{v}}- A_{g_0}^t)(x_0))\le
f(\lambda_{g_0}(- A_{g_0}^t)(x_0))\\
&\le& C\sigma_1(\lambda_{g_0}(- A_{g_0}^t)(x_0))\le C\max\limits_M(-R_g)=:C,
\end{array}
\]
so we have $\tilde{v}(x_0)\le C$.

{\bf Case 2.} If $x_0\in\partial M$, then $\psi (x_0)=0$. If not, then at $x_0$, the second equation in (\ref{c0e1}) implies that
\[
0\le \tilde{v}_{\nu_0}(x_0)=e^{\tilde{v}}(x_0)\psi(x_0)<0,
\]
which is a contradiction. Thus $\tilde{v}_{\nu_0}(x_0)=e^{\tilde{v}}(x_0)\psi(x_0)=0$, $\nabla_{g_0}\tilde{v}(x_0)=0$, and $\nabla^2_{g_0}\tilde {v}(x_0)\le 0$. We can proceed as in case 1 to obtain $v(x_0)\le C$.

Combining the above two cases, we have $\tilde {v}\le C$, which means $v\le C$. Lemma~\ref{c0l1} has been established. $\clubsuit$
\begin{lemma}
Let $(M^n,g)$ and $(f,\Gamma)$ be as in Theorem~\ref{generalresult}. For $t\le 1$, let $v$ be a $C^2$ solution of the equation (\ref{generalproblem}). Then there exists a universal constant $C>0$ depending only on $(M^n, g, t)$, $(f,\Gamma)$, $\phi$ and $\psi$ such that
\[
v\ge -C.
\]
\label{c0l2}
\end{lemma}
{\bf Proof of the Lemma~\ref{c0l2}.} Let $\bar {w}$ be a smooth function such that $\bar{w}$ is the distance function to $\partial M$ near the boundary and $\bar{w}$ takes value in $[0,1]$ in general. Then $\bar{w}_{\nu}|_{\partial M}\equiv -1$. Let $g_0=e^{2\epsilon_0 \bar{w}}g$ with $\epsilon_0>0$ being a constant to be chosen later. We have

\begin{equation}
h_{g_0}=(h_g+\epsilon_0\bar{w}_{\nu})e^{-\epsilon_0\bar{w}}\le -\epsilon_0 e^{-\epsilon_0} <0,
\label{c0l2e0}
\end{equation}
and
\[
\begin{array}{l c l}
-\lambda_{g_0}(A_{g_0}^t)&=&\lambda_{g}\Big( \epsilon_0[\nabla_g^2 \bar{w}+\frac{1-t}{n-2}(\Delta_g \bar{w})g \\
&&+\frac{2-t}{2}\epsilon_0|\nabla \bar{w}|_g^2 g-\epsilon_0 d \bar{w}\otimes d \bar{w}]-A_g^t\Big),
\end{array}
\]
so we can take $\epsilon_0\ll 1$ depending only on $(M^n, g,t,f,\Gamma)$ such that
\begin{equation}
-\lambda_{g_0}(A_{g_0}^t)\in\Gamma\quad\mbox{and}\quad f(-\lambda_{g_0}(A_{g_0}^t))\ge \frac 12\min\limits_M f(-\lambda_g(A_g^t)).
\label{c0l2e1}
\end{equation}
Let $\tilde {v}=v-\epsilon_0\bar{w}$. Then $e^{2v}g=e^{2\tilde{v}}g_0$ and $\tilde{v}$ solves
 \begin{equation}
\{
\begin{array}{l c l}
f(\lambda_{g_0}(W_{g_0}^{\tilde{v}}- A_{g_0}^t))&=&\phi (x) e^{2 \tilde{v}},\quad \lambda_{g_0}(W_{g_0}^{\tilde{v}}- A_{g_0}^t)\in\Gamma\quad\mbox{on}~~M\\
{\tilde{v}}_{\nu_0}+h_{g_0}&=&e^{\tilde{v}}\psi (x)\quad\mbox{on}~~\partial M,\\
\end{array}
\label{c0l2e2}
\end{equation}

Let $\tilde {v}(y_0)=\min\limits_{M}\tilde {v}$.

{\bf Case 1.} If $y_0$ is in the interior of $M$, then $\nabla_{g_0}\tilde{v}(y_0)=0$, $\nabla^2_{g_0}\tilde {v}(y_0)\ge 0$ and $W_{g_0}^{\tilde{v}}(y_0)\ge 0$. Hence by (\ref{f2}), (\ref{c0l2e1}) and (\ref{c0l2e2}),
\[
\begin{array}{l cl}
e^{2\tilde{v}(y_0)}\phi (y_0)&=& f(\lambda_{g_0}(W_{g_0}^{\tilde{v}}- A_{g_0}^t)(y_0))\ge
f(\lambda_{g_0}(- A_{g_0}^t)(y_0))\\
&\ge& \frac {1}{2}\min\limits_{M} f(\lambda_g(- A_g^t)),\\
\end{array}
\]
i.e.,
\[
\tilde {v}(y_0)\ge \frac 12\min\limits_M(\ln (\frac {1}{2\phi}\min\limits_M f(\lambda_g(- A_g^t))))\ge -C.
\]

{\bf Case 2.} If $y_0\in\partial M$, then $\tilde{v}_{\nu_0}(y_0)\le 0$. By (\ref{c0l2e0}) and (\ref{c0l2e2}),
\[
-\epsilon_0 e^{-\epsilon_0}\ge h_{g_0}(y_0)+\tilde{v}_{\nu_0}(y_0)=e^{\tilde{v}(y_0)}\psi (y_0)\ge -C e^{\tilde{v}(y_0)},
\]
so
\[
\tilde{v}(y_0)\ge \ln \frac{\epsilon_0}{C}-\epsilon_0\ge -C.
\]
Combining the above two cases, we know $\tilde {v}\ge -C$, hence $v\ge -C$. Lemma~\ref{c0l2} has been proved. $\clubsuit$

\begin{section}
{Tubular Neighborhood Normal Coordinates}
\end{section}

The main issue of the gradient and the Hessian estimates is the bounds on the boundary of $M$. For this reason, we need to introduce certain coordinates near $\partial M$. Let $g|_{\partial M}$ be the induced metric of $g$ on $\partial M$, and let $\delta_1>0$ be the minimum of the injectivity radius of $(M^n,g)$ and the injectivity radius of $(\partial M,g|_{\partial M})$. Consider the map $E:\partial M\times [0,\delta_1)\to M$ by $E(y,t)=\exp_y(-t\frac{\partial}{\partial\nu})$. Since $E(y,0)=y$ implies that, for any $y\in\partial M$, $dE|_{(y,0)}(X)=X$ for $X\in T_y(\partial M)$, and $dE|_{(y,0)}(\frac{d}{dt})=-\frac{\partial}{\partial\nu}\neq 0$. That is, $dE|_{(y,0)}$ is an isomorphism from $T_{(y,0)}(\partial M\times [0,\delta_1))\to T_y M$. By the Implicit Function Theorem, there exists some constant $\delta_y\in (0,\delta_1)$ such that $E$ is a smooth diffeomorphism on $(\partial M\cap B_{\delta_y}(y))\times [0,\delta_y)$, where $B_{\delta_y}(y)$ is the open geodesic ball of $(M^n,g)$ centered at $y$ with radius $\delta_y$. By shrinking $B_{\delta_y}(y)$, we can also assume the exponential map of $(\partial M,g|_{\partial M})$ at $y$ is a smooth diffeomorphism in $B_{\delta_y}(y)\cap\partial M$. Now we extend $\frac{\partial}{\partial\nu}$ to the interior of $M$, still denoted by $\frac{\partial}{\partial\nu}$ such that $\frac{\partial}{\partial\nu}|_{E(z,t)}=-\frac{dE}{dt}|_{(z,t)}$ for any $z\in \partial M\cap B_{\delta_y}(y)$. Then $\frac{\partial}{\partial\nu}$ is a smooth unit vector field in $E\Big((\partial M\cap B_{\delta_y}(y))\times [0,\delta_y)\Big)$.

\begin{prop}
For any $y_0\in\partial M,
$\[
B_{\frac{\delta_{y_0}}{2}}(y_0)\subset E\Big((\partial M\cap B_{\delta_{y_0}}(y_0))\times [0,\delta_{y_0})\Big),
\]
and for any $y\in B_{\frac{\delta_{y_0}}{2}}(y_0)$, there exists a unique $\bar{y}\in\partial M$ such that $d(y,\bar{y})=d(y,\partial M)$. Moreover $\bar{y}\in B_{\delta_{y_0}}(y_0)\cap\partial M$.
\label{c2p0}
\end{prop}
{\bf Proof of the Proposition~\ref{c2p0}} For any $y\in B_{\frac{\delta_{y_0}}{2}}(y_0)$,
\[
s:=d(y,\partial M)\le d(y,y_0)<\frac{\delta_{y_0}}{2}.
\]

For any $z\in\partial M\setminus B_{\delta_{y_0}}(y_0)$,
\[
d(y,z)\ge d(z,y_0)-d(y,y_0)>\delta_{y_0}-\frac{\delta_{y_0}}{2}=\frac{\delta_{y_0}}{2}.
\]

Thus if $d(y,\partial M)=d(y,\bar{y})$ for some $\bar{y}\in\partial M$, then $\bar{y}\in \partial M\cap B_{\delta_{y_0}}(y_0)$. Let $r(t)$ be the normalized geodesic connecting $y$ and $\bar{y}$ such that $r(0)=\bar{y}$ and $r(s)=y$. Then $\frac{dr}{dt}|_{t=0}=-\frac{\partial}{\partial\nu}|_{\bar{y}}$, that is $y=E(\bar{y},s)$. Therefore $y\in E\Big((\partial M\cap B_{\delta_{y_0}}(y_0))\times [0,\delta_{y_0})\Big)$,
and $B_{\frac{\delta_{y_0}}{2}}(y_0)\subset E\Big((\partial M\cap B_{\delta_{y_0}}(y_0))\times [0,\delta_{y_0})\Big)$. Recall $E$ is a smooth diffeomorphism in $(\partial M\cap B_{\delta_{y_0}}(y_0))\times [0,\delta_{y_0})$ and $\bar{y}\in \partial M\cap B_{\delta_{y_0}}(y_0)$. Thus $\bar{y}$ is uniquely determined by $y$. The Proposition~\ref{c2p0} has been proved. $\clubsuit$
\vskip 5pt

By the Proposition~\ref{c2p0}, $\frac{\partial}{\partial\nu}=-\frac{dE}{dt}$ is a smooth unit vector field in $B_{\frac{\delta_{y_0}}{2}}(y_0)$. Moreover, in $B_{\frac{\delta_{y_0}}{2}}(y_0)$, the parameter $t$ in $E(y,t)$ is the distance parameter to the boundary of $M$, which can be derived more precisely as in establishing (\ref{c2p1e1}). Let $\{y_j\}_{j=1}^{n-1}$ be the geodesic normal coordinates w.r.t. the metric $g|_{\partial M}$ at $y_0$. Then $\{y_j\}_{j=1}^{n-1}$ is smooth and well-defined in
$\partial M\cap B_{\delta_{y_0}}(y_0)$. For any $y\in B_{\frac{\delta_{y_0}}{2}}(y_0)$, there is a unique $\bar{y}\in\partial M$ such that $d(y,\bar{y})=d(y,\partial M)$. By the Proposition~\ref{c2p0}, $\bar{y}\in \partial M\cap B_{\delta_{y_0}}(y_0)$. Let $(y_1,\cdots,y_{n-1})$ be the geodesic normal coordinates of $\bar{y}$ w.r.t. the metric $g|_{\partial M}$ at $y_0$. Define $(y_1,\cdots,y_{n-1},y_n)$ as the coordinates of $y$ with $y_n=d(y,\partial M)$. Such coordinates are well-defined and smooth in $B_{\frac{\delta_{y_0}}{2}}(y_0)$. The reason is that $\bar{y}$ is uniquely determined and $\bar{y}\in \partial M\cap B_{\delta_{y_0}}(y_0)$, which implies that $y=E(\bar{y},y_n)$. Hence the map from $y$ to $(\bar{y},y_n)$ is the inverse of the smooth diffeomorphism $E$, therefore is also a smooth diffeomorphism, that is to say $(y_1,\cdots,y_n)$ is well-defined and smooth in $B_{\frac{\delta_{y_0}}{2}}(y_0)$. We call such coordinates the tubular neighborhood normal coordinates of $y$ at $y_0$. Observe that $g(\frac{\partial}{\partial y_i},\frac{\partial}{\partial y_j})=\delta_{ij}$ for $1\le i,j\le n-1$ at $y_0$. Moreover, such coordinates has the following proposition.

\begin{prop}
For $1\le j\le n-1$,
\[
\frac{\partial}{\partial y_n}=-\frac{\partial}{\partial\nu},\quad g(\frac{\partial}{\partial y_j},\frac{\partial}{\partial y_n})=0,\quad\mbox{in}~B_{\frac{\delta_{y_0}}{8}}(y_0).
\]
\label{c2p1}
\end{prop}
{\bf Proof of the Proposition~\ref{c2p1}} For any $y\in B_{\frac{\delta_{y_0}}{8}}(y_0)$ with $(a_1,\cdots,a_n)$ as its tubular neighborhood normal coordinates at $y_0$. Let $\bar{y}\in\partial M\cap B_{\delta_{y_0}}(y_0)$ be the unique point such that $d(y,\bar{y})=a_n<\frac{\delta_{y_0}}{8}$. Clearly $\bar{y}\in B_{\frac{\delta_{y_0}}{4}}(y_0)$ since
\[
d(y,z)\ge d(z,y_0)-d(y,y_0)>\frac{\delta_{y_0}}{4}-\frac{\delta_{y_0}}{8}=\frac{\delta_{y_0}}{8}\quad\mbox{for any}~ z\in\partial M\setminus B_{\frac{\delta_{y_0}}{4}}(y_0).
\]

Let $r(t)=E(\bar{y},t)$. Then $r$ is smooth and well-defined for $t\in [0,\delta_{y_0})$. For $t\in [0,\frac{\delta_{y_0}}{8})$, by
\[
d(r(t),y_0)\le d(r(t),\bar{y})+d(\bar{y},y_0)<\frac{\delta_{y_0}}{8}+\frac{\delta_{y_0}}{4}<\frac{\delta_{y_0}}{2},
\]
there exists a unique $\tilde{y}\in\partial M$ such that
\[
d(r(t),\tilde{y})=d(r(t),\partial M)=:d^t\le d(r(t),\bar{y})\le t.
\]

By the Proposition~\ref{c2p0}, $\tilde{y}\in B_{\delta_{y_0}}(y_0)\cap\partial M$ and $E(\tilde{y},d^t)=r(t)=E(\bar{y},t)$. Therefore $\tilde{y}=\bar{y}$ and $d^t=t$ since $E$ is a smooth diffeomorphism on $(\partial M\cap B_{\delta_{y_0}}(y_0))\times [0,\delta_{y_0})$. From which, we know
that $(a_1,\cdots,a_{n-1},t)$ is the tubular neighborhood normal coordinates of $r(t)$ at $y_0$ for $t\in [0,\frac{\delta_{y_0}}{8})$. Hence, for $t\in [0,\frac{\delta_{y_0}}{8})$,
\[
\frac{\partial}{\partial y_n}|_{r(t)}=\frac{dr}{dt}|_t=\frac{dE}{dt}|_{(\bar{y},t)}=-\frac{\partial}{\partial\nu}|_{E(\bar{y},t)}=-\frac{\partial}{\partial\nu}|_{r(t)}.
\]

In particular,

\[
\frac{\partial}{\partial y_n}|_y=\frac{\partial}{\partial y_n}|_{r(a_n)}=-\frac{\partial}{\partial\nu}|_{r(a_n)}=-\frac{\partial}{\partial\nu}|_y.
\]

To prove the second statement in the proposition, we consider the set
\[
{\cal S}:=\{z\in B_{\frac{\delta_{y_0}}{8}}(y_0)|~d(z,\partial M)=a_n \}.
\]

Clearly, $y\in{\cal S}\neq\emptyset$. For any $z\in {\cal S}$, let $r(t)=E(\bar{z},t)$ for some $\bar{z}\in B_{\frac{\delta_{y_0}}{4}}(y_0)\cap\partial M$ such that $r(0)=\bar{z}$ and $r(a_n)=z$. As derived earlier, $d(r(t),\partial M)=t$ for any $t\in [0,\frac{\delta_{y_0}}{8})$, which implies that $r([0,\frac{\delta_{y_0}}{8}))$ intersects ${\cal S}$ at a single point $z=r(a_n)$. Moreover, we claim that
\begin{equation}
d(r(t),{\cal S})=t-a_n,\quad\forall~t\in[a_n,\frac{\delta_{y_0}}{8}).
\label{c2p1e1}
\end{equation}

Notice that , for $t\in[a_n,\frac{\delta_{y_0}}{8})$, $d(r(t),{\cal S})\le d(r(t),r(a_n))\le t-a_n$. If (\ref{c2p1e1}) does not hold, then $d(r(t),{\cal S})<t-a_n$, which implies that there exists some $\tilde{z}\in{\cal S}$ such that $d(r(t),\tilde{z})<t-a_n$. Therefore
\[
t=d(r(t),\partial M)\le d(r(t),\tilde{z})+d(\tilde{z},\partial M)<t-a_n+a_n=t,
\]
which is a contradiction. Next, we claim
\begin{equation}
d(r(t),{\cal S})=a_n-t,\quad\forall~t\in[0,a_n).
\label{c2p1e2}
\end{equation}

If not, then $d(r(t),{\cal S})<a_n-t$ since $d(r(t),{\cal S})\le d(r(t),r(a_n))\le a_n-t$, so there exists some $\hat{z}\in {\cal S}$ such that $d(r(t),\hat{z})<a_n-t$, which implies that
\[
a_n=d(\hat{z},\partial M)\le d(r(t),\hat{z})+d(r(t),\partial M)<a_n-t+t=a_n,
\]
which is a contradiction.

By (\ref{c2p1e1}) and (\ref{c2p1e2}), we know $r(a_n)$ is a point in ${\cal S}$ such that $d(r(a_n),r(t))=d(r(t),{\cal S})$ for $t\in [0,\frac{\delta_{y_0}}{8})$, and $r$ is the normalized geodesic connecting $r(t)$ and $r(a_n)$, so $\frac{dr}{dt}|_{a_n}=\frac{dE}{dt}|_{(\bar{z},a_n)}$ is the unit normal vector of ${\cal S}$ at $r(a_n)=z$, i.e., $\frac{\partial}{\partial y_n}=-\frac{\partial}{\partial\nu}=\frac{dE}{dt}$ is the unit normal vector of ${\cal S}$ at $r(a_n)=z$. Let $(b_1,\cdots,b_{n-1},a_n)$ be the tubular neighborhood normal coordinates of $z$ at $y_0$. Observe that, for $1\le k\le n-1$, since $z$ is an interior point of $B_{\frac{\delta_{y_0}}{8}}(y_0)$, the curve
\[
\{(y_1,\cdots,y_k,\cdots,y_n)=(b_1,\cdots,b_{k-1},y_k,b_{k+1},\cdots,a_n)\}\quad\mbox{for $y_k$ near $b_k$}
\]
is contained in ${\cal S}$, which implies that $\{\frac{\partial}{\partial y_k}|_z\}\in T_z{\cal S}$. Hence $g(\frac{\partial}{\partial y_k},\frac{\partial}{\partial y_n})=0$ at $z$ for $1\le k\le n-1$ since $\frac{\partial}{\partial y_n}=-\frac{\partial}{\partial\nu}$ is the normal vector of ${\cal S}$ at $z$.
$z\in {\cal S}$ is arbitrary and $y\in{\cal S}$, so, at $y$, we also have
\[
g(\frac{\partial}{\partial y_k},\frac{\partial}{\partial y_n})=0\quad\mbox{for}~ 1\le k\le n-1.
\]

The Proposition~\ref{c2p1} has been proved. $\clubsuit$
\vskip 5pt

As a simple consequence, we have $g(\frac{\partial}{\partial y_i},\frac{\partial}{\partial y_j})=\delta_{ij}$ at $y_0$ for $1\le i,j\le n$.

 \begin{prop}
 \[
 \{(y_1,\cdots,y_n)|~\sqrt{y_1^2+\cdots+ y_n^2}<\frac{\delta_{y_0}}{16},\quad y_n\ge 0\}\subset B_{\frac{\delta_{y_0}}{8}}(y_0),
 \]
 where $(y_1,\cdots,y_n)$ is the tubular neighborhood normal coordinates at $y_0$.
 \label{c2p2}
 \end{prop}
{\bf Proof of the Proposition~\ref{c2p2}} For any $(y_1,\cdots,y_n)$ with $\sqrt{y_1^2+\cdots+ y_n^2}<\frac{\delta_{y_0}}{16}$ and $y_n\ge 0$, there exists a unique $\bar{y}\in B_{\frac{\delta_{y_0}}{16}}(y_0)$ such that $(y_1,\cdots,y_{n-1})$ is the geodesic normal coordinates of $\bar{y}$ w.r.t. the metric $g|_{\partial M}$ at $y_0$. Consider $r(t)=E(\bar{y},t)$. Then $r(t)$ is smooth for $t\in [0,\frac{\delta_{y_0}}{16})$ and $r([0,\frac{\delta_{y_0}}{16}))\subset B_{\frac{\delta_{y_0}}{8}}(y_0)$. Moreover $d(r(t),\partial M)=t$ for $t\in [0,\frac{\delta_{y_0}}{16})$ as shown earlier. In particular, by $y_n<\frac{\delta_{y_0}}{16}$, $y=E(\bar{y},y_n)$ has $(y_1,\cdots,y_n)$ as its tubular neighborhood normal coordinates at $y_0$. The Proposition~\ref{c2p2} has been proved. $\clubsuit$
\vskip 5pt

Denote $B^T_{\frac{\delta_{y_0}}{16}}(y_0):=\{(y_1,\cdots,y_n)|~\sqrt{y_1^2+\cdots+ y_n^2}<\frac{\delta_{y_0}}{16},\quad y_n\ge 0\}$, which is different from the geodesic ball $B_{\frac{\delta_{y_0}}{16}}(y_0)$. The Proposition~\ref{c2p2} says
that $B^T_{\frac{\delta_{y_0}}{16}}(y_0)\subset B_{\frac{\delta_{y_0}}{8}}(y_0)$. Since $\cup_{y_0\in\partial M}B^T_{\frac{\delta_{y_0}}{64}}(y_0)=\partial M$ and $\partial M$ is compact, we can find $\{y^i\}_{i=1}^N\subset\partial M$ such that $\cup_{i=1}^N\Big( B^T_{\frac{\delta_{y^i}}{64}}(y^i)\cap\partial M\Big)=\partial M$.

 \begin{section}
{Gradient estimates}
\end{section}
\begin{lemma}
Under the same assumptions as in Theorem~\ref{generalresult}, for $t<1$, let $v$ be a $C^3$ solution of the equation (\ref{generalproblem}). Then there exists a universal constant $C>0$ depending only on $(M^n,g,t)$, $(f,\Gamma)$, $\phi$, and $\psi$, such that
\[
|\nabla v|_g\le C\qquad\mbox{on}~~\partial M.
\]
\label{c1l1}
\end{lemma}
{\bf Proof of the Lemma~\ref{c1l1}.} Extend $h_g$ to a smooth function on $M$, and $\psi$ to a $C^{3,\alpha_0}$ function on $M$. More explanation is given in section 7. We still denote the extended functions by $\psi,h_g$ respectively. For each $1\le i_0\le N$, Let $\{y_j\}_{j=1}^n$ be the tubular neighborhood normal coordinates at $y^{i_0}$. Let $\rho=\rho (y_1^2+\cdots+ y_n^2)$ be a smooth cut-off function satisfying
\[
\rho(y)=\Big\{
\begin{array}{l c l}
&1&,\quad\mbox{if}~y\in \overline{B^T_{\frac{\delta_{y^{i_0}}}{32}}(y^{i_0})} \\
&\in [0,1]&,\quad\mbox{if}~y\in B^T_{\frac{\delta_{y^{i_0}}}{16}}(y^{i_0})\setminus B^T_{\frac{\delta_{y^{i_0}}}{32}}(y^{i_0})\\
&0&,\quad\mbox{otherwise},\\
\end{array}
\]
and let $\beta(y)$ be a smooth function in $B^T_{\frac{\delta_{y^{i_0}}}{16}}(y^{i_0})$ satisfying
\[
\beta (y)=\{
\begin{array}{lcl}
y_n, \quad\mbox{if}~~y_n<\delta_0, &&\\
\in[0,2\delta_0],\quad\mbox{o.w.},&&\\
\end{array}
\]
where $0<\delta_0<\frac{\delta_{y^{i_0}}}{32}$ is a very small constant such that $1+2\delta_0\psi e^v>\frac 12$ on $M$ and to be chosen later. Then in $B^T_{\frac{\delta_{y^{i_0}}}{16}}(y^{i_0})\cap\partial M$,
\begin{equation}
\beta\equiv 0,\quad \beta_{\nu}\equiv -1.
\label{c1l1e1}
\end{equation}

Let
\[
\gamma:=(\psi e^v-h_g)\beta,
\]

In the following, we use subindices to denote the covariant derivatives w.r.t. $\frac{\partial}{\partial y_j}$, e.g.,
\[
(v+\gamma)_k=\Big(\nabla (v+\gamma)\Big)(\frac{\partial}{\partial y_k}),\qquad (v+\gamma)_{k\nu}=\Big(\nabla^2 (v+\gamma)\Big) (\frac{\partial}{\partial y_k},\frac{\partial}{\partial\nu}).
\]

Consider
\[
G:=\rho \sum\limits_k(v+\gamma)_k^2\alpha (\frac{v+\gamma+L}{L^2}),
\]
where $L>0$ is a constant satisfying $1<v+\gamma+L<2L$ and $\alpha:~\Bbb{R}^+\to\Bbb{R}^+$ is a smooth positive function to be chosen later.

Notice that $\frac{\partial}{\partial\nu}=-\frac{\partial}{\partial y_n}$ in $B^T_{\frac{\delta_{y^{i_0}}}{16}}(y^{i_0})$ and
\[
B^T_{\frac{\delta_{y^{i_0}}}{16}}(y^{i_0})\cap\partial M=\{y\in B^T_{\frac{\delta_{y^{i_0}}}{16}}(y^{i_0})|~y_n=0~\mbox{and}~\sqrt{y_1^2+\cdots y_{n-1}^2}< \frac{\delta_{y^{i_0}}}{16}\}.
\]
Hence in $B^T_{\frac{\delta_{y^{i_0}}}{16}}(y^{i_0})\cap\partial M$,
\begin{equation}
\rho_{\nu}=-\frac{\partial\rho}{\partial y_n}|_{y_n=0}=0.
\label{c2l1e1}
\end{equation}

\begin{claim}
In $B^T_{\frac{\delta_{y^{i_0}}}{16}}(y^{i_0})\cap\partial M$, $G_{\nu}\equiv 0$.
\label{c1c1}
\end{claim}
{\bf Proof of the Claim~\ref{c1c1}.} In $B^T_{\frac{\delta_{y^{i_0}}}{16}}(y^{i_0})\cap\partial M$, by (\ref{c1l1e1}) and the second equation in (\ref{generalproblem}),
\begin{equation}
\begin{array}{l cl}
(v+\gamma)_{\nu}&=& v_{\nu}+((\psi e^v-h_g)\beta)_{\nu}=(\psi e^v-h_g)+(\psi e^v-h_g)_{\nu}\beta+
(\psi e^v-h_g)\beta_{\nu}\\
&=&(\psi e^v-h_g)-(\psi e^v-h_g)=0\\
\end{array}
\label{c1c1e0}
\end{equation}
Therefore in $B^T_{\frac{\delta_{y^{i_0}}}{16}}(y^{i_0})\cap\partial M$,
\begin{equation}
(v+\gamma)_{k,\nu }=-(v+\gamma)_{k,n }=-(v+\gamma)_{n, k}=(v+\gamma)_{\nu, k}=0,\quad\forall\quad 1\le k\le n-1,
\label{c1c1e00}
\end{equation}
where in the last equality, we used the fact that $\frac{\partial}{\partial y_k}$ is a tangent vector field of $B^T_{\frac{\delta_{y^{i_0}}}{16}}(y^{i_0})\cap\partial M$.

In $B^T_{\frac{\delta_{y^{i_0}}}{16}}(y^{i_0})\cap\partial M$, by (\ref{c2l1e1}) and (\ref{c1c1e0}),
\[
\begin{array}{l cl}
G_{\nu}&=& 2\rho\alpha (\frac{v+\gamma+L}{L^2})\sum\limits_{k=1}^{n}(v+\gamma)_k(v+\gamma)_{k,\nu}\\
&=&2\rho\alpha (\frac{v+\gamma+L}{L^2})(v+\gamma)_n(v+\gamma)_{n,\nu}\quad\mbox{by}~~\ref{c1c1e00})\\
&=&2\rho\alpha (\frac{v+\gamma+L}{L^2})(v+\gamma)_{\nu}(v+\gamma)_{\nu,\nu}=0.\\
\end{array}
\]

Claim~\ref{c1c1} has been proved. $\clubsuit$.
\vskip 5 pt

Let $G(x_0)=\max\limits_{\overline{B^T_{\frac{\delta_{y^{i_0}}}{16}}(y^{i_0})}} G$ for some $x_0\in B^T_{\frac{\delta_{y^{i_0}}}{16}}(y^{i_0})$. W.l.o.g., $G(x_0)\ge 1$. By the Claim~\ref{c1c1}, we have
\[
\nabla G(x_0)=0,\quad\nabla^2 G (x_0)\le 0.
\]

In $B^T_{\frac{\delta_{y^{i_0}}}{16}}(y^{i_0})$

\[
\begin{array}{lcl}
G_i&=&\rho_i\alpha\sum\limits_k(v+\gamma)_k^2+2\rho \alpha (v+\gamma)_{k,i}(v+\gamma)_k+\frac{\rho\alpha'}{L^2}(v+\gamma)_i \sum\limits_k(v+\gamma)_k^2\\
&=&2\rho \alpha (v+\gamma)_{k,i}(v+\gamma)_k+\Big(\frac{\rho_i}{\rho}+\frac{\alpha'}{L^2\alpha}(v+\gamma)_i \Big)G ,\\
\end{array}
\]
so at $x_0$,
\begin{equation}
(v+\gamma)_k(v+\gamma)_{k, i}=-\frac{\alpha'}{2L^2\alpha}\sum\limits_{k}(v+\gamma)^2_k (v+\gamma)_i-\frac{\rho_i}{2\rho}\sum\limits_{k}(v+\gamma)^2_k,
\label{c1e0}
\end{equation}
and
\[
\begin{array}{lcl}
G_{ij}(x_0)&=&2\rho\alpha (v+\gamma)_{k,i}(v+\gamma)_{k,j}+2\rho\alpha (v+\gamma)_k(v+\gamma)_{k,ij}\\
&&+\Big(\frac{\alpha(\rho\rho_{ij}-2\rho_i\rho_j)}{\rho}+\frac{\rho(\alpha\alpha''-2(\alpha')^2)}{L^4\alpha}(v+\gamma)_i(v+\gamma)_j+\frac{\alpha'\rho}{L^2}(v+\gamma)_{ij}\\
&&-\frac{\alpha'}{L^2}(\rho_j(v+\gamma)_i+\rho_i(v+\gamma)_j)\Big)\sum\limits_k (v+\gamma)^2_k\\
\end{array}
\]

Let $F(\bar{A})=f(\lambda (\bar{A}))$ for any symmetric matrix $\bar{A}$ with $\lambda(\bar{A})\in\Gamma$. Let $\{e_i\}_{i=1}^n$ be an orthonormal basis of $T_x M$. Denote $\bar {W}:=
W_g^v-A_g^t$. Let $\bar {W}(e_i,e_j)=\bar{w}_{i j}$ and let $F^{i j}=\frac{\partial F}{\partial\bar {w}_{i j}}$. (\ref{f1}) implies
$(F^{i j})>0$. Denote $\bar{L}^{ij}:=F^{ir}g^{rj}+\frac{1-t}{n-2}(\sum\limits_{i=1}^n F^{ll})g^{ij}$. At $x_0$, assume $e_i=a^j_i\frac{\partial}{\partial y_j}$. Then $g(e_i,e_j)=\delta_{ij}$ is to say that $A^TA={\cal G}^{-1}$, where $A=(a_i^j)$, ${\cal G}^{-1}=(g_{ij})^{-1}$, and $g_{ij}=g(\frac{\partial}{\partial y_i},\frac{\partial}{\partial y_j})$. Denote $B=(F^{ij})$ and $D=(G_{ij})$. By $(\nabla^2 G(e_i,e_j))\le 0$, we have $\sum\limits_i\nabla^2 G(e_i,e_i)=g^{ij}G_{ij}(x_0)\le 0$ and
\[
\begin{array}{lcl}
0&\ge&F^{i j}\nabla^2 G(e_i,e_j)(x_0)=F^{i j}a_i^r a_j^s G_{rs}=tr(B A^T D A)=tr(B A^TA D)\\
&=&tr(B {\cal G}^{-1}D)=F^{ir}g^{rj}G_{ij},\\
\end{array}
\]
i.e., we have $\bar{L}^{ij}G_{ij}(x_0)\le 0$.

In the following, we use $C_1>0$ to denote a universal constant depending only on $(M^n,g,t)$, $\phi$, $\psi$, $\delta_{y^{i_0}}$, $L$, $\alpha$, and we use $C_2>0$ to denote a universal constant depending only on $(M^n,g,t)$, $\phi$, $\psi$, $\delta_{y^{i_0}}$, $L$, $\delta_0$, $\alpha$, $\beta$. We also use $O_1(1)$ to denote a quantity bounded by $C_1$, and $O_2(1)$ to denote a quantity bounded by $C_2$. Observe that $\frac{1}{C_1}(\delta^{ij})\le{\cal G}^{-1}\le C_1(\delta^{ij})$ in $\overline{B^T_{\frac{\delta_{y^{i_0}}}{16}}(y^{i_0})}$. We will use this fact without mentioning.

At $x_0$,
\begin{equation}
\begin{array}{lcl}
0&\ge& \bar{L}^{i j} G_{i j}=2\rho\alpha \bar{L}^{ij}(v+\gamma)_{k,i}(v+\gamma)_{k,j}+2\rho\alpha \bar{L}^{ij}(v+\gamma)_k(v+\gamma)_{k,ij}\\
&&+\bar{L}^{i j}\Big(\frac{\alpha(\rho\rho_{ij}-2\rho_i\rho_j)}{\rho}+\frac{\rho(\alpha\alpha''-2(\alpha')^2)}{L^4\alpha}(v+\gamma)_i(v+\gamma)_j\\
&&+\frac{\alpha'\rho}{L^2}(v+\gamma)_{ij}-\frac{\alpha'}{L^2}(\rho_j(v+\gamma)_i+\rho_i(v+\gamma)_j)\Big)\sum\limits_k (v+\gamma)^2_k\\
&\ge&2\rho\alpha \bar{L}^{i j}(v+\gamma)_{k,i}(v+\gamma)_{k,j}+2\rho\alpha \bar{L}^{ij}(v+\gamma)_k(v+\gamma)_{k,ij}\\
&&+\frac{\rho(\alpha\alpha''-2(\alpha')^2)}{L^4\alpha}\sum\limits_k (v+\gamma)^2_k\bar{L}^{ij}(v+\gamma)_i(v+\gamma)_j\\
&&+\frac{\alpha'\rho}{L^2}\sum\limits_k (v+\gamma)^2_k\bar {L}^{ij}(v+\gamma)_{ij}-C_1\sqrt{\rho}\sum\limits_{k,l}F^{ll} |(v+\gamma)_k|^3-C_1\sum\limits_{k,l}F^{ll} (v+\gamma)^2_k\\
&\ge&2\rho\alpha \bar{L}^{i j}(v+\gamma)_{k,i}(v+\gamma)_{k,j}+2\rho\alpha \bar{L}^{ij}(v+\gamma)_k(v+\gamma)_{k,i j}\\
&&+\frac{\rho(\alpha\alpha''-2(\alpha')^2)}{L^4\alpha}\sum\limits_k (v+\gamma)^2_k\bar{L}^{ij}(v+\gamma)_i(v+\gamma)_j\\
&&+\frac{\alpha'\rho}{L^2}\sum\limits_k (v+\gamma)^2_k\bar {L}^{ij}(v+\gamma)_{ij}-C_1\sqrt{\rho}\sum\limits_{k,l}F^{ll} |(v+\gamma)_k|^3,\\
\end{array}
\label{c1en1}
\end{equation}
where in the last inequality, we used $G(x_0)\ge 1$, therefore $\sqrt{\rho}\sum\limits_k |(v+\gamma)_k|\ge \frac{1}{C_1}$.

In general,
\[
\begin{array}{lcl}
(v+\gamma)_{i j ,k}&=&\frac{\partial}{\partial y_k}\Big((v+\gamma)_{i,j}-\Gamma_{ji}^l(v+\gamma)_l\Big)\\
&=&(v+\gamma)_{i,j,k}-\Gamma_{ji}^l(v+\gamma)_{l,k}-\frac{\partial\Gamma_{ji}^l}{\partial y_k}(v+\gamma)_l,
\end{array}
\]
so
\[
\begin{array}{lcl}
(v+\gamma)_{k, i j}&=&(\nabla_{\frac{\partial}{\partial y_j}}\nabla_{\frac{\partial}{\partial y_i}}-\Gamma_{ji}^l\frac{\partial}{\partial y_l})((v+\gamma)_k)\\
&=&(v+\gamma)_{k,i,j}-\Gamma_{ji}^l (v+\gamma)_{k,l}\\
&=&(v+\gamma)_{i j ,k}+\frac{\partial\Gamma_{ji}^l}{\partial y_k}(v+\gamma)_l,\\
\end{array}
\]
and
\begin{equation}
\begin{array}{l c l}
0&\ge& \bar{L}^{i j} G_{i j}(x_0)\\
&\ge&2\rho\alpha \bar{L}^{i j}(v+\gamma)_{k,i}(v+\gamma)_{k,j}+2\rho\alpha \bar{L}^{ij}(v+\gamma)_k(v+\gamma)_{i j,k}\\
&&+\frac{\rho(\alpha\alpha''-2(\alpha')^2)}{L^4\alpha}\sum\limits_k (v+\gamma)^2_k\bar{L}^{ij}(v+\gamma)_i(v+\gamma)_j\\
&&+\frac{\alpha'\rho}{L^2}\sum\limits_k (v+\gamma)^2_k\bar {L}^{i j}
(v+\gamma)_{i j}-C_1\sqrt{\rho}\sum\limits_{k,l}F^{ll} |(v+\gamma)_k|^3.\\
\end{array}
\label{c1e1}
\end{equation}

Recall that $\gamma=(\psi e^v-h_g)\beta$. At $x_0$,
\begin{equation}
\begin{array}{lcl}
(v+\gamma)_k&=&(1+\psi\beta e^v)v_k+e^v\beta\psi_k+\psi e^v\beta_k-(h_g\beta)_k\\
&=&a v_k+O_2(1)\quad\mbox{with}~~a:=1+\psi\beta e^v.\\
\end{array}
\label{c1e2}
\end{equation}
\begin{equation}
\begin{array}{lcl}
\sum\limits_k (v+\gamma)_k^2&=& a^2 \sum\limits_k v_k^2+O_2(1)\sum\limits_k |v_k|.\\
\end{array}
\label{c1e3}
\end{equation}
\begin{equation}
\begin{array}{l c l}
(v+\gamma)_{i j}&=&av_{i j}+e^v\beta (v_i\psi_j+\psi_iv_j )+e^v\psi (v_i\beta_j+ v_j\beta_i)
+e^v(\psi_i\beta_j+\psi_j\beta_i)\\
&&+e^v\psi\beta v_i v_j+e^v\beta\psi_{i j}+e^v\psi\beta_{i j}-(h_g\beta)_{i j}\\
&=&a v_{i j}+O_2(1)\sum\limits_k |v_k|+O_1(1) \beta\sum\limits_k v_k^2.\\
\end{array}
\label{c1e4}
\end{equation}

The above identity (\ref{c1e4}) also holds for $(v+\gamma)_{i,j}$ after a slight modification, i.e., we only need to change $v_{ij},\psi_{ij},\beta_{ij},(h_g\beta)_{ij}$ to $v_{i,j},\psi_{i,j},\beta_{i,j},(h_g\beta)_{i,j}$ respectively.
\begin{equation}
\begin{array}{l c l}
(v+\gamma)_{i j, k}&=&av_{i j, k}+e^v(\psi_k\beta+\psi\beta_k)v_{i j}+e^v\psi\beta v_k v_{i j}
+e^v\psi\beta (v_j v_{i, k}+v_i v_{j ,k})\\
&&+\Big(e^v(\psi\beta_j+\psi_j\beta)v_{i, k}+e^v(\psi\beta_i+\psi_i\beta)v_{j, k} \Big)+e^v\psi\beta v_i v_j v_k\\
&&+e^v\Big(\psi v_i\beta_{j ,k}+\psi_i\beta_{j, k}+\psi_j\beta_{i, k}+\psi v_j\beta_{i, k}+\psi_k\beta_{i j}+\psi v_k\beta_{i j}+\psi\beta_{i j, k} \Big)\\
&&+e^v\beta\Big(v_i\psi_{j, k}+\psi_j v_i v_k+\psi_k v_i v_j+\psi_{i j, k}+\psi_{i j}v_k+\psi_{i, k}v_j+\psi_i v_j v_k \Big)\\
&&e^v\Big(\psi_j v_i\beta_k+\psi_k v_i\beta_j+\psi v_i v_k\beta_j+\psi v_i v_j\beta_k+\psi_{i j}\beta_k+\psi_i v_j\beta_k \\
&&\psi_{i, k}\beta_j+\psi_i v_k\beta_j+\psi_{j, k}\beta_i+\psi_j v_k\beta_i+\psi_k v_j\beta_i+\psi v_j v_k\beta_i\Big)-(h_g\beta)_{i j, k}\\
&=&av_{i j, k}+e^v(\psi_k\beta+\psi\beta_k)v_{i j}+e^v\psi\beta v_k v_{i j}+e^v\psi\beta v_i v_j v_k\\
&&+\Big(e^v(\psi\beta_j+\psi_j\beta+\psi\beta v_j)v_{i, k}+e^v(\psi\beta_i+\psi_i\beta+\psi\beta v_i)v_{j, k} \Big)\\
&&+O_2(1) \sum\limits_k v_k^2.\\
\end{array}
\label{c1e5}
\end{equation}

By (\ref{c1e0}) and (\ref{c1e2}-\ref{c1e4}), at $x_0$,
\[
\begin{array}{l c l}
&&(v+\gamma)_k\Big(a v_{k, i }+O_2(1)\sum\limits_l |v_l|+O_1(1)\beta \sum\limits_l v_l^2\Big)\\
&=&-\frac{\alpha'}{2 L^2\alpha}\Big(a^2\sum\limits_l v_l^2+O_2(1)\sum\limits_l |v_l| \Big)(a v_i+O_2(1))\\
&&-\frac{\rho_i}{2\rho}\Big(a^2\sum\limits_l v_l^2+O_2(1)\sum\limits_l |v_l| \Big),\\
\end{array}
\]
therefore
\[
\begin{array}{l c l}
&&a(v+\gamma)_k v_{k, i }+(av_k+O_2(1))\Big(O_2(1)\sum\limits_l |v_l|+O_1(1)\beta \sum\limits_l v_l^2\Big)\\
&=&-\frac{\alpha'}{2 L^2\alpha}a^3 \sum\limits_l v_l^2 v_i+O_2(1)\frac{1}{\sqrt{\rho}}\sum\limits_l v_l^2,\\
\end{array}
\]
which implies that
\begin{equation}
(v+\gamma)_k v_{k, i}=-\frac{\alpha'}{2 L^2\alpha}a^2\sum\limits_l v_l^2 v_i+O_1(1)\beta \sum\limits_{l} |v_l|^3 +O_2(1)\frac{1}{\sqrt{\rho}}\sum\limits_l v_l^2 ,
\label{c1e6}
\end{equation}
where we used $a=1+\psi\beta e^v\in [\frac 12, 1]$.

Combine (\ref{c1e2}-\ref{c1e6}). At $x_0$,

\begin{eqnarray}
&&2\alpha \rho(v+\gamma)_k\bar{L}^{i j}(v+\gamma)_{i j, k}\ge\nonumber{}\\
&&2\alpha a\rho(v+\gamma)_k \bar{L}^{i j} v_{i j, k}+2\alpha \rho e^v(v_k+\gamma_k)(\psi_k\beta+\psi\beta_k+\psi\beta v_k)\bar{L}^{i j}v_{i j}\nonumber{}\\
&&+4\alpha \rho e^v \bar{L}^{i j}(\psi\beta_j+\psi_j\beta)(v_k+\gamma_k)v_{i, k}+4\alpha \rho e^v\psi\beta \bar{L}^{i j}v_j (v_k+\gamma_k)v_{i, k}\nonumber{}\\
&&+2\alpha \rho e^v\psi\beta (v_k+\gamma_k)\bar{L}^{i j} v_i v_j v_k-C_2 \rho\sum\limits_{k, l}F^{ll}|v_k|^3\nonumber{}\\
&&\ge 2\alpha a\rho(v+\gamma)_k \bar{L}^{i j} v_{i j,k}+2\alpha \rho e^v(v+\gamma)_k(\psi_k\beta+\psi\beta_k+\psi\beta v_k)\bar{L}^{i j}v_{i j}\nonumber{}\\
&&+4\alpha \rho e^v \bar{L}^{i j}(\psi\beta_j+\psi_j\beta)\Big( -\frac{\alpha'}{2 L^2\alpha}a^2\sum\limits_l v_l^2 v_i+O_1(1)\beta \sum\limits_l|v_l|^3\nonumber{}\\
&& +O_2(1)\sum\limits_l v_l^2\Big)\nonumber{}\\
&&+4\alpha \rho e^v\psi\beta \bar{L}^{i j}v_j \Big( -\frac{\alpha'}{2 L^2\alpha}a^2\sum\limits_l v_l^2 v_i+O_1(1)\beta \sum\limits_l |v_l|^3
+O_2(1)\sum\limits_l v_l^2\Big)\nonumber{}\\
&&+2\alpha \rho e^v\psi\beta (av_k+O_2(1))\bar{L}^{i j} v_i v_j v_k-C_2 \sum\limits_{k,l}F^{ll}|v_k|^3\nonumber{}\\
&&\ge 2\alpha a\rho (v+\gamma)_k\bar{L}^{i j} v_{i j,k}+2\alpha \rho e^v(v+\gamma)_k(\psi_k\beta+\psi\beta_k+\psi\beta v_k)\bar{L}^{i j}v_{i j}\nonumber{}\\
&&-C_1\beta \rho\sum\limits_{k,l}F^{ll}v_k^4-C_2 \rho\sum\limits_{k,l}F^{ll}|v_k|^3,\nonumber{}\\
\label{c1e7}
\end{eqnarray}

Recall the Laplace-Beltrami operator $\Delta_g=\frac{1}{\sqrt{|g|}}\frac{\partial}{\partial y_k}(\sqrt{|g|}g^{km}\frac{\partial}{\partial y_m})$, where $|g|=\det (g_{km})$.
\[
\Delta_g v=g^{km}v_{m,k}+\frac{1}{\sqrt{|g|}}(\sqrt{|g|}g^{km})_k v_m=g^{km}v_{m k}+g^{k m}\Gamma_{k m}^lv_l+\frac{1}{\sqrt{|g|}}(\sqrt{|g|}g^{km})_k v_m.
\]

Since $f$ is homogeneous of degree $1$, the equation $F(\bar {W}_{ij}g^{jr})=\phi e^{2v}$ implies that
\[
\begin{array}{lcl}
\phi e^{2v}&=& F^{i r}g^{j r}\Big(v_{i j}+\frac{1-t}{n-2}(\Delta_g v)g_{i j}+\frac{2-t}{2}|\nabla v|_g^2 g_{i j}-v_iv_j-(A_g^t)_{ij} \Big)\\
&=&F^{i r}g^{j r}v_{i j}+\frac{1-t}{n-2}(\Delta_g v)\sum\limits_l F^{ll}+\frac{2-t}{2}|\nabla v|_g^2 \sum\limits_l F^{ll}-F^{ir}g^{jr}v_iv_j-F^{ir}g^{jr}(A_g^t)_{ij}\\
&=&F^{i r}g^{j r}v_{i j}+\frac{1-t}{n-2}\Big( g^{km}v_{km}+g^{k m}\Gamma_{k m}^r v_r+\frac{1}{\sqrt{|g|}}(\sqrt{|g|}g^{km})_k v_m)\Big)\sum\limits_l F^{ll}\\
&&+\frac{2-t}{2}v_k v_lg^{kl} \sum\limits_l F^{ll}-F^{ir}g^{jr}v_iv_j-F^{ir}g^{jr}(A_g^t)_{ij}\\
&=&\bar{L}^{ij}v_{ij}+\frac{1-t}{n-2}g^{k m}\Gamma_{k m}^r v_r\sum\limits_l F^{ll}+\frac{1-t}{n-2}\frac{1}{\sqrt{|g|}}(\sqrt{|g|}g^{km})_k v_m\sum\limits_l F^{ll}\\
&&+\frac{2-t}{2}v_k v_lg^{kl} \sum\limits_l F^{ll}-F^{ir}g^{jr}v_iv_j-F^{ir}g^{jr}(A_g^t)_{ij},\\
\end{array}
\]
that is,
\[
\bar{L}^{ij}v_{ij}=F^{ir}g^{jr}v_i v_j-\frac{2-t}{2}v_k v_l g^{kl}\sum\limits_{l}F^{ll}+O_1(1)\sum\limits_{k,l}F^{ll}|v_k|.
\]

From which, we have
\begin{eqnarray}
&&\frac{\alpha'\rho}{L^2}\sum\limits_k (v+\gamma)_k^2 \bar{L}^{i j}(v+\gamma)_{i j}\nonumber{}\\
&=&\frac{\alpha'\rho}{L^2}\Big(a^2\sum\limits_k v_k^2+O_2(1)\sum\limits_k |v_k|
 \Big)\bar{L}^{i j}\Big( a v_{i j}+O_2(1)\sum\limits_k |v_k|\nonumber{}\\
&&+O_1(1)\beta\sum\limits_k v_k^2 \Big)~~~~~~~~~~~~~~~~~~~~\qquad\mbox{by}~~(\ref{c1e3})~~\mbox{and}~~(\ref{c1e4})\nonumber{}\\
&\ge& \frac{a\alpha'\rho}{L^2}\Big( a^2\sum\limits_k v_k^2+O_2(1)\sum\limits_k |v_k|
 \Big)\bar{L}^{ij} v_{i j}\nonumber{}\\
&&-C_1\beta \rho \sum\limits_{k,l}F^{ll}v_k^4
-C_2\rho\sum\limits_{k,l}F^{ll}|v_k|^3\nonumber{}\\
&\ge& \frac{a\alpha'\rho}{L^2}\Big( a^2\sum\limits_k v_k^2+O_2(1)\sum\limits_k |v_k|
 \Big)\Big( F^{ir}g^{jr}v_i v_j-\frac{2-t}{2}v_k v_l g^{kl}\sum\limits_{i}F^{ii}\nonumber{}\\
&&+O_1(1)\sum\limits_{k,l}F^{ll}|v_k|\Big)-C_1\beta \rho \sum\limits_{k,l}F^{ll}v_k^4-C_2\rho\sum\limits_{k,l}F^{ll}|v_k|^3\nonumber{}\\
&\ge& \frac{a^3\alpha'\rho}{L^2}\sum\limits_k v_k^2 F^{ir}g^{jr}v_i v_j-\frac{(2-t)a^3\alpha'\rho}{2 L^2}v_k v_l g^{kl}\sum\limits_{i,j}F^{ii}v_j^2-C_1\beta \rho \sum\limits_{k,l}F^{ll}v_k^4\nonumber{}\\
&&-C_2\rho\sum\limits_{k,l}F^{ll}|v_k|^3\nonumber{}\\
\label{c1e8}
\end{eqnarray}
and
\[
2\alpha \rho e^v(v+\gamma)_k(\psi_k\beta+\psi\beta_k+\psi\beta v_k)\bar{L}^{i j}v_{i j}\ge-C_1\beta \rho \sum\limits_{k,l}F^{ll}v_k^4 -C_2\rho \sum\limits_{k,l}F^{ll}|v_k|^3.\\
\]
which implies, by (\ref{c1e7}), that
\begin{equation}
\begin{array}{l c l}
&&2\alpha \rho(v+\gamma)_k \bar{L}^{i j}(v+\gamma)_{i j,k}\ge 2\alpha a \rho(v+\gamma)_k \bar{L}^{i j} v_{i j,k}\\
&&-C_1\beta \rho\sum\limits_{k,l}F^{ll}v_k^4-C_2 \rho\sum\limits_{k,l}F^{ll}|v_k|^3,\nonumber{}\\
\end{array}
\label{c1e9}
\end{equation}

Differentiate the equation $F(\bar {W}_{ij}g^{jr})=\phi e^{2v}$ along the $y_k-th$ direction and evaluate at $x_0$.

\[
\begin{array}{lcl}
&&\phi_k e^ {2v}+2\psi e^{2v} v_k=F^{ir}\Big(g^{jr}\bar{W}_{ij} \Big)_k=F^{ir} g^{jr}(\bar{W}_{ij})_k+\frac{\partial g^{jr}}{\partial y_k}F^{ir}\bar{W}_{ij}\\
&=&F^{ir}g^{jr}\Big( v_{i j,k}+\frac{1-t}{n-2}(\Delta_g v)_k g_{i j}+\frac{1-t}{n-2}\frac{\partial g_{ij}}{\partial y_k}(\Delta_g v)+\frac{2-t}{2}(2v_{m,k}v_lg^{ml}\\
&&+v_m v_l\frac{\partial g^{ml}}{\partial y_k})g_{i j}+\frac{2-t}{2}v_m v_l g^{ml}\frac{\partial g_{i j}}
{\partial y_k}-2v_{i,k} v_j-(A_g^t)_{i j,k}\Big)\\
&&+\frac{\partial g^{jr}}{\partial y_k}F^{ir} \Big(v_{i j}+\frac{1-t}{n-2}(\Delta_g v)g_{i j}+\frac{2-t}{2}|\nabla v|_g^2 g_{i j}-v_iv_j-(A_g^t)_{ij} \Big)\\
&=&F^{ir}g^{jr} v_{i j,k}+\frac{1-t}{n-2}(\Delta_g v)_k \sum\limits_l F^{ll}+\frac{1-t}{n-2}F^{ir}g^{jr}\frac{\partial g_{ij}}{\partial y_k}(\Delta_g v)\\
&&+(2-t)v_{m,k}v_lg^{ml}\sum\limits_i F^{ii}-2 F^{ir}g^{jr}v_{i,k} v_j+\frac{\partial g^{jr}}{\partial y_k}F^{ir} \Big(v_{i j}\\
&&+\frac{1-t}{n-2}(\Delta_g v)g_{i j}\Big)+O_1(1)\sum\limits_{i,j}F^{ii}v_j^2\\
&=&F^{ir}g^{jr} v_{i j,k}+\frac{1-t}{n-2}\Big( g^{lm}v_{lm}+g^{l m}\Gamma_{l m}^r v_r+\frac{1}{\sqrt{|g|}}(\sqrt{|g|}g^{lm})_l v_m\Big)_k \sum\limits_i F^{ii}\\
&&+\frac{1-t}{n-2}F^{ir}g^{jr}\frac{\partial g_{ij}}{\partial y_k}\Big( g^{l m}v_{l m}+g^{l m}\Gamma_{l m}^r v_r+\frac{1}{\sqrt{|g|}}(\sqrt{|g|}g^{lm})_l v_m\Big)\\
&&+(2-t)v_{m,k}v_lg^{ml}\sum\limits_i F^{ii}-2 F^{ir}g^{jr}v_{i,k} v_j\\
&&+\frac{\partial g^{jr}}{\partial y_k}F^{ir} \Big(v_{i j}+\frac{1-t}{n-2}( g^{lm}v_{lm}+g^{l m}\Gamma_{l m}^r v_r+\frac{1}{\sqrt{|g|}}(\sqrt{|g|}g^{l m})_l v_m)g_{i j}\Big)\\
&&+O_1(1)\sum\limits_{i,j}F^{ii}v_j^2\\
&=&F^{ir}g^{jr} v_{i j,k}+\frac{1-t}{n-2}g^{l m}v_{l m,k} \sum\limits_i F^{ii}+(2-t)v_{m,k}v_lg^{ml}\sum\limits_i F^{ii}\\
&&-2 F^{ir}g^{jr}v_{i,k} v_j+O_1(1)\sum\limits_{i,j,l}F^{ii}|v_{jl}|+O_1(1)\sum\limits_{i,j}F^{ii}(v_j^2+1)\\
&=&\bar{L}^{i j}v_{i j,k}+(2-t)v_{m,k}v_lg^{ml}\sum\limits_i F^{ii}-2 F^{i r}g^{j r}v_{i,k} v_j\\
&&+O_1(1)\sum\limits_{i,j,l}F^{ii}|v_{j, l}+\Gamma_{lj}^rv_r|+O_1(1)\sum\limits_{i,j}F^{ii}v_j^2\\
&=&\bar{L}^{i j}v_{i j,k}+(2-t)v_{m,k}v_lg^{ml}\sum\limits_i F^{ii}-2 F^{i r}g^{j r}v_{i,k} v_j\\
&&+O_1(1)\sum\limits_{i,j,l}F^{ii}|v_{j, l}|+O_1(1)\sum\limits_{i,j}F^{ii}v_j^2.\\
\end{array}
\]

Multiply both sides by $2\alpha a\rho(v+\gamma)_k$ and solve it for $2\alpha a\rho(v+\gamma)_k\bar{L}^{ij}v_{i j,k}$.
\begin{eqnarray*}
&&2\alpha a\rho(v+\gamma)_k \bar{L}^{i j}v_{i j,k}\\
&=&-2(2-t)\alpha a\rho(v+\gamma)_k v_{m,k}v_lg^{ml}\sum\limits_i F^{ii}+4\alpha a\rho(v+\gamma)_k v_{i,k} F^{i r}g^{j r} v_j\\
&&+O_1(1)\rho\sum\limits_{i,j,k,l}F^{ii}|v_{j, l}||(v+\gamma)_k|+O_1(1)\rho\sum\limits_{i,j,k}F^{ii} v_j^2|(v+\gamma)_k|\\
&=&-2(2-t)\alpha a\rho\Big( -\frac{\alpha'}{2 L^2\alpha}a^2\sum\limits_j v_j^2 v_m+O_1(1)\beta \sum\limits_{j} |v_j|^3 \\
&&+O_2(1)\frac{1}{\sqrt{\rho}}\sum\limits_j v_j^2\Big) v_lg^{ml}\sum\limits_i F^{ii}+4\alpha a\rho\Big(-\frac{\alpha'}{2 L^2\alpha}a^2\sum\limits_l v_l^2 v_i\\
&&+O_1(1)\beta \sum\limits_{l} |v_l|^3 +O_2(1)\frac{1}{\sqrt{\rho}}\sum\limits_l v_l^2\Big) F^{i r}g^{j r} v_j~~\mbox{by}~~(\ref{c1e6})\\
&&+O_1(1)\rho\sum\limits_{i,j,k,l}F^{ii}|v_{j, l}||v_k|+O_2(1)\sqrt{\rho}\sum\limits_{i,j}F^{ii}|v_j|^3\\
&\ge&\frac{(2-t)a^3\alpha'\rho}{ L^2}\sum\limits_{i,j} F^{i i}v_j^2 v_m v_lg^{ml}-\frac{2 a^3\rho\alpha'}{ L^2}\sum\limits_l v_l^2 F^{i r}g^{j r} v_i v_j\\
&&-C_1\beta\rho\sum\limits_{i,j}F^{ii}v_j^4-C_1\rho\sum\limits_{i,j,k,l}F^{ii}|v_{j, l}||v_k|-C_2\sqrt{\rho}\sum\limits_{i,j}F^{ii}|v_j|^3.\\
\end{eqnarray*}

Substitute the above inequality into (\ref{c1e9}).

\begin{equation}
\begin{array}{lcl}
&&2\alpha (v+\gamma)_k \bar{L}^{i j}(v+\gamma)_{i j,k}\\
&\ge& \frac{(2-t)a^3\alpha'\rho}{ L^2}\sum\limits_{i,j} F^{i i}v_j^2 v_m v_lg^{ml}-\frac{2 a^3\rho\alpha'}{ L^2}\sum\limits_l v_l^2 F^{i r}g^{j r} v_i v_j\\
&&-C_1\beta\rho\sum\limits_{i,j}F^{ii}v_j^4-C_1\rho\sum\limits_{i,j,k,l}F^{ii}|v_{j, l}||v_k|-C_2\sqrt{\rho}\sum\limits_{i,j}F^{ii}|v_j|^3.\\
\end{array}
\label{c1e10}
\end{equation}

By (\ref{c1e2}) and (\ref{c1e3}),
\begin{equation}
\begin{array}{l c l}
&&\frac{\rho(\alpha\alpha''-2(\alpha')^2)}{L^4\alpha}\sum\limits_k(v+\gamma)^2_k \bar{L}^{i j}(v+\gamma)_i(v+\gamma)_j
-C_1 \sqrt{\rho}\sum\limits_{k,l} F^{ll}|(v+\gamma)_k|^3\\
&\ge& \frac{\rho(\alpha\alpha''-2(\alpha')^2)}{L^4\alpha}\Big(a^2 v_k^2+O_2(1)|v_k| \Big)
\bar{L}^{ij}(av_i+O_2(1))(av_j+O_2(1))\\
&&-C_1 \sqrt{\rho}\sum\limits_{k,l} F^{ll}\Big(a^3 |v_k|^3 +O_2(1)v_k^2 \Big)\\
&\ge&\frac{\alpha''\alpha-2(\alpha')^2}{L^4\alpha} \rho a^4 \sum\limits_k v_k^2\bar{L}^{i j}v_i v_j-C_2 \sqrt{\rho}\sum\limits_{k,l}F^{ll}|v_k|^3,\\
\end{array}
\label{c1e11}
\end{equation}
and

\begin{equation}
\begin{array}{lcl}
&&2\rho\alpha\bar{L}^{i j}(v+\gamma)_{k,i}(v+\gamma)_{k,j}\\
&=&2\rho\alpha\bar{L}^{i j}\Big(av_{k,i}+O_1(1)\beta\sum\limits_l v_l^2+O_2(1)\sum\limits_l |v_l| \Big)\Big(av_{k,j}\\
&&+O_1(1)\beta\sum\limits_l v_l^2+O_2(1)\sum\limits_l |v_l| \Big)\\
&\ge& 2 a^2\alpha\rho\bar{L}^{i j} v_{k,i}v_{k,j}-C_1\beta\rho\sum\limits_{i,j,k,l}F^{ii}v_j^2 |v_{k,l}|-C_2\rho\sum\limits_{i,j,k,l}F^{ii}|v_j||v_{k,l}|\\
&&-C_1\beta \rho\sum\limits_{i,j}F^{ii}v_j^4-C_2 \rho\sum\limits_{i,j}F^{ii}|v_j|^3.\\
\end{array}
\label{c1e11-1}
\end{equation}

Substitute (\ref{c1e8}), (\ref{c1e10}), (\ref{c1e11}), and (\ref{c1e11-1}) into (\ref{c1e1}). We have,

\[
\begin{array}{l c l}
&&0\ge \bar{L}^{i j}G_{i j}(x_0)\ge \frac{(2-t)a^3\alpha'\rho}{2 L^2}\sum\limits_{i,j} F^{i i}v_j^2 v_m v_lg^{ml}-\frac{ a^3\rho\alpha'}{ L^2}\sum\limits_l v_l^2 F^{i r}g^{j r} v_i v_j\\
&&+\frac{\alpha''\alpha-2(\alpha')^2}{L^4\alpha} \rho a^4 \sum\limits_k v_k^2\bar{L}^{i j}v_i v_j+2 a^2\alpha\rho\bar{L}^{ij} v_{k,i}v_{k,j}-C_1\beta\rho\sum\limits_{i,j}F^{ii}v_j^4\\
&&-C_1\beta\rho\sum\limits_{i,j,k,l}F^{ii}v_j^2|v_{k, l}|
-C_2\rho\sum\limits_{i,j,k,l}F^{ii}|v_{j, l}||v_k|-C_2\sqrt{\rho}\sum\limits_{i,j}F^{ii}|v_j|^3.\\
\end{array}
\]

Recall $a=1+\psi\beta e^v$. We can replace it by $1+O_1(1)\beta$. Meanwhile we replace $\bar{L}^{ij}$ by $F^{ir}g^{rj}+\frac{1-t}{n-2}(\sum\limits_l F^{ll})g^{ij}$ in the above inequality. We have,
\begin{eqnarray}
0&\ge& \bar{L}^{i j}G_{i j}(x_0)\ge\frac{(2-t)\alpha'\rho}{2 L^2}\sum\limits_{i,j} F^{i i}v_j^2 v_m v_lg^{ml}-\frac{ \rho\alpha'}{ L^2}\sum\limits_l v_l^2 F^{i r}g^{j r} v_i v_j\nonumber{}\\
&&+\frac{\rho(\alpha''\alpha-2(\alpha')^2)}{L^4\alpha} \sum\limits_k v_k^2 \Big( F^{ir}g^{rj}+\frac{1-t}{n-2}(\sum\limits_l F^{ll})g^{ij}\Big)v_i v_j\nonumber{}\\
&&+2 a^2\alpha\rho \Big( F^{ir}g^{rj}+\frac{1-t}{n-2}(\sum\limits_l F^{ll})g^{ij}\Big)v_{k,i}v_{k,j}-C_1\beta\rho\sum\limits_{i,j}F^{ii}v_j^4\nonumber{}\\
&&-C_1\beta\rho\sum\limits_{i,j,k,l}F^{ii}v_j^2|v_{k, l}|
-C_2\rho\sum\limits_{i,j,k,l}F^{ii}|v_{j, l}||v_k|-C_2\sqrt{\rho}\sum\limits_{i,j}F^{ii}|v_j|^3\nonumber{}\\
&\ge&\rho\Big(\frac{(2-t)\alpha'}{2 L^2}+ \frac{(1-t)(\alpha\alpha''-2(\alpha')^2)}{(n-2)L^4\alpha}\Big)\sum\limits_{i,j} F^{i i}v_j^2 v_m v_lg^{ml}\nonumber{}\\
&&+\rho\Big(\frac{\alpha''\alpha-2(\alpha')^2}{L^4\alpha}-\frac{ \alpha'}{ L^2} \Big)\sum\limits_l v_l^2 F^{i r}g^{j r} v_i v_j\nonumber{}\\
&&+\frac{2(1-t) a^2\alpha\rho}{n-2}(\sum\limits_l F^{ll})(\frac{1}{C_1}\delta^{ij})v_{k,i}v_{k,j}-C_1\beta\rho\sum\limits_{i,j}F^{ii}v_j^4\nonumber{}\\
&&-C_1\beta\rho\sum\limits_{i,j,k,l}F^{ii}v_j^2|v_{k, l}|
-C_2\rho\sum\limits_{i,j,k,l}F^{ii}|v_{j, l}||v_k|-C_2\sqrt{\rho}\sum\limits_{i,j}F^{ii}|v_j|^3\nonumber{}\\
&\ge&\rho\Big(\frac{(2-t)\alpha'}{2 L^2}+ \frac{(1-t)(\alpha\alpha''-2(\alpha')^2)}{(n-2)L^4\alpha}\Big)\sum\limits_{i,j} F^{i i}v_j^2 v_m v_lg^{ml}\nonumber{}\\
&&+\rho\Big(\frac{\alpha''\alpha-2(\alpha')^2}{L^4\alpha}-\frac{ \alpha'}{ L^2} \Big)\sum\limits_l v_l^2 F^{i r}g^{j r} v_i v_j\nonumber{}\\
&&+\frac{(1-t) a^2\alpha\rho}{(n-2)C_1}(\sum\limits_l F^{ll})v_{k,i}^2-C_1\beta\rho\sum\limits_{i,j}F^{ii}v_j^4-C_2\sqrt{\rho}\sum\limits_{i,j}F^{ii}|v_j|^3\nonumber{}\\
&\ge&\rho\Big(\frac{(2-t)\alpha'}{2 L^2}+ \frac{(1-t)(\alpha\alpha''-2(\alpha')^2)}{(n-2)L^4\alpha}\Big)\sum\limits_{i,j} F^{i i}v_j^2 v_m v_lg^{ml}\nonumber{}\\
&&+\rho\Big(\frac{\alpha''\alpha-2(\alpha')^2}{L^4\alpha}-\frac{ \alpha'}{ L^2} \Big)\sum\limits_l v_l^2 F^{i r}g^{j r} v_i v_j\nonumber{}\\
&&-C_1\beta\rho\sum\limits_{i,j}F^{ii}v_j^4-C_2\sqrt{\rho}\sum\limits_{i,j}F^{ii}|v_j|^3.\nonumber{}\\
\label{c1e12}
\end{eqnarray}

It is enough to find a smooth function $\alpha:~[\frac{1}{L^2},\frac{2}{L}]\to\Bbb{R}^+$ satisfying
\begin{equation}
\{
\begin{array}{l c l}
\alpha'>0&&\\
\alpha\alpha''-2(\alpha')^2-L^2\alpha\alpha'>0.&&\\
\end{array}
\label{c1e13}
\end{equation}
since the above inequalities imply that
\[
\frac{\alpha''\alpha-2(\alpha')^2}{L^4\alpha}-\frac{\alpha'}{L^2}=\frac{1}{L^4\alpha}(\alpha\alpha''-2(\alpha')^2-L^2\alpha\alpha')>0,
\]
and
\[
\alpha\alpha''-2(\alpha')^2> L^2\alpha\alpha'>0,
\]
and
\[
\begin{array}{lcl}
\frac{(2-t)\alpha'}{2 L^2}+\frac{1-t}{n-2}\frac{\alpha''\alpha-2(\alpha')^2}{L^4\alpha}&>&0,\\
\end{array}
\]
i.e., the coefficients of the two leading terms in the inequality (\ref{c1e12}) are both positive, which will lead the preferred gradient bound.

Let $\alpha=e^{\eta}$. The two inequalities in (\ref{c1e13}) are equivalent to
\[
\{
\begin{array}{l c l}
\eta'>0&&\\
\eta''-(\eta')^2-L^2\eta'>0.&&\\
\end{array}
\]

To find $\alpha$, let $\eta (s)=s^r$ with $r\gg 1$ being chosen later. Clearly, $\eta'>0$ and
\begin{eqnarray*}
\eta''-(\eta')^2-L^2\eta'&=&rs^{r-2}\Big((r-1)-rs^r-L^2 s \Big)\\
&\ge& rs^{r-2}\Big((r-1)-r(\frac{2}{L})^r-L^2 (\frac{2}{L}) \Big)\\
&=&rs^{r-2}\Big((r-1)-r(\frac{2}{L})^r
-2 L \Big)\\
&\ge& rs^{r-2}\Big((r-1)-\frac{r}{2} -2 L \Big)\quad\mbox{by choosing}~~L>4\\
&=&rs^{r-2}\Big(\frac{r}{2}-1 -2 L \Big)\ge rs^{r-2}>0\quad\mbox{by choosing}~~r>4+4 L.
\end{eqnarray*}

Pick $L>|v+\gamma|+4$ and $r>4+4L$. Then we have $\frac{v+\gamma+L}{L^2}\in [\frac{1}{L^2},\frac{2}{L}]$ and there exists a universal constant $C_3>0$ independent of $\beta$ such that (\ref{c1e13}) holds. By (\ref{c1e12}),

\[
\begin{array}{l c l}
0&\ge& \bar{L}^{i j}G_{i j}(x_0)\ge C_3\rho\sum\limits_{i,j,k,l} v_l^2 F^{k k}(\frac{1}{C_1}\delta^{i j}) v_i v_j-C_1\beta\rho\sum\limits_{i,j}F^{ii}v_j^4-C_2\sqrt{\rho}\sum\limits_{i,j}F^{ii}|v_j|^3\\
&\ge&C_3\rho\sum\limits_{k,l} v_l^4 F^{k k} -C_1\beta\rho\sum\limits_{i,j}F^{ii}v_j^4-C_2\sqrt{\rho}\sum\limits_{i,j}F^{ii}|v_j|^3\\
&\ge&\frac{C_3}{2}\rho\sum\limits_{k,l} v_l^2 F^{k k} -C_2\sqrt{\rho}\sum\limits_{i,j}F^{ii}|v_j|^3,\\
\end{array}
\]
where in the last inequality, we used $\beta\in [0,2\delta_0]$, therefore we can pick $\delta_0\ll 1$ such that $C_1\beta\le\frac{C_3}{2}$.

We conclude that
\[
\begin{array}{lcl}
0&\ge&\frac{C_3}{2}\rho\sum\limits_{k,l} v_l^4 F^{k k}  -C_2\sum\limits_{i,j}F^{j j}|v_i|^3\ge
C_3\rho\sum\limits_k F^{k k}(\sum\limits_{l} v_l^2)^2  -C_2\sum\limits_j F^{j j}(\sum\limits_{i}v_i^2)^{\frac{3}{2}}\\
&=&\rho\sum\limits_{k,l} F^{k k}v_l^2\Big(C_3\sum\limits_i v_i^2-C_2\sum\limits_{i}v_i^2)^{\frac{1}{2}} \Big),\\
\end{array}
\]
which implies that $\sum\limits_iv_i^2\le C$, therefore $G(x_0)\le C$. In particular $\sum\limits_{i}v_i^2\le C$ in $B^T_{\frac{\delta_{y^{i_0}}}{32}}(y^{i_0})$. From which, we have, in $B^T_{\frac{\delta_{y^{i_0}}}{32}}(y^{i_0})$,
\[
|\nabla v|_g^2=v_k v_lg^{k l}\le C\sum\limits_k v_k^2\le C.
\]

By $\cup_{i_0=1}^N\Big(B^T_{\frac{\delta_{y^{i_0}}}{64}}(y^{i_0})\cap\partial M\Big)=\partial M$, $|\nabla v|_g^2\le C$ on $\partial M$. The Lemma~\ref{c1l1} has been
established. $\clubsuit$
\vskip 5 pt
\begin{remark}
When the manifold $(M^n,g)$ is umbilic on the boundary, the above lemma and therefore the next lemma also hold for $t=1$. The above proof still works after a slight modification.
\end{remark}
\begin{lemma}
Under the same assumptions as in Theorem~\ref{generalresult}, for $t<1$, let $v$ be a $C^3$ solution of the equation (\ref{generalproblem}). Then there exists a universal constant $C>0$ depending only on $(M^n,g,t)$, $(f,\Gamma)$, $\phi$, and $\psi$, such that
\[
|\nabla v|_g\le C\qquad\mbox{on}~~M.
\]
\label{c1l2}
\end{lemma}
{\bf Proof of the Lemma~\ref{c1l2}.} Consider
\[
\bar{G}:=|\nabla v|_g^2\bar{\alpha} (\frac{v+L}{L^2}),
\]
where $L>0$ is a constant satisfying $1<v+L<2L$ and $\alpha:~\Bbb{R}^+\to\Bbb{R}^+$ is a smooth positive function to be chosen later. Let $\bar{G}(x_0)=\max\limits_{M}G$. Let $\{x_j\}_{j=1}^n$ be a geodesic normal coordinates w.r.t. the metric $g$ at $x_0$. W.l.o.g., we can assume $x_0$ is an interior point of $M$. In the following, subindices are taken w.r.t. $\frac{\partial}{\partial x_j}$. Repeat the arguments in the proof of the Lemma~\ref{c1l1}. We arrive at
\[
\begin{array}{lcl}
0&\ge& \bar{L}^{ij}\bar{G}_{ij}(x_0)\ge\Big(\frac{(2-t)\alpha'}{2 L^2}+ \frac{(1-t)(\alpha\alpha''-2(\alpha')^2)}{(n-2)L^4\alpha}\Big)|\nabla v|_g^4\sum\limits_{i} F^{i i}\\
&&+\Big(\frac{\alpha''\alpha-2(\alpha')^2}{L^4\alpha}-\frac{ \alpha'}{ L^2} \Big)|\nabla v|_g^2 F^{i j} v_i v_j-C|\nabla v|^3_g\sum\limits_{i}F^{ii}.\\
\end{array}
\]

Choose the same $\alpha$ as in the proof of the Lemma~\ref{c1l1}. We conclude that there exists some universal constant $C_3>0$ such that
\[
\begin{array}{lcl}
0&\ge& \bar{L}^{i j}\bar{G}_{i j}(x_0)\ge C_3 |\nabla v|_g^4\sum\limits_{i} F^{i i}
-C|\nabla v|^3_g\sum\limits_{i}F^{ii}\\
&\ge& |\nabla v|_g^3\sum\limits_i F^{i i}(C_3|\nabla v|_g-C),\\
\end{array}
\]
which implies that $|\nabla v|_g(x_0)\le C$ and therefore $G(x_0)\le C$. The Lemma~\ref{c1l2} has been proved. $\clubsuit$
\begin{section}
{Hessian Estimates}
\end{section}

The main issue of the Hessian estimates is to bound the Hessian of the solutions on the boundary of $M$.

\begin{lemma}
Under the same assumptions as in Theorem~\ref{generalresult}, for $t<1$, let $v$ be a $C^4$ solution of the equation (\ref{generalproblem}). For any $1\le i_0\le N$, there exists a universal constant $C>0$ depending only on $(M^n,g,t)$, $(f,\Gamma)$, $\phi$, $\psi$, and $\delta_{y^{i_0}}$ such that in $B^T_{\frac{\delta_{y^{i_0}}}{32}}(y^{i_0})$,
\[
|v_{\tau\tau}|<C,\quad\mbox{for any unit direction $\frac{\partial}{\partial\tau}$ satisfying}~g(\frac{\partial}{\partial\tau},\frac{\partial}{\partial\nu})=0.
\]
\label{c2l1}
\end{lemma}
{\bf Proof of the Lemma~\ref{c2l1}.} Consider
\[
\bar {H}(y):=\rho e^{\beta_0 y_n}\Big(\{\max\limits_{\begin{array}{ccc}
\tau\in T_y M,~\|\tau\|_g=1,&&\\
g(\frac{\partial}{\partial\nu},\frac{\partial}{\partial\tau})=0&&\\
\end{array}}
(\nabla^2 v+a|\nabla (v+\gamma)|_g^2 g)(\tau,\tau)\}-s_0 v_{\nu}(y)\Big),
\]
where $(y_1,\cdots,y_n)$ is the tubular neighborhood normal coordinates of $y\in B^T_{\frac{\delta_{y^{i_0}}}{16}}(y^{i_0})$ at $y^{i_0}$, $\gamma, \rho$ are the same as in the proof of Lemma~\ref{c1l1}, and $a>0,~\beta_0>0,s_0>0$ are constants to be chosen later.

Let $\bar {H}(x_0)=\max\limits_{\overline{B^T_{\frac{\delta_{y^{i_0}}}{16}}(y^{i_0})}}\bar {H}$ for some $x_0\in\overline{ B^T_{\frac{\delta_{y^{i_0}}}{16}}(y^{i_0})}$.

\begin{claim}
Either $\bar{H}(x_0)< C$ or $x_0$ is an interior point of $B^T_{\frac{\delta_{y^{i_0}}}{16}}(y^{i_0})$ by choosing $\beta_0,s_0\gg 1$.
\label{c2c1}
\end{claim}
{\bf Proof of the Claim~\ref{c2c1}.} If not, we assume $H(x_0)\ge 1$ and $x_0\in B^T_{\frac{\delta_{y^{i_0}}}{16}}(y^{i_0})\cap\partial M$. Let $\{\bar{x}_1,\cdots \bar{x}_n\}$ be a tubular neighborhood normal coordinates at $x_0$. Then $\{\bar{x}_1,\cdots,\bar{x}_n\}$ is well-defined and smooth near $x_0$. Meanwhile, $y_n=\bar{x}_n$ near $x_0$ since they both represent the distance parameter to the boundary $\partial M$, which is to say $\frac{\partial}{\partial\nu}$ has the same definition near $x_0$. Recall that $g(\frac{\partial}{\partial \bar{x}_i},\frac{\partial}{\partial \bar{x}_j})=\delta_{ij}$ at $x_0$. W.l.o.g., we can assume $\bar {H}(x_0):=\rho e^{\beta_0  y_n}\big(v_{11}+a|\nabla (v+\gamma)|_g^2-s_0v_{\nu}\big)(x_0)$, where and in the following subindices denote the covariant derivatives w.r.t. $\frac{\partial}{\partial \bar{x}_i}$. Let
\[
H(x):=\rho e^{\beta_0 y_n}(\frac{v_{11}}{g_{11}}+a|\nabla (v+\gamma)|_g^2-s_0v_{\nu}).
\]

By $g(\frac{\partial}{\partial \bar{x}_k},\frac{\partial}{\partial\nu})=0$ near $x_0$, we know $x_0$ is a local minimum point of $H$. Moreover $\frac{\partial}{\partial \bar{x}_n}=-\frac{\partial}{\partial\nu}$ near $x_0$ implies that
\begin{eqnarray*}
(|\nabla(v+\gamma)|^2_g)_{\nu}(x_0)&=&\Big((v+\gamma)_k(v+\gamma)_lg^{kl}\Big)_{\nu}\\
&=&2(v+\gamma)_k(v+\gamma)_{k,\nu}+(v+\gamma)_k(v+\gamma)_l g^{kl}_{\nu}\\
&=&-2(v+\gamma)_k(v+\gamma)_{k,n}+(v+\gamma)_k(v+\gamma)_l g^{kl}_{\nu}\\
&=&-2(v+\gamma)_k\frac{\partial^2(v+\gamma)}{\partial \bar{x}_n\partial \bar{x}_k}+(v+\gamma)_k(v+\gamma)_l g^{kl}_{\nu}\\
&=&-2(v+\gamma)_k\frac{\partial^2(v+\gamma)}{\partial \bar{x}_k\partial \bar{x}_n}+(v+\gamma)_k(v+\gamma)_l g^{kl}_{\nu}\\
&=&-2(v+\gamma)_k(v+\gamma)_{n,k}+(v+\gamma)_k(v+\gamma)_l g^{kl}_{\nu}\\
&=&2(v+\gamma)_k(v+\gamma)_{\nu,k}+(v+\gamma)_k(v+\gamma)_l g^{kl}_{\nu}.\\
\end{eqnarray*}

Since $(v+\gamma)_{\nu}|_{\partial M}=0$ by (\ref{c1c1e0}), we have $(v+\gamma)_{\nu,k}(x_0)=0$ for $k\le n-1$, which implies that

\[
\sum\limits_{k=1}^n (v+\gamma)_k(v+\gamma)_{\nu,k}(x_0)=(v+\gamma)_n(v+\gamma)_{\nu,n}=-(v+\gamma)_{\nu}(v+\gamma)_{\nu,n}=0.
\]

Thus $(|\nabla(v+\gamma)|^2_g)_{\nu}(x_0)=(v+\gamma)_k(v+\gamma)_l g^{kl}_{\nu}$. By (\ref{c2l1e1}), $y_n=0$ at $x_0$, and $\frac{\partial y_n}{\partial \nu}=-\frac{\partial y_n}{\partial y_n}=-1$ in $B^T_{\frac{\delta_{y^{i_0}}}{16}}(y^{i_0})$,
\begin{equation}
\begin{array}{l c l}
0\le H_{\nu}(x_0)&=&\Big(v_{11,\nu}+a(v+\gamma)_k(v+\gamma)_l g^{kl}_{\nu}-v_{11}g_{11,\nu}-s_0 v_{\nu,\nu}\Big)\rho e^{\beta_0 y_n}\\
&&-\beta_0 (v_{11}+a |\nabla (v+\gamma)|_g^2-s_0 v_{\nu})\rho e^{\beta_0 y_n}\\
&=&\rho \Big(v_{11,\nu}-\beta_0 (v_{11}+a |\nabla (v+\gamma)|_g^2-s_0 v_{\nu})-v_{11}g_{11,\nu}\\
&&-s_0 v_{n,n}+a(v+\gamma)_k(v+\gamma)_l g^{kl}_{\nu}\Big)\\
\end{array}
\label{c2c1e0}
\end{equation}

We need to interchange $v_{11,\nu}$ to $\frac{\partial^2(v_{\nu})}{\partial \bar{x}_1\partial \bar{x}_1}$ in the above equation so that we can use the boundary condition. Recall $\{\bar{x}_1,\cdots,\bar{x}_{n-1}\}$ is the geodesic normal coordinates w.r.t. the metric $g|_{\partial M}$ at $x_0$. Then $\bar{\nabla}_{\frac{\partial}{\partial \bar{x}_k}}^{\frac{\partial}{\partial \bar{x}_l}}(x_0)=0$ for $1\le k,l\le n-1$, where $\bar{\nabla}$ is the covariant derivative of $\partial M$ induced by $g|_{\partial M}$. For $1\le k,l\le n-1$,
\[
\begin{array}{lcl}
\Gamma_{kl}^i(x_0)\frac{\partial}{\partial \bar{x}_i}&=&\nabla_{\frac{\partial}{\partial \bar{x}_k}}^{\frac{\partial}{\partial \bar{x}_l}}(x_0)=\bar{\nabla}_{\frac{\partial}{\partial \bar{x}_k}}^{\frac{\partial}{\partial \bar{x}_l}}(x_0)+II(\frac{\partial}{\partial \bar{x}_k},\frac{\partial}{\partial \bar{x}_l})(x_0)\frac{\partial}{\partial \nu}\\
&=&-II(\frac{\partial}{\partial \bar{x}_k},\frac{\partial}{\partial \bar{x}_l})(x_0)\frac{\partial}{\partial \bar{x}_n}.\\
\end{array}
\]

Comparing both sides of the above equation, we have, at $x_0$,
\begin{equation}
\Gamma_{kl}^i=0\quad\mbox{for}~1\le i,k,l\le n-1,\quad \Gamma_{kl}^n=-II(\frac{\partial}{\partial \bar{x}_k},\frac{\partial}{\partial \bar{x}_l}).
\label{c2c1e00}
\end{equation}

Hence at $x_0$,
\[
\begin{array}{lcl}
v_{11,\nu}&=&-v_{11,n}=-\frac{\partial}{\partial \bar{x}_n}(\frac{\partial^2 v}{\partial \bar{x}_1\partial \bar{x}_1}-\Gamma_{11}^lv_l)\\
&=&-\frac{\partial^3 v}{\partial \bar{x}_n\partial \bar{x}_1\partial \bar{x}_1}+\frac{\partial}{\partial \bar{x}_n}(\Gamma_{11}^l)v_l+\Gamma_{11}^lv_{l,n}\\
&=&-\frac{\partial^3 v}{\partial \bar{x}_n\partial \bar{x}_1\partial \bar{x}_1}+\frac{\partial}{\partial \bar{x}_n}(\Gamma_{11}^l)v_l+\Gamma_{11}^n v_{n,n}\quad\mbox{by}~~(\ref{c2c1e00})\\
&=&-\frac{\partial^3 v}{\partial \bar{x}_1\partial \bar{x}_1\partial \bar{x}_n}+\frac{\partial}{\partial \bar{x}_n}(\Gamma_{11}^l)v_l+\Gamma_{11}^n v_{n,n}\\
&=&\frac{\partial^2 (v_{\nu})}{\partial \bar{x}_1\partial \bar{x}_1}+\frac{\partial}{\partial \bar{x}_n}(\Gamma_{11}^l)v_l+\Gamma_{11}^n v_{n,n}\\
&=&\frac{\partial^2 (\psi e^v-h_g)}{\partial \bar{x}_1\partial \bar{x}_1}+\frac{\partial}{\partial \bar{x}_n}(\Gamma_{11}^l)v_l+\Gamma_{11}^n v_{n,n}\\
&=&e^v\psi \frac{\partial^2 v}{\partial \bar{x}_1\partial \bar{x}_1}+e^v\psi v_1^2+2 e^v\psi_1 v_1+e^v \frac{\partial^2 \psi}{\partial \bar{x}_1\partial \bar{x}_1}+\frac{\partial^2 h_g}{\partial \bar{x}_1\partial \bar{x}_1}\\
&&+\frac{\partial}{\partial \bar{x}_n}(\Gamma_{11}^l)v_l+\Gamma_{11}^n v_{n,n},\\
\end{array}
\]
where in the second to last equality, we used the fact that $\frac{\partial}{\partial \bar{x}_1}$ is a tangent vector field of $\partial M$ near $x_0$, so we can replace $v_{\nu}$ by $e^v\psi-h_g$.

In the following, we use $C>0$ to denote a universal constant independent of $\beta_0$. Substitute the above equation into the inequality (\ref{c2c1e0}).
\begin{eqnarray}
0&\le& H_{\nu}(x_0)=\rho \Big(e^v\psi \frac{\partial^2 v}{\partial \bar{x}_1\partial \bar{x}_1}-\beta_0 (v_{11}+a |\nabla (v+\gamma)|_g^2-s_0 v_{\nu})\nonumber{}\\
&&-v_{11}g_{11,\nu}-(s_0-\Gamma_{11}^n) v_{n,n}+C\Big)\nonumber{}\\
&=&\rho \Big(e^v\psi v_{11}+e^v\psi\Gamma_{11}^l v_l-\beta_0 (v_{11}+a |\nabla (v+\gamma)|_g^2-s_0 v_{\nu})\nonumber{}\\
&&-v_{11}g_{11,\nu}-(s_0-\Gamma_{11}^n) v_{n,n}+C\Big)\nonumber{}\\
&\le&\rho \Big((e^v\psi-\beta_0-g_{11,\nu}) (v_{11}+a |\nabla (v+\gamma)|_g^2-s_0 v_{\nu})\nonumber{}\\
&&-(s_0-\Gamma_{11}^n) v_{n,n}+C\Big).\nonumber{}\\
\label{c2c1e5}
\end{eqnarray}

Since $\frac{\partial}{\partial\nu}$ is the tangent vector of geodesic curves, we have $\nabla_{\frac{\partial}{\partial\nu}}^{\frac{\partial}{\partial\nu}}=0$ near $x_0$. In particular, we have
\[
v_{n,n}(x_0)=v_{\nu,\nu}(x_0)=v_{\nu\nu}+(\nabla_{\frac{\partial}{\partial\nu}}^{\frac{\partial}{\partial\nu}})v=v_{\nu\nu}=v_{n n}(x_0).
\]

Recall $\Gamma_{11}^n=-II(\frac{\partial}{\partial \bar{x}_1},\frac{\partial}{\partial \bar{x}_1})$ by (\ref{c2c1e00}). We can pick $s_0\gg 1$ such that $\frac{s_0}{2}$ is bigger than the largest absolute value of the principle curvatures of the second fundamental form on $\partial M$. Then we have $\frac{3s_0}{2}\ge s_0-\Gamma_{11}^n\ge\frac{s_0}{2}>0$ at $x_0$. By $\Gamma\subset\Gamma_1$, we have
\[
(1+\frac{(1-t)n}{n-2})\Delta_g v+(\frac{(2-t)n}{2}-1)|\nabla v|^2_g-\frac{(2-t)n-2}{2(n-1)(n-2)}R_g>0,
\]
which implies that $\Delta_g v(x_0)\ge -C$. W.l.o.g., we assume $v_{11}(x_0)>1$ and $v_{kk}(x_0)\le Cv_{11}(x_0)$ for $1\le k\le n-1$. Then
\[
-v_{n,n}(x_0)=-v_{n n}(x_0)\le C+\sum\limits_{k=1}^{n-1}v_{k k}(x_0)\le C v_{11}(x_0),
\]
and
\[
-(s_0-\Gamma_{11}^n) v_{n,n}(x_0)\le C(s_0-\Gamma_{11}^n)v_{11}\le \frac{3Cs_0}{2}v_{11}\le Cs_0 v_{11}.
\]

Substitute the above inequality into (\ref{c2c1e5}).

\[
\begin{array}{l c l}
0&\le& H_{\nu}(x_0)\\
&\le&\rho \Big((e^v\psi-\beta_0-g_{11,\nu}) (v_{11}+a |\nabla (v+\gamma)|_g^2-s_0 v_{\nu})+Cs_0 v_{11}+C\Big)\\
&\le&\rho \Big((e^v\psi+Cs_0-g_{11,\nu}-\beta_0) (v_{11}+a |\nabla (v+\gamma)|_g^2-s_0 v_{\nu})+C\Big)\\
&\le&\rho \Big((C-\beta_0) (v_{11}+a |\nabla (v+\gamma)|_g^2-s_0 v_{\nu})+C\Big)\\
&\le&\rho \Big(-(v_{11}+a |\nabla (v+\gamma)|_g^2-s_0 v_{\nu})+C\Big)\quad\mbox{by choosing $\beta_0>C+1$},\\
\end{array}
\]
which implies that $\Big(v_{11}+a |\nabla (v+\gamma)|_g^2-s_0 v_{\nu}\Big)(x_0)<C$ and $H(x_0)<C$. The Claim~\ref{c2c1} has been proved. $\clubsuit$.
\vskip 5 pt

Due to the above claim, we assume $x_0$ is an interior point of $B^T_{\frac{\delta_{y^{i_0}}}{16}}(y^{i_0})$. To continue the proof of the Lemma~\ref{c2l1}, we need to introduce a new coordinates near $x_0$. Let $d_0=d(x_0,\partial M)$, and let ${\cal S}_0:=\{y\in B^T_{\frac{\delta_{y^{i_0}}}{16}}(y^{i_0})|~y_n=d_0\}$. As shown in the proof of the Proposition~\ref{c2p1}, $\frac{\partial}{\partial\nu}$ is still the unit normal vector field of ${\cal S}_0$. For any $x\in B^T_{\frac{\delta_{y^{i_0}}}{16}}(y^{i_0})$ but near $x_0$ with $(y_1,\cdots,y_n)$ as its tubular neighborhood normal coordinates of $x$ at $y^{i_0}$, then $\sqrt{\sum\limits_{j=1}^{n}y_j^2}<\frac{\delta_{y^{i_0}}}{16}$. We conclude that there exists a unique $\tilde{x}\in {\cal S}_0$ such that $d(x,\tilde{x})=d(x,{\cal S}_0)$. In fact for such $x$, let $\bar{x}=(y_1,\cdots,y_{n-1},0)$. Then $\bar{x}$ is the unique point on $\partial M$ such that $d(\bar{x},x)=d(x,\partial M)=y_n$. Consider $r(t)=E(\bar{x},t)$. Then $r(t)$ is smooth and well defined for $t\in [0,\delta_{y^{i_0}})$, $r(y_n)=x$, and

\[
\sqrt{\sum\limits_{j=1}^{n-1}y_j^2+(\max\{d_0,y_n\})^2}<\frac{\delta_{y^{i_0}}}{16}
\]
as long as $x$ is close to $x_0$ enough since $x$ is an interior point of
$B^T_{\frac{\delta_{y^{i_0}}}{16}}(y^{i_0})$. Moreover for $t\in [0,\max\{d_0,y_n\}]$, the tubular neighborhood normal coordinates of $r(t)$ at $y^{i_0}$ is $(y_1,\cdots,y_{n-1},t)$, which implies that the curve $r([0,\max\{d_0,y_n\}])\subset B^T_{\frac{\delta_{y^{i_0}}}{16}}(y^{i_0})$, therefore intersects with ${\cal S}_0$ at a unique point $r(d_0)$. As shown in the proof of the Proposition~\ref{c2p1}, i.e., by (\ref{c2p1e1}) and (\ref{c2p1e2}), $d(r(t),{\cal S}_0)=|t-d_0|$ for $t\in [0,\frac{\delta_{y^{i_0}}}{2})$. In particular,
\begin{equation}
d(x,{\cal S}_0)=d(r(y_n),{\cal S}_0)=d(r(y_n),r(d_0))=|y_n-d_0|.
\label{c2n1}
\end{equation}

Next, we want to show that there exists only one point $\tilde{x}\in{\cal S}_0$ such that $d(x,\tilde{x})=d(x,{\cal S}_0)$. This is because if $(a_1,\cdots,a_n)$ is the tubular neighborhood normal coordinates of $\tilde{x}$ at $y^{i_0}$, then $\hat{x}:=(a_1,\cdots,a_{n-1},0)\in B^T_{\frac{\delta_{y^{i_0}}}{16}}(y^{i_0})\cap\partial M\subset B_{\frac{\delta_{y^{i_0}}}{8}}(y^{i_0})\cap\partial M$ by the Proposition~\ref{c2p2}. Thus $\hat{r}(t):=E(\hat{x},t)$ is smooth and well-defined for $t\in [0,\delta_{y^{i_0}})$ and $\hat{r}(a_n)=\tilde{x}$. Let $\tilde{r}(t)$ be the shortest normalized geodesic connecting $x$ with $\tilde{x}$. Then $\tilde{r}$ has $\frac{\partial}{\partial\nu}$ as its tangent vector at $\tilde{x}$. Since $\frac{\partial}{\partial\nu}$ is also the tangent vector of $\hat{r}$ at $\hat{r}(a_n)=\tilde{x}$, we know $\tilde{r}$ and $\hat{r}$ coincide. Hence $\hat{r}(y_n)=x$, which implies that $E(\hat{x},y_n)=x=E(\bar{x},y_n)$, and $\hat{x}=\bar{x}$ since $E$ is a diffeomorphism in $B_{\delta_{y^{i_0}}}(y^{i_0})$. Therefore $\tilde{x}=r(d_0)$ is uniquely determined by $x$. Clearly $r(d_0)\in{\cal S}_0$ is near $x_0$ as long as $x$ is near $x_0$. Let $\{x_1,\cdots,x_{n-1}\}$ be the geodesic normal coordinates w.r.t. the metric $g|_{{\cal S}_0}$ at $x_0$. Then $\{x_1,\cdots,x_{n-1}\}$ is smooth and well-defined near $x_0$ in ${\cal S}_0$. For any $x\in B^T_{\frac{\delta_{y^{i_0}}}{16}}(y^{i_0})$ and near $x_0$, there exists a unique $\tilde{x}\in{\cal S}_0$ such that $d(\tilde{x},x)=d(x,{\cal S}_0)$. We assume $x$ is close enough to $x_0$ such that the geodesic normal coordinates of $\tilde{x}$ w.r.t. the metric $g|_{{\cal S}_0}$ at $x_0$ is smooth and well-defined. Let $(x_1,\cdots,x_{n-1})$ be such geodesic normal coordinates of $\tilde{x}$ w.r.t. the metric $g|_{{\cal S}_0}$ at $x_0$. Define $(x_1,\cdots,x_n)$ to be the new coordinates of $x$ such that $x_n=y_n-d_0$. Then $\{x_j\}_{j=1}^n$ is smooth and well-defined for $x$ near $x_0$, and $d(x,{\cal S}_0)=|y_n-d_0|=|x_n|$ by (\ref{c2n1}). As shown in the proof of the Proposition~\ref{c2p1}, for $x$ near $x_0$ and for $1\le k\le n-1$,
\begin{equation}
\frac{\partial}{\partial x_n}=\frac{\partial}{\partial y_n}=-\frac{\partial}{\partial\nu},\qquad g(\frac{\partial}{\partial x_k},\frac{\partial}{\partial x_n})=0.
\label{c2l1jm0}
\end{equation}

Let $II_0$ denote the second fundamental form of $g$ w.r.t. $\frac{\partial}{\partial\nu}$ on ${\cal S}_0$ and let $\tilde{\nabla}$ be the Levi-Civita connection induced by $g|_{{\cal S}_0}$. Recall on ${\cal S}_0$, $\{x_j\}_{j=1}^{n-1}$ is the geodesic normal coordinates w.r.t. the metric $g|_{{\cal S}_0}$ at $x_0$. Therefore $g_{lm}(x_0):=g(\frac{\partial}{\partial x_l},\frac{\partial}{\partial x_m})(x_0)=\delta_{lm}$ for $1\le l,m\le n$,
and $\tilde{\nabla}_{\frac{\partial}{\partial x_i}}^{\frac{\partial}{\partial x_j}}(x_0)=0$ for $1\le i,j\le n-1$, which implies that
\[
\nabla_{\frac{\partial}{\partial x_i}}^{\frac{\partial}{\partial x_j}}(x_0)=\tilde{\nabla}_{\frac{\partial}{\partial x_i}}^{\frac{\partial}{\partial x_j}}+II_0(\frac{\partial}{\partial x_i},\frac{\partial}{\partial x_j})\frac{\partial}{\partial\nu}=II_0(\frac{\partial}{\partial x_i},\frac{\partial}{\partial x_j})\frac{\partial}{\partial\nu}.
\]

Thus for $1\le i,j,k\le n-1$,

\begin{equation}
\begin{array}{lcl}
\frac{\partial}{\partial x_k}g_{ij}(x_0)&=&g(\nabla_{\frac{\partial}{\partial x_k}}^{\frac{\partial}{\partial x_i}},\frac{\partial}{\partial x_j})+g(\frac{\partial}{\partial x_i},\nabla_{\frac{\partial}{\partial x_k}}^{\frac{\partial}{\partial x_j}})\\
&=&g(II_0(\frac{\partial}{\partial x_k},\frac{\partial}{\partial x_i})\frac{\partial}{\partial\nu}
,\frac{\partial}{\partial x_j})+g(\frac{\partial}{\partial x_i},II_0(\frac{\partial}{\partial x_k},\frac{\partial}{\partial x_j})\frac{\partial}{\partial\nu})\\
&=&0\quad\mbox{by}~(\ref{c2l1jm0}).\\
\end{array}
\label{c2l1jm1}
\end{equation}

Also by (\ref{c2l1jm0}), we have, for $1\le i\le n-1$, $g_{i n}=0$ near $x_0$, which implies that
\begin{equation}
\begin{array}{lcl}
\frac{\partial}{\partial x_k}g_{i n}(x_0)&=&\frac{\partial}{\partial x_k}g_{n i}(x_0)=0\quad\mbox{for}~1\le k\le n.\\
\end{array}
\label{c2l1jm2}
\end{equation}

Notice that $\frac{\partial}{\partial x_n}=-\frac{\partial}{\partial\nu}$ is a unit vector field. Therefore $g_{nn}\equiv 1$ near $x_0$, and for $1\le k\le n$,

\begin{equation}
\begin{array}{lcl}
\frac{\partial}{\partial x_k}g_{nn}(x_0)&=&0.\\
\end{array}
\label{c2l1jm3}
\end{equation}

Combine (\ref{c2l1jm1})-(\ref{c2l1jm3}). We have
\begin{equation}
\frac{\partial}{\partial x_k}g_{ij}(x_0)=0\quad\mbox{for}~ 1\le i,j\le n\quad\mbox{and}~ 1\le k\le n-1.
\label{c2l1jm4}
\end{equation}

Recall ${\cal G}=(g_{ij})$. ${\cal G}{\cal G}^{-1}=I_{n\times n}$ implies that $\frac{\partial{\cal G}^{-1}}{\partial x_k}=-{\cal G}^{-1}\frac{\partial{\cal G}}{\partial x_k} {\cal G}^{-1}$. For $1\le k\le n-1$,
\begin{equation}
\begin{array}{lcl}
\frac{\partial}{\partial x_k}g^{ij}(x_0)&=&-g^{ir}\Big(\frac{\partial}{\partial x_k}g_{rs}\Big ) g^{sj}=0.\\
\end{array}
\label{c2l1jm}
\end{equation}

In the following, subindices denote the covariant derivatives w.r.t. $\frac{\partial}{\partial x_i}$. Notice that $g_{ij}(x_0)=\delta_{ij}$. W.l.o.g., we assume
\[
\bar{H}(x_0)=\rho e^{\beta_0 d_0}\Big(v_{11}+a |\nabla (v+\gamma)|_g^2-s_0 v_{\nu}\Big),
\]
and $v_{1,1}(x_0)\gg 1$.

Let
\[
\tilde{H}=\rho e^{\beta_0 (x_n+d_0)}\Big(\frac{v_{11}}{g_{11}}+a |\nabla (v+\gamma)|_g^2-s_0 v_{\nu}\Big).
\]

By (\ref{c2l1jm0}), $x_0$ is a local maximum point of $\tilde{H}$. Near $x_0$,
\[
\begin{array}{l c l}
\tilde{H}_i&=&\rho e^{\beta_0 (x_n+d_0)}\Big(\frac{v_{11,i}}{g_{11}}-\frac{v_{11}}{g_{11}^2}g_{11,i}+2 a g^{kl}(v+\gamma)_k (v+\gamma)_{l,i}\\
&&+a g^{kl}_{,i}(v+\gamma)_k (v+\gamma)_l-s_0 v_{\nu,i}\Big)+(\frac{\rho_i}{\rho}+\delta_{n i}\beta_0)H.\\
\end{array}
\]

At $x_0$,
\begin{equation}
\begin{array}{lcl}
&&v_{11,i}-g_{11,i}v_{11}+2 a (v_k+\gamma_k) (v_{k,i}+\gamma_{k,i})
+a g^{kl}_{,i}(v_k+\gamma_k) (v_l+\gamma_l)\\
&&-s_0 v_{\nu,i}=-(\frac{\rho_i}{\rho}+\delta_{n i}\beta_0)(v_{11}+a |\nabla (v+\gamma)|_g^2-s_0 v_{\nu}),\\
\end{array}
\label{c2l1e2}
\end{equation}
and
\begin{eqnarray*}
&&\tilde{H}_{i j}(x_0)=\\
&&\rho e^{\beta_0 d_0}\Big(v_{11,i j}-g_{11,j}v_{11,i}-g_{11,i}v_{11,j}+2 g_{11,i}g_{11,j}v_{11}-g_{11,i j}v_{11}\\
&&+2 a (v_{k,i}+\gamma_{k,i}) (v_{k,j}+\gamma_{k,j})+2 a (v_k+\gamma_k) (v_{k,i j}+\gamma_{k,i j})\\
&&+2 a g^{kl}_{,j}(v_k+\gamma_k) (v_{l,i}+\gamma_{l,i})+2 a g^{kl}_{,i}(v_k+\gamma_k) (v_{l,j}+\gamma_{l,j})\\
&&+a g^{kl}_{,i j}(v_k+\gamma_k) (v_l+\gamma_l)-s_0 v_{\nu,i j}\Big)\\
&&+\beta_0 \delta_{n j}\rho e^{\beta_0 d_0}\Big(v_{11,i}-g_{11,i}v_{11}+2 a (v_k+\gamma_k) (v_{k,i}+\gamma_{k,i})\\
&&+a g^{kl}_{,i}(v_k+\gamma_k) (v_l+\gamma_l)-s_0 v_{\nu,i}  \Big)\\
&&+\rho_j e^{\beta_0 d_0}\Big(v_{11,i}-g_{11,i}v_{11}+2 a (v_k+\gamma_k) (v_{k,i}+\gamma_{k,i})\\
&&+a g^{kl}_{,i}(v_k+\gamma_k) (v_l+\gamma_l)-s_0 v_{\nu,i}  \Big)\\
&&+(\frac{\rho_{i j}\rho-\rho_i\rho_j}{\rho})(v_{11}+a |\nabla (v+\gamma)|_g^2-s_0 v_{\nu}) e^{\beta_0 d_0},\\
\end{eqnarray*}
so by (\ref{c2l1e2}),
\[
\begin{array}{l c l}
&&\rho^{-1}e^{-\beta_0 d_0}\tilde{H}_{i j}(x_0)=\\
&&\Big(v_{11,i j}-g_{11,j}v_{11,i}-g_{11,j}v_{11,i}+2 g_{11,i}g_{11,j}v_{11}-g_{11,i j}v_{11}\\
&&+2 a (v_{k,i}+\gamma_{k,i}) (v_{k,j}+\gamma_{k,j})+2 a (v_k+\gamma_k) (v_{k,i j}+\gamma_{k,i j})\\
&&+2 a g^{k l}_{,j}(v_k+\gamma_k) (v_{l,i}+\gamma_{l,i})+2 a g^{k l}_{,i}(v_k+\gamma_k) (v_{l,j}+\gamma_{l,j})\\
&&+a g^{kl}_{,i j}(v_k+\gamma_k) (v_l+\gamma_l)-s_0 v_{\nu,i j}\Big)\\
&&+\Big(\frac{\rho_{i j}\rho-2\rho_i\rho_j}{\rho^2}-\frac{\beta_0(\rho_i \delta_{n j}+\rho_j \delta_{n i})}{\rho}-\beta_0^2 \delta_{n i} \delta_{n j} \Big)(v_{11}+a |\nabla (v+\gamma)|_g^2-s_0 v_{\nu}),\\
\end{array}
\]

Recall in the proof of the Claim~\ref{c2c1}, the choice of $\beta_0$ depends on $a$. We need to prove the choice of $a$ is independent of $\beta_0$. For this reason, we let $C_1$ denote the universal constant depending only on $(M^n,g,t)$, $(f,\Gamma)$, $\phi$, $\psi$ and $\delta_{y^{i_0}}$, but independent of $a,\beta_0$, and let $C_2$ denote the universal constant depending on $(M^n,g,t)$, $(f,\Gamma)$, $\phi$, $\psi$, $\delta_{y^{i_0}}$, and $a,\beta_0$.

Notice that $g_{ij}(x_0)=\delta_{ij}$ and $\bar{L}^{ij}(x_0)=F^{ij}+\frac{1-t}{n-2}(\sum\limits_l F^{ll})\delta^{ij}$.
\begin{equation}
\begin{array}{l c l}
0&\ge&\rho^{-1}e^{-\beta_0 d_0} \bar{L}^{i j}\tilde{H}_{i j}(x_0)\\
&\ge&\bar{L}^{i j}\Big(v_{11,i j}-2 g_{11,j}v_{11,i}+2 a (v_{k,j}+\gamma_{k,j}) (v_{k,i}+\gamma_{k,i})\\
&&+2 a (v_k+\gamma_k) (v_{k,i j}+\gamma_{k,i j})-s_0 v_{\nu,i j}\Big)-C_2\rho^{-1}\sum\limits_{l}F^{ll}|v_{k,i}|,\\
\end{array}
\label{c2l1e3}
\end{equation}
where we used $|\nabla \rho|<C_1\sqrt{\rho}$, $|\nabla^2 \rho|<C_1$, and $v_{1,1}(x_0)\ge 1$.

By (\ref{c2l1e2}),
\begin{equation}
\begin{array}{l c l}
v_{11,i}(x_0)&=&-2 a (v_k+\gamma_k) (v_{k,i}+\gamma_{k,i})
-a g^{kl}_{,i}(v_k+\gamma_k) (v_l+\gamma_l)-s_0 v_{\nu,i}\\
&&+g_{11,i}v_{11}-(\frac{\rho_i}{\rho}+\delta_{n i}\beta_0)(v_{11}+a |\nabla (v+\gamma)|_g^2-s_0 v_{\nu}),\\
\end{array}
\label{c2l1a1}
\end{equation}

Substitute the above into (\ref{c2l1e3}). Since $v_{\nu,i}=-v_{n,i}$
and $v_{11}=v_{1,1}-\Gamma_{11}^lv_l$,
\begin{eqnarray}
0&\ge&\rho^{-1}e^{-\beta_0 d_0} \bar{L}^{i j}\tilde{H}_{i j}(x_0)\nonumber{}\\
&\ge&\bar{L}^{i j}\Big(v_{11,i j}+2 a (v_{k,j}+\gamma_{k,j}) (v_{k,i}+\gamma_{k,i})\nonumber{}\\
&&+2 a (v_k+\gamma_k) (v_{k,i j}+\gamma_{k,i j})-s_0 v_{\nu,i j}\Big)-C_2\rho^{-1}\sum\limits_{l}F^{ll}|v_{k,i}|\nonumber{}\\
&\ge&\bar{L}^{i j}\Big(v_{11,i j}+2 a (1+\psi\beta e^v)^2 v_{k,j}v_{k,i}+2 a (1+\psi\beta e^v)(v_k+\gamma_k) v_{k,i j}\nonumber{}\\
&&-s_0 v_{\nu,i j}\Big)-C_2\rho^{-1}\sum\limits_{l}F^{ll}|v_{k,i}|\nonumber{}\\
\label{c2l1e03}
\end{eqnarray}

At $x_0$,
\[
\begin{array}{l c l}
v_{i j,l}&=&\frac{\partial}{\partial x_l}(\nabla^2 v (\frac{\partial}{\partial x_i},\frac{\partial}{\partial x_j}))=\frac{\partial}{\partial x_l}(\frac{\partial^2 v}{\partial x_j\partial x_i}-\Gamma_{j i}^k v_k)\\
&=&\frac{\partial^3 v}{\partial x_l\partial x_j\partial x_i}-\Gamma_{j i}^k v_{k,l}-\frac{\partial(\Gamma_{j i}^k)}{\partial x_l}v_k,\\
\end{array}
\]
so
\begin{equation}
\begin{array}{l c l}
v_{l,i j}&=&(\nabla^2 v_l)(\frac{\partial}{\partial x_i},\frac{\partial}{\partial x_j})=(\nabla_{\frac{\partial}{\partial x_j}}\nabla_{\frac{\partial}{\partial x_i}}-\Gamma_{j i}^k\frac{\partial}{\partial x_k})(v_l)\\
&=&\frac{\partial^3 v}{\partial x_j\partial x_i\partial x_l}-\Gamma_{j i}^k v_{l,k}\\
&=&v_{i j,l}+\frac{\partial(\Gamma_{j i}^k)}{\partial x_l}v_k,\\
\end{array}
\label{c2l1e4}
\end{equation}
and
\[
\begin{array}{l c l}
v_{i j,11}&=&(\nabla^2 v_{i j})(\frac{\partial}{\partial x_1},\frac{\partial}{\partial x_1})=(\nabla_{\frac{\partial}{\partial x_1}}\nabla_{\frac{\partial}{\partial x_1}}-\Gamma_{11}^l\frac{\partial}{\partial x_l}) v_{i j}\\
&=&\frac{\partial^2}{\partial x_1\partial x_1}(\frac{\partial^2 v}{\partial x_j\partial x_i}-\Gamma_{j i}^l v_l)-\Gamma_{11}^l v_{i j,l}\\
&=&\frac{\partial^4 v}{\partial x_1\partial x_1\partial x_j\partial x_i}-\Gamma_{j i}^l v_{l,1,1}-2\frac{\partial (\Gamma_{j i}^l)}{\partial x_1} v_{l, 1}-\frac{\partial^2(\Gamma_{i j}^l)}{\partial x_1\partial x_1}v_l-\Gamma_{11}^l v_{i j,l}\\
&=&\frac{\partial^4 v}{\partial x_1\partial x_1\partial x_j\partial x_i}-\Gamma_{j i}^l (v_{1 1}+\Gamma_{1 1}^k v_k)_l-2\frac{\partial (\Gamma_{j i}^l)}{\partial x_1} v_{l, 1}\\
&&-\frac{\partial^2(\Gamma_{i j}^l)}{\partial x_1\partial x_1}v_l-\Gamma_{11}^l( v_{i ,j}-\Gamma_{j i}^k v_k)_l\\
&=&\frac{\partial^4 v}{\partial x_1\partial x_1\partial x_j\partial x_i}-\Gamma_{j i}^l v_{1 1,l}-\Gamma_{j i}^l\Gamma_{1 1}^k v_{k,l}-\frac{\partial \Gamma_{1 1}^k}{\partial x_l}\Gamma_{j i}^l v_k\\
&&-2\frac{\partial (\Gamma_{j i}^l)}{\partial x_1} v_{l, 1}-\frac{\partial^2(\Gamma_{i j}^l)}{\partial x_1\partial x_1}v_l-\Gamma_{11}^l v_{i ,j,l}-\Gamma_{1 1}^l\Gamma_{j i}^k v_{k,l}-\frac{\partial\Gamma_{j i}^k}{\partial x_l}\Gamma_{1 1}^l v_k,\\
\end{array}
\]
therefore
\begin{equation}
\begin{array}{lcl}
v_{11,i j}&=&(\frac{\partial^2}{\partial x_j\partial
x_i}-\Gamma_{j i}^k\frac{\partial}{\partial
x_k})(v_{11})=\frac{\partial^2(v_{11})}{\partial x_j\partial
x_i}-\Gamma_{j i}^k v_{11,k}\nonumber{}\\
&=&\frac{\partial^2}{\partial x_j\partial x_i}(\frac{\partial^2 v}{\partial x_1\partial x_1}-\Gamma_{11}^l v_l)
-\Gamma_{j i}^k v_{11,k}\nonumber{}\\
&=&\frac{\partial^4 v}{\partial x_j\partial x_i\partial x_1\partial
x_1}-\Gamma_{11}^l v_{l,i,j}-\frac{\partial (\Gamma_{11}^l)}{\partial x_j} v_{l, i}
-\frac{\partial (\Gamma_{11}^l)}{\partial x_i} v_{l, j}\nonumber{}\\
&&-\frac{\partial^2(\Gamma_{11}^l)}{\partial x_j\partial x_i}v_l-\Gamma_{j i}^k v_{11,k}\nonumber{}\\
&=&v_{i j,11}+\Gamma_{j i}^l\Gamma_{1 1}^k v_{k,l}+\frac{\partial \Gamma_{1 1}^k}{\partial x_l}\Gamma_{j i}^l v_k+2\frac{\partial (\Gamma_{j i}^l)}{\partial x_1} v_{l, 1}+\frac{\partial^2(\Gamma_{i j}^l)}{\partial x_1\partial x_1}v_l\nonumber{}\\
&&+\Gamma_{1 1}^l\Gamma_{j i}^k v_{k,l}+\frac{\partial\Gamma_{j i}^k}{\partial x_l}\Gamma_{1 1}^l v_k-\frac{\partial (\Gamma_{11}^l)}{\partial x_j} v_{l, i}-\frac{\partial (\Gamma_{11}^l)}{\partial x_i} v_{l, j}-\frac{\partial^2(\Gamma_{11}^l)}{\partial x_j\partial x_i}v_l\nonumber{}\\
\end{array}
\label{c2l1e5}
\end{equation}

Substitute (\ref{c2l1e4}) and (\ref{c2l1e5}) into (\ref{c2l1e3}). At $x_0$,

\begin{equation}
\begin{array}{lcl}
0&\ge&\rho^{-1}e^{-\beta_0 d_0} \bar{L}^{i j}\tilde{H}_{i j}(x_0)\\
&\ge&\bar{L}^{i j}\Big(v_{i j,11}+2 a (1+\psi\beta e^v)^2 v_{k,j}v_{k,i}+2 a (1+\psi\beta e^v)(v_k+\gamma_k) v_{i j,k}\\
&&-s_0 v_{ij,\nu}\Big)-C_2\rho^{-1}\sum\limits_{l}F^{ll}|v_{k,i}|\\
\end{array}
\label{c2l1e6}
\end{equation}

Differentiate the equation $F(\bar {W}_{ij}g^{jr})=\phi e^{2v}$ along the $x_l-$th direction.
\[
\begin{array}{lcl}
(\phi e^{2v})_l&=&F^{ir}\Big(g^{jr}\bar{W}_{ij,l}+\bar {W}_{ij}g^{jr}_{,l}\Big)\\
&=&F^{ir}g^{jr}(v_{ij,l}+\frac{1-t}{n-2}(\Delta_g v)_lg_{ij})+F^{ir}g^{jr}\Big(\frac{1-t}{n-2}(\Delta_g v)g_{ij,l}\\
&&+(2-t)v_mv_{k,l}g^{km}g_{ij}+\frac{2-t}{2}v_kv_mg^{km}_{,l}g_{ij}+\frac{2-t}{2}|\nabla v|_g^2g_{ij,l}\\
&&-v_{i,l}v_j-v_{j,l}v_i-(A_g^t)_{ij,l}\Big)+F^{ir}g^{jr}_l\bar{W}_{ij},\\
\end{array}
\]
which implies that, at $x_0$, by $g_{ij}=\delta_{ij}$
\begin{equation}
|F^{ij}v_{ij,l}+\frac{1-t}{n-2}(\Delta_g v)_l\sum\limits_i F^{ii}|\le C_1\sum\limits_{i,j,r}F^{rr}|v_{i,j}|,
\label{c2l1e06}
\end{equation}
where we used $|\Delta_g v|(x_0)=|\sum\limits_{k}v_{kk}|\le C_1 \sum\limits_{k,m}|v_{k,m}|$ and
\[
\begin{array}{lcl}
|\bar{W}_{ij}(x_0)|&=&|v_{ij}+\frac{1-t}{n-2}(\Delta_g v)g_{ij}+\frac{2-t}{2}|\nabla v|_g^2g_{ij}-v_iv_j-(A_g^t)_{ij}|\\
&\le& C_1\sum\limits_{k,m}|v_{k,m}|.\\
\end{array}
\]

Recall the Laplace-Beltrami operator $\Delta_g=\frac{1}{\sqrt{|g|}}\frac{\partial}{\partial x_k}(\sqrt{|g|}g^{km}\frac{\partial}{\partial x_m})$.
\[
\begin{array}{lcl}
(\Delta_g v)_l&=&\Big(\frac{1}{\sqrt{|g|}}\frac{\partial}{\partial x_k}(\sqrt{|g|}g^{km}v_m) \Big)_l\\
&=&\Big(\frac{1}{\sqrt{|g|}}(\sqrt{|g|}g^{km}) v_{m,k}+\frac{1}{\sqrt{|g|}}(\sqrt{|g|}g^{km})_k v_m \Big)_l\\
&=&g^{km}v_{m,k,l}+g^{km}_{,l}v_{m,k}+\Big(\frac{1}{\sqrt{|g|}}(\sqrt{|g|}g^{km})_k\Big)_l v_m+\frac{1}{\sqrt{|g|}}(\sqrt{|g|}g^{km})_k v_{m,l}\\
&=&g^{km}(v_{mk}+\Gamma_{km}^sv_s)_l+g^{km}_{,l}v_{m,k}+\Big(\frac{1}{\sqrt{|g|}}(\sqrt{|g|}g^{km})_k\Big)_l v_m\\
&&+\frac{1}{\sqrt{|g|}}(\sqrt{|g|}g^{km})_k v_{m,l}\\
&=&g^{km}v_{mk,l}+g^{km}\Gamma_{km}^sv_{s,l}+g^{km}\frac{\partial\Gamma_{km}^s}{\partial x_l}v_s+g^{km}_{,l}v_{m,k}\\
&&+\Big(\frac{1}{\sqrt{|g|}}(\sqrt{|g|}g^{km})_k\Big)_l v_m+\frac{1}{\sqrt{|g|}}(\sqrt{|g|}g^{km})_k v_{m,l}.\\
\end{array}
\]

Substitute the above identity into (\ref{c2l1e06}). At $x_0$,

\begin{equation}
\begin{array}{lcl}
C_1\sum\limits_{i,j,r}F^{rr}|v_{i,j}|&\ge&
|F^{ij}v_{ij,l}+\frac{1-t}{n-2}g^{km}v_{mk,l}\sum\limits_i F^{ii}|\\
&=&|F^{ij}v_{ij,l}+\frac{1-t}{n-2}\sum\limits_{i,k}F^{ii} v_{kk,l}|\\
&=&|\bar{L}^{ij}v_{ij,l}|\\
\end{array}
\label{c2l1e7}
\end{equation}

Differentiate the equation $F(\bar {W}_{il}g^{lj})=\phi e^{2v}$ along the $x_1-$th direction twice and evaluate it at $x_0$.
\[
\begin{array}{lcl}
F^{ij}(\bar {W}_{il}g^{lj})_{11}+F^{ij,rs}(\bar {W}_{il}g^{lj})_1(\bar {W}_{rk}g^{ks})_1&=&e^{2v}(\phi_{11}+4\phi v_1^2+4\phi_1 v_1+2\phi v_{11})\\
&\ge&e^{2 v}(\phi_{11}+4\phi v_1^2+4\phi_1 v_1)\\
\end{array}
\]
since we have already assumed $v_{1,1}(x_0)\gg 1$.

By the concavity of $f$ in $\Gamma$, we have $F^{ij,rs}(\bar {W}_{il}g^{lj})_1(\bar {W}_{rk}g^{ks})_1<0$, and
\[
\begin{array}{lcl}
-C_1&\le& F^{ij}(\bar {W}_{il}g^{lj})_{11}(x_0)\\
&=&F^{ij}\Big(\bar {W}_{ij,11}+2\bar {W}_{il,1}g^{l j}_{,1}+\bar {W}_{ij}g^{lj}_{,11} \Big)\\
&=&F^{ij}\Big(\bar {W}_{ij,11}+\bar {W}_{ij}g^{lj}_{,11} \Big)\quad\mbox{by}~~(\ref{c2l1jm})\\
&\le&F^{ij}\bar {W}_{ij,11}+C_1\sum\limits_{i,j,k}F^{ii}|v_{j,k}|\\
&=&F^{ij}\Big( v_{ij,11}+\frac{1-t}{n-2}(\Delta_g v)_{11}\delta_{ij}
+\frac{2(1-t)}{n-2}(\Delta_g v)_1 g_{ij,1}+\frac{1-t}{n-2}(\Delta_g v)g_{ij,11}\\
&&+(2-t)v_{k,1}^2\delta_{ij}+(2-t)v_kv_{k,11}\delta_{ij}+2(2-t)v_kv_{l,1}g^{kl}_{,1}\delta_{ij}\\
&&+\frac{2-t}{2}v_kv_lg^{kl}_{,11}\delta_{ij}+\frac{2-t}{2}|\nabla v|^2_g g_{ij,11}-2 v_{i,1}v_{j,1}-2 v_i v_{j,11}-(A_g^t)_{i j,11}\Big)\\
&&+C_1\sum\limits_{i,j,k}F^{ii}|v_{j,k}|\\
&=&F^{ij}\Big( v_{ij,11}+\frac{1-t}{n-2}(\Delta_g v)_{11}\delta_{ij}
+\frac{1-t}{n-2}(\Delta_g v)g_{ij,11}+(2-t)v_{k,1}^2\delta_{ij}\\
&&+(2-t)v_kv_{k,11}\delta_{ij}+\frac{2-t}{2}v_kv_lg^{kl}_{,11}\delta_{ij}+\frac{2-t}{2}|\nabla v|^2_g g_{ij,11}-2v_{i,1}v_{j,1}\\
&&-2 v_i v_{j,11}-(A_g^t)_{i j,11}\Big)+C_1\sum\limits_{i,j,k}F^{ii}|v_{j,k}|,~~~~\mbox{by}~(\ref{c2l1jm4})~~\mbox{and}~~(\ref{c2l1jm}).\\
\end{array}
\]

Thus at $x_0$,
\[
\begin{array}{l c l}
-C_1\sum\limits_{i,j,k}F^{ii}|v_{j,k}|&\le& F^{ij}\Big( v_{ij,11}+\frac{1-t}{n-2}(\Delta_g v)_{11}\delta_{ij}
+(2-t)v_kv_{k,11}\delta_{ij}\\
&&-2 v_iv_{j,11}+(2-t)v^2_{k,1}\delta_{ij}-2v_{i,1}v_{j,1}\Big)\\
&\le&F^{ij}\Big( v_{ij,11}+\frac{1-t}{n-2}(\Delta_g v)_{11}\delta_{ij}
+(2-t)v_kv_{11,k}\delta_{ij}\\
&&-2 v_iv_{11,j}+(2-t)v^2_{k,1}\delta_{ij}-2v_{i,1}v_{j,1}\Big)\\
&&+C_1\sum\limits_{i,j,k}F^{ii}|v_{j,k}|\quad\mbox{by}~~(\ref{c2l1e4})\\
&\le&F^{ij}\Big( v_{ij,11}+\frac{1-t}{n-2}(\Delta_g
v)_{11}\delta_{ij}+(2-t)v^2_{k,1}\delta_{ij}-2v_{i,1}v_{j,1}\Big)\\
&&+C_2\frac{1}{\sqrt{\rho}}\sum\limits_{i,j,k}F^{ii}|v_{j,k}|\quad\mbox{by}~~(\ref{c2l1a1})\\
\end{array}
\]
i.e., at $x_0$
\begin{equation}
\begin{array}{lcl}
&&F^{ij}\Big( v_{ij,11}+\frac{1-t}{n-2}(\Delta_g
v)_{11}\delta_{ij}\Big)\ge-(2-t)\sum\limits_{k,i}F^{ii}v^2_{k,1}+2F^{ij}v_{i,1}v_{j,1}\nonumber{}\\
&&-C_2\frac{1}{\sqrt{\rho}}\sum\limits_{i,j,k}F^{ii}|v_{j,k}|\ge-C_1\sum\limits_{i,j,k}F^{ii}v^2_{j,k}-C_2\frac{1}{\sqrt{\rho}}\sum\limits_{i,j,k}F^{ii}|v_{j,k}|\nonumber{}\\
\label{c2l1a2}
\end{array}
\end{equation}

For the term $(\Delta_g v)_{11}$ in the above inequality, we need to replace it by $\sum\limits_k v_{kk,11}$. For this reason, recall $\Delta_g=\frac{1}{\sqrt{|g|}}\frac{\partial}{\partial x_k}(\sqrt{|g|}g^{km}\frac{\partial}{\partial x_m})$,
\begin{eqnarray*}
(\Delta_g v)_{11}(x_0)&=&\Big(\frac{1}{\sqrt{|g|}}\frac{\partial}{\partial x_k}(\sqrt{|g|}g^{km}v_m)\Big)_{11}\\
&=&\Big(g^{km} v_{m,k}+\frac{1}{\sqrt{|g|}}(\sqrt{|g|}g^{km})_k v_m \Big)_{11}\\
&=&\Big(g^{km} (v_{m k}+\Gamma_{km}^lv_l)+\frac{1}{\sqrt{|g|}}(\sqrt{|g|}g^{km})_k v_m \Big)_{11}\\
&=&\Big(g^{km} v_{m k}+g^{km}\Gamma_{km}^l v_l+\frac{1}{\sqrt{|g|}}(\sqrt{|g|}g^{km})_k v_m \Big)_{11}\\
&=&g^{km} v_{m k,11}+2g^{km}_{,1} v_{m k,1}+g^{km}_{,11} v_{m k}+g^{km}\Gamma_{km}^l v_{l,11}\\
&&+2 (g^{km}\Gamma_{km}^l)_1 v_{l,1}+(g^{km}\Gamma_{km}^l)_{11} v_l+\frac{1}{\sqrt{|g|}}(\sqrt{|g|}g^{km})_k v_{m,11}\\
&&+2(\frac{1}{\sqrt{|g|}}(\sqrt{|g|}g^{km})_k)_1 v_{m,1}+(\frac{1}{\sqrt{|g|}}(\sqrt{|g|}g^{km})_k)_{11} v_m \\
&=&v_{k k,11}+g^{km}_{,11} v_{m k}+\Gamma_{kk}^l v_{l,11}+2 (g^{km}\Gamma_{km}^l)_1 v_{l,1}\\
&&+(g^{km}\Gamma_{km}^l)_{11} v_l+(\sqrt{|g|}g^{km})_k v_{m,11}+2(\frac{1}{\sqrt{|g|}}(\sqrt{|g|}g^{km})_k)_1 v_{m,1}\\
&&+(\frac{1}{\sqrt{|g|}}(\sqrt{|g|}g^{km})_k)_{11} v_m \quad\mbox{by}~~(\ref{c2l1jm}).\\
\end{eqnarray*}

Plug the above equation into (\ref{c2l1a2}). At $x_0$,
\begin{eqnarray*}
&&-C_1\sum\limits_{i,j,k}F^{ii}v^2_{j,k}-C_2\frac{1}{\sqrt{\rho}}\sum\limits_{i,j,k}F^{ii}|v_{j,k}|\le\\
&&F^{ij}\Big( v_{ij,11}+\frac{1-t}{n-2}v_{k k,11}\delta_{ij}+\frac{1-t}{n-2}\Gamma_{kk}^l v_{l,11}\delta_{ij}\\
&&+\frac{1-t}{n-2}(\sqrt{|g|}g^{km})_k v_{m,11}\delta_{ij}\Big)+C_1\sum\limits_{i,j,k}F^{ii}|v_{j,k}|\\
&\le& F^{ij}\Big( v_{ij,11}+\frac{1-t}{n-2}v_{k k,11}\delta_{ij}+\frac{1-t}{n-2}\Gamma_{kk}^l v_{11,l}\delta_{ij}\\
&&+\frac{1-t}{n-2}(\sqrt{|g|}g^{km})_k v_{11,m}\delta_{ij}\Big)+C_1\sum\limits_{i,j,k}F^{ii}|v_{j,k}|\quad\mbox{by}~~(\ref{c2l1e4})\\
&\le& F^{ij}\Big( v_{ij,11}+\frac{1-t}{n-2}v_{k k,11}\delta_{ij}\Big)+C_2\frac{1}{\sqrt{\rho}}\sum\limits_{i,j,k}F^{ii}|v_{j,k}|\quad\mbox{by}~~(\ref{c2l1a1})\\
&=&\bar{L}^{ij}v_{ij,11}+C_2\frac{1}{\sqrt{\rho}}\sum\limits_{i,j,k}F^{ii}|v_{j,k}|,\\
\end{eqnarray*}
that is,
\begin{equation}
\begin{array}{l c l}
\bar{L}^{ij}v_{ij,11}(x_0)&\ge&-C_1\sum\limits_{i,j,k}F^{ii}v^2_{j,k}-C_2\frac{1}{\sqrt{\rho}}\sum\limits_{i,j,k}F^{ii}|v_{j,k}|\\
\end{array}
\label{c2l1e8}
\end{equation}

Substitute (\ref{c2l1e7}) and (\ref{c2l1e8}) into (\ref{c2l1e6}).
Notice that $v_{ij,\nu}=-v_{ij,n}$ and $1+\psi\beta e^v\in [\frac
12,1]$. At $x_0$,
\[
\begin{array}{l c l}
0&\ge& 2a(1+\psi\beta e^v)^2\bar{L}^{i j}v_{k,j} v_{k,i}-C_2\rho^{-1}\sum\limits_{l,k,i}F^{ll}|v_{k,i}|-C_1\sum\limits_{l,k,i} F^{ll}v_{k,i}^2\\
&\ge&\frac{a}{2}\bar{L}^{i j}v_{k,j} v_{k,i}-C_2\rho^{-1}\sum\limits_{l,k,i}F^{ll}|v_{k,i}|-C_1\sum\limits_{l,k,i} F^{ll}v_{k,i}^2\\
&\ge& \frac{a (1-t)}{2(n-2)}\sum\limits_l F^{ll} v^2_{k,i}-C_2\rho^{-1}\sum\limits_{l,k,i}F^{ll}|v_{k,i}|-C_1\sum\limits_{l,k,i} F^{ll}v_{k,i}^2\\
&\ge&\sum\limits_{l,k,i} F^{ll}
v^2_{k,i}-C_2\rho^{-1}\sum\limits_{l,k,i}F^{ll}|v_{k,i}|\quad\mbox{by
taking}~~
a>\frac{2(n-2)(C_1+1)}{(1-t)}.\\
\end{array}
\]

Multiply the above inequality by $\rho^2$. At $x_0$,
\[
\begin{array}{l c l}
0&\ge&  \sum\limits_{l,k,i} F^{ll} \Big( (\rho v_{k,i})^2-C_2\rho|v_{k,i}|\Big),\\
\end{array}
\]
which implies that $(\rho |v_{1,1}|)(x_0)<C_2$, therefore
$H(x_0)<C_2$. Lemma~\ref{c2l1} has been established. $\clubsuit$.

\vskip 5 pt

\begin{remark}
As a consequence of the Lemma~\ref{c2l1}, $\forall$ $y\in
B^T_{\frac{\delta_{y^{i_0}}}{32}}(y^{i_0})$, let $(e_1,\cdots,e_n)$ be an
orthonormal basis of $T_yM$ with $e_n=\frac{\partial}{\partial\nu}$,
and let subindices denote the covariant derivatives w.r.t. $e_j$. By
$\Gamma\subset\Gamma_1$, we have $\Delta_g v(y)>-C$, which implies
that $v_{\nu\nu}(y)>-C$,  and for $1\le k\le n-1$,
\[
\begin{array}{lcl}
v_{kk}(y)=\Delta_g v-\sum\limits_{l\neq k,n}v_{ll}-v_{nn}
\ge-C-v_{nn}=-C-v_{\nu\nu}.\\
\end{array}
\]

If $v_{\nu\nu}(y)\ge 0$, then $v_{\nu\nu}+C\ge C>v_{kk}(y)>-C-v_{\nu\nu}$ implies that $|v_{kk}(y)|\le C+v_{\nu\nu}(y)$ for $1\le k\le n-1$. If $v_{\nu\nu}(y)<0$, then $C>v_{kk}(y)>-C-v_{\nu\nu}>-C$ implies that $|v_{kk}(y)|\le C\le C+v_{\nu\nu}(y)$ for $1\le k\le n-1$ since $v_{\nu\nu}(y)\ge -C$. Hence, for any two vectors
$X,Y\in T_yM$ with
$g(X,\frac{\partial}{\partial\nu})=g(Y,\frac{\partial}{\partial\nu})=0$,
\begin{equation}
\begin{array}{lcl}
|\nabla^2_g v(X,Y)|(y)&=&|\frac 12(\nabla^2_g v(X+Y,X+Y)-\nabla^2_g v(X,X)-\nabla^2_g v(Y,Y))|\\
&\le& \frac 12(|\nabla^2_g v(X+Y,X+Y)|+|\nabla^2_g v(X,X)|+|\nabla^2_g v(Y,Y)|\\
&\le& \frac 12(|X+Y|_g^2+|X|_g^2+|Y|_g^2)(v_{\nu\nu}+C)\\
&\le& \frac 32(|X|_g^2+|Y|_g^2)(v_{\nu\nu}+C)\le 2(|X|_g^2+|Y|_g^2)(v_{\nu\nu}+C)\\
&\le& (|X|_g^2+|Y|_g^2)(2 v_{\nu\nu}+C).\\
\end{array}
\label{c2l1r1e1}
\end{equation}
\label{c2l1r1}
\end{remark}
\begin{lemma}
Under the same assumptions as in Theorem~\ref{generalresult}, for $t<1$, let $v$ be a $C^4$ solution of the equation (\ref{generalproblem}). Then there exists a universal constant $C>0$ depending only on $(M^n,g,t)$, $(f,\Gamma)$, $\phi$, $\psi$, $\delta_{y^{i_0}}$ such that in $B^T_{\frac{\delta_{y^{i_0}}}{64}}(y^{i_0})\cap\partial M$,
\[
v_{\nu\nu}<C.
\]
\label{c2l2}
\end{lemma}
{\bf Proof of the Lemma~\ref{c2l2}.}  Let $\{y_1,\cdots,y_n\}$ be
the tubular neighborhood normal coordinates of $y\in
B^T_{\frac{\delta_{y^{i_0}}}{32}}(y^{i_0})$ at $y^{i_0}$. Let
$\{e_1,\cdots, e_n\}$ be a smooth orthonormal frame of $TM$ in
$B^T_{\frac{\delta_{y^{i_0}}}{32}}(y^{i_0})$ with
$e_n=\frac{\partial}{\partial\nu}$. In fact, we can obtain such
frame by moving an orthonormal basis of $T_{y^{i_0}}(\partial M)$
parallelly along the geodesic of $(\partial M,g|_{\partial M})$ to
get an orthonormal frame of $T(\partial M)$ in
$B^T_{\frac{\delta_{y^{i_0}}}{32}}(y^{i_0})$, then moving such frame
parallelly along the geodesic
$r(t)=E(\frac{\partial}{\partial\nu},t)$. In this way, we can get
smooth orthonormal vector fields $\{e_j\}_{j=1}^{n-1}$ in
$B^T_{\frac{\delta_{y^{i_0}}}{32}}(y^{i_0})$ with
$g(e_j,\frac{\partial}{\partial\nu})=0$, and $\{e_j\}_{j=1}^n$ with
$e_n=\frac{\partial}{\partial\nu}$ will be an orthonormal frame of
$TM$ in $B^T_{\frac{\delta_{y^{i_0}}}{32}}(y^{i_0})$.

Observe $\frac{\partial}{\partial \nu}$ is the unit tangent vector of the geodesic. We have
\begin{equation}
\nabla_{\frac{\partial}{\partial\nu}}^{\frac{\partial}{\partial\nu}}=0\qquad\mbox{in}~B^T_{\frac{\delta_{y^{i_0}}}{32}}(y^{i_0}).
\label{c2tn1}
\end{equation}

In the following, subindices denote the covariant derivatives w.r.t.
$\{e_1,\cdots,e_n\}$. Differentiate the equation $F(\bar
{W}_{ij})=\phi e^{2v}$ along the normal direction $e_n$,
\begin{equation}
F^{ij}\Big( v_{ij,\nu}+\frac{1-t}{n-2}(\Delta v)_{\nu}\delta_{ij}
+(2-t)v_k v_{k,\nu}\delta_{ij}-2v_i v_{j,\nu}-(A^t_g)_{ij,\nu}\Big)=e^{2v}(\phi_{\nu}+2\phi v_{\nu})
\label{c2l2e1}
\end{equation}

We need to interchange $v_{ij,\nu}$ to $v_{\nu,ij}$. For this reason, let $e_i=a_i^j\frac{\partial}{\partial y_j}$. Then $a_i^j\in C^{\infty}(B^T_{\frac{\delta_{y^{i_0}}}{32}}(y^{i_0}))$, and

\begin{equation}
g(e_i,e_j)=\delta_{ij}\Longleftrightarrow a_i^k g_{kl}a_j^l=\delta_{i j}.
\label{c2l2e2}
\end{equation}

Notice $e_n=\frac{\partial}{\partial\nu}=-\frac{\partial}{\partial y_n}$ and $g(e_i,e_n)=\delta_{i n}$. We have
\begin{equation}
 a_n^k=-\delta_{n}^k,\qquad a_i^n=-\delta_i^n.
\label{c2l2e3}
 \end{equation}

In $B^T_{\frac{\delta_{y^{i_0}}}{32}}(y^{i_0})$,
\begin{equation}
\begin{array}{l c l}
v_{i,\nu}&=&-\frac{\partial}{\partial y_n}(a_i^r\frac{\partial v}{\partial y_r})
=-a_i^r\frac{\partial^2 v}{\partial y_n\partial y_r}-\frac{\partial a_i^r}{\partial y_n}\frac{\partial v}{\partial y_r}.\\
\end{array}
\label{c2l2e4}
\end{equation}
\[
\begin{array}{l c l}
v_{\nu,ij}&=&(\nabla^2 v_{\nu})(e_i,e_j)=a_i^r a_j^s(\nabla^2 v_{\nu})(\frac{\partial}{\partial y_r},\frac{\partial}{\partial y_s})\\
&=&a_i^r a_j^s(\frac{\partial^2}{\partial y_s\partial y_r}-\Gamma_{sr}^l\frac{\partial}{\partial y_l})(-\frac{\partial v}{\partial y_n})\\
&=&-a_i^r a_j^s \frac{\partial^3 v}{\partial y_s\partial y_r\partial y_n}+a_i^ra_j^s\Gamma_{sr}^l\frac{\partial^2 v}{\partial y_l\partial y_n},\\
\end{array}
\]
so

\begin{equation}
\begin{array}{lcl}
v_{ij,\nu}&=&\frac{\partial}{\partial\nu}\Big(\nabla^2 v(e_i,e_j)\Big)=-\frac{\partial}{\partial y_n}\Big(a_i^ra_j^s\nabla^2 v (\frac{\partial}{\partial y_r},\frac{\partial}{\partial y_s})\Big)\\
&=&-a_i^r a_j^s\frac{\partial}{\partial y_n}\Big(\nabla^2 v
(\frac{\partial}{\partial y_r},\frac{\partial}{\partial
y_s})\Big)-\frac{\partial (a_i^ra_j^s)}{\partial y_n}\Big(\nabla^2 v
(\frac{\partial}{\partial y_r},\frac{\partial}{\partial y_s})\Big)\\
&=&-a_i^r a_j^s\frac{\partial}{\partial y_n}\Big(\frac{\partial^2
v}{\partial y^s\partial y^r}-\Gamma_{sr}^l\frac{\partial v}{\partial
y_l}\Big)-\frac{\partial (a_i^ra_j^s)}{\partial
y_n}\Big(\frac{\partial^2 v}{\partial y^s\partial
y^r}-\Gamma_{sr}^l\frac{\partial v}{\partial y_l}\Big)\\
&=&-a_i^r a_j^s\frac{\partial^3 v}{\partial y_n\partial y_s\partial y_r}+a_i^r a_j^s\Gamma_{sr}^l\frac{\partial^2 v}{\partial y_n\partial y_l}+a_i^r a_j^s\frac{\partial v}{\partial y_l}\frac{\partial\Gamma_{sr}^l}{\partial y_n}\\
&&-\frac{\partial (a_i^r a_j^s)}{\partial y_n}\frac{\partial^2 v}{\partial y_s\partial y_r}+\frac{\partial (a_i^r a_j^s)}{\partial y_n}\Gamma_{sr}^l\frac{\partial v}{\partial y_l}\\
&=&v_{\nu,ij}+a_i^r a_j^s\frac{\partial v}{\partial
y_l}\frac{\partial\Gamma_{sr}^l}{\partial y_n} -\frac{\partial
(a_i^r a_j^s)}{\partial y_n}\frac{\partial^2 v}{\partial y_s\partial
y_r}
+\frac{\partial (a_i^r a_j^s)}{\partial y_n}\Gamma_{sr}^l\frac{\partial v}{\partial y_l}\\
&:=&v_{\nu,ij}+\Omega_{ij}^{rs}\frac{\partial^2 v}{\partial y_r\partial y_s}+\Theta_{ij}^l\frac{\partial v}{\partial y_l},
\end{array}
\label{c2l2e5}
\end{equation}
where
\[
\begin{array}{l c l}
\Omega_{ij}^{rs}=-\frac{\partial (a_i^r a_j^s)}{\partial
y_n},\quad\mbox{ and}\quad \Theta_{ij}^l=a_i^r
a_j^s\frac{\partial\Gamma_{sr}^l}{\partial y_n}
+\frac{\partial (a_i^r a_j^s)}{\partial y_n}\Gamma_{sr}^l,&&\\
\end{array}
\]
depend only on $(M^n,g)$, and are smooth and bounded in $B^T_{\frac{\delta_{y^{i_0}}}{32}}(y^{i_0})$.

In particular,

\begin{equation}
(\Delta v)_{\nu}=\Delta (v_{\nu})+\Omega_{kk}^{rs}\frac{\partial^2 v}{\partial y_r\partial y_s}+\Theta_{kk}^l\frac{\partial v}{\partial y_l}.
\label{c2l2e6}
\end{equation}

Substitute (\ref{c2l2e4}), (\ref{c2l2e5}), and (\ref{c2l2e6}) into (\ref{c2l2e1}). We have
\[
\begin{array}{l c l}
F^{ij}\Big\{v_{\nu,ij}+\Omega_{ij}^{rs}\frac{\partial^2 v}{\partial y_r\partial y_s}
+\Theta_{ij}^l\frac{\partial v}{\partial y_l}+\frac{1-t}{n-2}\Delta (v_{\nu})\delta_{ij}
+\frac{1-t}{n-2}\Omega_{kk}^{rs}\frac{\partial^2 v}{\partial y_r\partial y_s}\delta_{ij}&&\\
+\frac{1-t}{n-2}\Theta_{kk}^l\frac{\partial v}{\partial y_l}\delta_{ij}
+(2-t)v_k (-a_k^r\frac{\partial^2 v}{\partial y_n\partial y_r}-\frac{\partial a_k^r}{\partial y_n}\frac{\partial v}{\partial y_r} )\delta_{ij}&&\\
-2 v_i(-a_j^r\frac{\partial^2 v}{\partial y_n\partial y_r}
-\frac{\partial a_j^r}{\partial y_n}\frac{\partial v}{\partial y_r} )-(A^t_g)_{ij,\nu}\Big\}&&\\
=e^{2v}(\phi_{\nu}+2\phi v_{\nu}),&&\\
\end{array}
\]
which can be written as
\begin{equation}
\begin{array}{l c l}
L^{ij}(v_{\nu})_{ij}+\Lambda^{rs}\frac{\partial^2 v}{\partial y_r\partial y_s}=\Pi,&&\\
\end{array}
\label{c2l2e7}
\end{equation}
where
\[
L^{ij}=F^{ij}+\frac{1-t}{n-2}\sum\limits_l F^{ll}\delta^{ij},
\]
\[
\begin{array}{l c l}
\Lambda^{rs}&=&F^{ij}\Big\{\Omega_{ij}^{rs}+\frac{1-t}{n-2}\Omega_{kk}^{rs}\delta_{ij}-(2-t)v_k
a_k^r\delta_{ij}\delta^{sn}+2 v_i a_j^r\delta^{s n}\Big\},\\
\end{array}
\]
and
\[
\begin{array}{l c l}
\Pi &=& e^{2v}(\phi_{\nu}+2\phi v_{\nu})-F^{i j}\Big\{\Theta_{ij}^l\frac{\partial v}{\partial y_l}
+\frac{1-t}{n-2}\Theta_{kk}^l\frac{\partial v}{\partial y_l}\delta_{ij}\\
&&-(2-t)\frac{\partial a_k^r}{\partial y_n}\frac{\partial
v}{\partial y_r} v_k\delta_{i j}
+2\frac{\partial a_j^r}{\partial y_n}\frac{\partial v}{\partial y_r} v_i-(A^t_g)_{ij,\nu}\Big\}\\
\end{array}
\]
depend only on $(M^n,g)$, $|\nabla v|_{C^1(M^,g)}$, $t$, and $\phi$, and are $C^3$ and bounded by $C\sum\limits_l F^{ll}$.

For $\frac{\partial^2 v}{\partial y_r\partial y_s}$, we need to replace it by the partial derivatives of $v$ w.r.t. $e_i$. Recall that $e_i=a_i^j\frac{\partial}{\partial y_j}$. Hence $\frac{\partial}{\partial y_i}=b_i^j e_j$ with $(b_i^j)=(a_i^j)^{-1}$, which is also smooth in $B^T_{\frac{\delta_{y^{i_0}}}{32}}(y^{i_0})$. In $B^T_{\frac{\delta_{y^{i_0}}}{32}}(y^{i_0})$,
\[
\begin{array}{l c l}
\frac{\partial^2 v}{\partial y_r\partial y_s}&=&\nabla^2 v(\frac{\partial}{\partial y_s},\frac{\partial}{\partial y_r})+\Gamma_{rs}^l\frac{\partial v}{\partial y_l}=b_s^i b_r^j \nabla^2 v(e_i,e_j)+\Gamma_{rs}^l\frac{\partial v}{\partial y_l}\\
&=&b_s^i b_r^j v_{ij}+\Gamma_{rs}^l\frac{\partial v}{\partial y_l},\\
\end{array}
\]
therefore (\ref{c2l2e7}) implies that
\[
\begin{array}{l c l}
L^{ij}(v_{\nu})_{ij}+\Lambda^{rs}b_s^i b_r^j v_{i j}&=&\Pi-\Lambda^{rs}\Gamma_{rs}^l\frac{\partial v}{\partial y_l},\\
\end{array}
\]
or
\begin{equation}
\begin{array}{l c l}
L^{i j}(v_{\nu})_{i j}+\sum\limits_{j=1}^n\Lambda^{rs}b_s^n b_r^j v_{\nu j}+\sum\limits_{i=1}^{n-1}\Lambda^{rs}b_s^i b_r^n v_{i \nu}&=&\Pi-\Lambda^{rs}\Gamma_{rs}^l\frac{\partial v}{\partial y_l}-\sum\limits_{i,j=1}^{n-1}\Lambda^{rs}b_s^i b_r^j v_{i j}.\\
\end{array}
\label{c2l2e8}
\end{equation}

By
\[
\begin{array}{l c l}
v_{\nu j}&=&(\nabla_{e_j}\nabla_{\nu}) (v)-(\nabla_{e_j}^{\nu}) (v)=v_{\nu,j}-(\nabla_{e_j}^{\nu})(v),\\
\end{array}
\]
and
\[
\begin{array}{l c l}
v_{i\nu}&=&v_{\nu i}=v_{\nu,i}-(\nabla_{e_i}^{\nu})(v),\\
\end{array}
\]
(\ref{c2l2e8}) implies that
\[
\begin{array}{l c l}
L^{i j}(v_{\nu})_{i j}+\sum\limits_{j=1}^n\Lambda^{rs}b_s^n b_r^j v_{\nu, j}+\sum\limits_{i=1}^{n-1}\Lambda^{rs}b_s^i b_r^n v_{\nu,i}=\Pi-\Lambda^{rs}\Gamma_{rs}^l\frac{\partial v}{\partial x_l}&&\\
+\sum\limits_{j=1}^n\Lambda^{rs}b_s^n b_r^j (\nabla_{e_j}^{\nu})(v)
+\sum\limits_{i=1}^{n-1}\Lambda^{rs}b_s^i b_r^n (\nabla_{e_i}^{\nu})(v)-\sum\limits_{i,j=1}^{n-1}\Lambda^{rs}b_s^i b_r^j v_{i j}.&&\\
\end{array}
\]

Define an elliptic 2nd order linear differential operator in
$B^T_{\frac{\delta_{y^{i_0}}}{32}}(y^{i_0})$ as follows.
\[
L(w)=L^{i j}w_{i j}+(b^i-\bar{s}\sum\limits_l F^{ll}\delta^{n i}) w_i,
\]
where $b^i=\{
\begin{array}{l c l}
\Lambda^{rs}b_s^n b_r^i+\Lambda^{rs}b_s^i b_r^n\quad&\mbox{if}&~~1\le i\le n-1\\
\Lambda^{rs}b_s^n b_r^n\quad&\mbox{if}&~~i=n,\\
\end{array}
$ and $\bar{s}>0$ is some constant to be  determined later. Then
$|b^i|\le C\sum\limits_l F^{ll}$ in
$B^T_{\frac{\delta_{y^{i_0}}}{32}}(y^{i_0})$, and
\begin{equation}
\begin{array}{lcl}
L(v_{\nu})&=&\Pi-\Lambda^{rs}\Gamma_{rs}^l\frac{\partial v}{\partial x_l}+\sum\limits_{j=1}^n\Lambda^{rs}b_s^n b_r^j (\nabla_{e_j}^{\nu})(v)
+\sum\limits_{i=1}^{n-1}\Lambda^{rs}b_s^i b_r^n (\nabla_{e_i}^{\nu})(v)\\
&&-\sum\limits_{i,j=1}^{n-1}\Lambda^{rs}b_s^i b_r^j v_{i j}-\bar{s}\sum\limits_l F^{ll}v_{\nu,n}\\
&\le& C\sum\limits_l F^{ll}-\sum\limits_{i,j=1}^{n-1}\Lambda^{rs}b_s^i b_r^j v_{i j}-\bar{s}\sum\limits_l F^{ll}v_{\nu,\nu}\\
&\le& C\sum\limits_l F^{ll}+C\sum\limits_l F^{ll}\sum\limits_{i,j=1}^{n-1} |v_{i j}|-\bar{s}\sum\limits_l F^{ll}\Big(v_{\nu\nu}+(\nabla_{\frac{\partial}{\partial\nu}}^{\frac{\partial}{\partial\nu}})v\Big)\\
&\le& C\sum\limits_l F^{ll}+C\sum\limits_l F^{ll}\sum\limits_{i,j=1}^{n-1} |v_{i j}|-\bar{s}\sum\limits_l F^{ll} v_{\nu\nu}\quad\mbox{by}~(\ref{c2tn1})\\
&\le& C\sum\limits_l F^{ll}+C\sum\limits_l F^{ll}\sum\limits_{i,j=1}^{n-1} 2(2v_{\nu\nu}+C)-\bar{s}\sum\limits_l F^{ll} v_{\nu\nu}\quad\mbox{by}~(\ref{c2l1r1e1})\\
&\le& C\sum\limits_l F^{ll}+C\sum\limits_l F^{ll}v_{\nu\nu}-\bar{s}\sum\limits_l F^{ll} v_{\nu\nu}\\
&=& C\sum\limits_l F^{ll}-(\bar{s}-C)\sum\limits_l F^{ll} v_{\nu\nu}\\
&\le& C\sum\limits_l F^{ll}\quad\mbox{by taking}~\bar{s}>C,\\
\end{array}
\label{c2l2e9}
\end{equation}
where, in the last inequality, we used $v_{\nu\nu}>-C$, therefore $-v_{\nu\nu}<C$.

From the equation $F(\bar {W})=\phi e^{2v}$, we know
\[
L^{ij}v_{ij}=\phi e^{2v}+F^{ij}v_iv_j-\frac{2-t}{2}|\nabla
v|_g^2\sum\limits_l F^{ll}+F^{ij} (A_g^t)_{ij},
\]
hence
\begin{equation}
|L (v)|\le C\sum\limits_l F^{ll}\quad\mbox{in}~B^T_{\frac{\delta_{y^{i_0}}}{32}}(y^{i_0}).
\label{c2l2e10}
\end{equation}

For any $y_0\in B^T_{\frac{\delta_{y^{i_0}}}{64}}(y^{i_0})\cap\partial M$, let $(a_1,\cdots,a_{n-1},0)$ be the tubular neighborhood normal coordinates of $y_0$ at $y^{i_0}$, and let
\[
D:=\{(y_1,\cdots,y_n)|~y_n\ge
0,~~\sqrt{(y_1-a_1)^2+\cdots+(y_{n-1}-a_{n-1})^2+y_n^2}<\frac{\delta_{y^{i_0}}}{64}\}.
\]

Then

\[
\begin{array}{lcl}
\sqrt{y_1^2+\cdots+y_n^2}&\le&\sqrt{a_1^2+\cdots+ a_{n-1}^2}+\sqrt{(y_1-a_1)^2+\cdots+(y_{n-1}-a_{n-1})^2+y_n^2}\\
&<&\frac{\delta_{y^{i_0}}}{64}+\frac{\delta_{y^{i_0}}}{64}=\frac{\delta_{y^{i_0}}}{32},\\
\end{array}
\]
i.e., $D\subset B^T_{\frac{\delta_{y^{i_0}}}{32}}(y^{i_0})$.

Extend $h_g,\psi$ to a smooth and $C^{3,\alpha}$ function in
$B^T_{\frac{\delta_{y^{i_0}}}{32}}(y^{i_0})$ independently, still
denoted by $h_g,\psi$. In $D$, consider
\[
\bar{w}=v_{\nu}-\psi e^v +h_g + a(1-e^{-b y_n})+\bar{c}
\Big((y_1-a_1)^2+\cdots+(y_{n-1}-a_{n-1})^2+y_n^2\Big),
\]
where $a,b,\bar{c}$ are positive constants to be determined later.

Pick $\bar{c}>0$ such that
\[
v_{\nu}-\psi e^v+h_g+\bar{c}\Big(\frac{\delta_{y^i}}{64}\Big)^2\ge
0\quad\mbox{in}~~B^T_{\frac{\delta_{y^{i_0}}}{32}}(y^{i_0}).
\]

Then
\begin{equation}
\bar{w}(x_0)=0\quad\mbox{and}\quad \bar{w}\ge 0\quad\mbox{on}~~\partial D.
\label{c2l2e12}
\end{equation}

Denote ${\cal R}=(y_1-a_1)^2+\cdots+(y_{n-1}-a_{n-1})^2+y_n^2$. By
(\ref{c2l2e10}),
\begin{equation}
\begin{array}{l c l}
|L(-\psi e^v +h_g +\bar{c} {\cal R})|&=&|L(h_g+\bar{c} {\cal R})-\psi e^v L(v) -e^v L(\psi)\\
&&-2 L^{i j}\psi_i v_j-\psi e^v L^{i j} v_i v_j|\le C\sum\limits_l F^{ll}.\\
\end{array}
\label{c2l2e11}
\end{equation}

To estimate $L(e^{-by_n})$, recall $e_i=a_i^j\frac{\partial }{\partial y_j}$. In $D$,

\[
|b^i(e^{-b y_n})_i|=|b^i a_i^j\frac{\partial }{\partial y_j}(e^{-b y_n})|
=|-b b^i a_i^j e^{-b y_n}\delta_{j n}|\le C b e^{-b y_n}\sum\limits_l F^{ll},
\]
where and in the following, $C>0$ denotes a universal constants
independent of $a$ and $b$.
\[
\begin{array}{l c l}
L^{ij} (e^{-b y_n})_{ij}&=&L^{ij}\Big(\nabla^2 (e^{-b y_n}) (e_i,e_j)\Big)=a_i^r a_j^s L^{i j}\Big(\nabla^2 (e^{-b y_n})(\frac{\partial}{\partial y_r},\frac{\partial}{\partial y_s})\Big)\\
&=& a_i^r a_j^s L^{i j}\Big( \frac{\partial^2}{\partial y_s\partial y_r}(e^{-b y_n})-\Gamma_{sr}^l\frac{\partial }{\partial y_l}(e^{-b y_n})\Big)\\
&=&a_i^r a_j^s  L^{i j}\Big(b^2 e^{-b y_n}\delta_{rn}\delta_{sn}+b e^{-b y_n}\Gamma_{sr}^l\delta_{nl} \Big)\\
&=&b^2 e^{-b y_n} L^{ij}a_i^n a_j^n +b e^{-b y_n}\Gamma_{sr}^n L^{ij} a_i^r a_j^s\\
&\ge& b^2 e^{-b y_n} L^{ij}a_i^n a_j^n -Cb e^{-b y_n}\sum\limits_l F^{ll} \\
&\ge& b^2 e^{-b y_n} a_i^n a_j^n (F^{ij}+\frac{1-t}{n-2}\sum\limits_l F^{ll}\delta_{ij}) -Cb e^{-b y_n}\sum\limits_l F^{ll} \\
&\ge& \frac{1-t}{n-2}b^2 e^{-b y_n} (a_i^n)^2\sum\limits_l F^{ll}   -Cb e^{-b y_n}\sum\limits_l F^{ll} \\
&=& \frac{1-t}{n-2}b^2 e^{-b y_n} \sum\limits_l F^{ll}   -Cb e^{-b y_n}\sum\limits_l F^{ll} \quad\mbox{by}~~(\ref{c2l2e3}).\\
\end{array}
\]

Thus in $D$,
\[
\begin{array}{l c l}
L(e^{-b y_n})&\ge&\frac{1-t}{n-2}b^2 e^{-b y_n} \sum\limits_l F^{ll}   -Cb e^{-b y_n}\sum\limits_l F^{ll}\\
&\ge&b e^{-b y_n}\Big(\frac{(1-t)b}{n-2}-C\Big)\sum\limits_l F^{ll}\\
&\ge& b e^{-b y_n}\sum\limits_l F^{ll},\\
\end{array}
\]
by choosing $b\gg 1$ such that $\frac{(1-t)b}{n-2}-C\ge 1$.

Back to $L(\bar {w})$, we have in $D$,
\[
\begin{array}{l c l}
L(\bar{w})&=& L(v_{\nu}-\psi e^v +h_g +c {\cal R})-a L(e^{-b y_n})\\
&\le& -ab e^{-b y_n}\sum\limits_l F^{ll}+C\sum\limits_l F^{ll}\le 0,\\
\end{array}
\]
by choosing $a\gg 1$ such that $a b e^{-\frac{b\delta_{y^i}}{64}}>C$. 

Hence (\ref{c2l2e12}) implies that
\[
\bar {w}\ge 0\quad\mbox{in}~~D,
\]
therefore we have $\bar {w}_{\nu}(y_0)\le 0$, i.e., $v_{\nu\nu}(y_0)<C$. Since $y_0\in B^T_{\frac{\delta_{y^{i_0}}}{64}}(y^{i_0})\cap\partial M$ is arbitrary, Lemma~\ref{c2l2} has been established. $\clubsuit$
\vskip 5 pt

\begin{remark}
By the Lemma~\ref{c2l1} and the Lemma~\ref{c2l2} and $\cup_{i_0=1}^N (B^T_{\frac{\delta_{y^{i_0}}}{64}}(y^{i_0})\cap\partial M)=\partial M$, we know the Hessian of $v$ on $\partial M$ is upper bounded w.r.t. the metric $g$. Thus $\Gamma\subset\Gamma_1$ implies that
\[
|\nabla^2_g v|_g\le C\qquad\mbox{on}~~\partial M.
\]
\label{c2l2r1}
\end{remark}
\begin{lemma}
Under the same assumptions as in Theorem~\ref{generalresult}, for $t<1$, let $v$ be a $C^4$ solution of the equation (\ref{generalproblem}). Then there exists a universal constant $C>0$ depending only on $(M^n,g,t)$, $(f,\Gamma)$, $\phi$, $\psi$ such that on $M$,
\[
|\nabla ^2 v|<C.
\]
\label{c2l3}
\end{lemma}
{\bf Proof of the Lemma~\ref{c2l3}.} Consider
\[
E(x)=\max\limits_{e\in T_x M,~g(e,e)=1}(\nabla^2 v+a|\nabla v|^2_g g)(e,e).
\]

Let $E(x_0)=\max\limits_M E$, and let $\{x_j\}_{j=1}^n$ be a geodesic normal coordinates w.r.t. the metric $g$ at $x_0$. In the following, subindices denote the covariant derivatives w.r.t. $\frac{\partial}{\partial x_j}$. W.l.o.g, we assume $x_0$ is an interior point of $M$, and $E(x_0)=v_{11}+a|\nabla v|^2_g$. Consider $\bar{E}=\frac{v_{11}}{g_{11}}+a|\nabla v|^2_g$. Then $x_0$ is a local maximum point of $\bar{E}$. We can proceed as in the proof of the Lemma~\ref{c2l1} to finish the proof of the Lemma~\ref{c2l3}. $\clubsuit$
\vskip 5 pt
\begin{section}
{Proof of the Theorem~\ref{generalresult}}
\end{section}
Consider the homotopy equation $H_s$, for $0\le s\le 1$,
\begin{equation}
\Big\{
\begin{array}{l c l}
f\Big(-\lambda_g(s\bar {W}+(1-s)\sigma_1(\bar{W})g)\Big)-s\phi e^{2v}-(1-s)e^{2v}=0\quad\mbox{on}~~M,&&\\
v_{\nu}+h_g-s e^v\psi=0\quad\mbox{on}~~\partial M,&&\\
\end{array}
\label{existe1}
\end{equation}
where $\bar{W}=W_g^v-A_g^t$.

By the uniform $C^2$ estimates we established and the result of Lieberman and Trudinger (\cite{LT}), we have
the uniform $C^{2,\alpha_0}$ bounds for the solutions of the above equation. $C^{4,\alpha_0}$ estimates follow from the Schauder estimates. By the direct computation, the linearized operator ${\cal L}_s(w)$ at a solution $v$ is given by
\begin{equation}
\Big\{
\begin{array}{l c l}
\Big(s L^{i j}+(1-s)L^{ll}\delta_{i j}\Big)w_{i j}+\bar {b}^i w_i-2(s\phi+(1-s))e^{2v}w\quad\mbox{on}~~M,&&\\
w_{\nu}-s\psi e^v w\quad\mbox{on}~~\partial M,&&\\
\end{array}
\end{equation}
where
\[
\bar {b}^i=s(2-t)F^{ll}v_i-2 s F^{ij}v_j+(2n-nt-2)(1-s)F^{ll} v_i.
\]

By $\phi>0$ $\psi\le 0$ and the maximum principle, the linearized operator is an elliptic invertible operator: $C^{2,\alpha}\to C^{\alpha}$. Hence the equation of (\ref{existe1}) for $s=1$ is uniquely
solvable in $C^{4,\alpha_0}$ if and only if the equation of (\ref{existe1}) for $s=0$ is uniquely solvable in $C^{4,\alpha_0}$. When $s=0$, the uniqueness and the existence of the solution has been confirmed in \cite{Cher}. $\clubsuit$
\begin{section}
{Proof of the Theorem~\ref{existence}}
\end{section}
Take an arbitrary Riemannian metric $g$ on $M^n$. For instance, let $\{U_i,x^{(i)}_j\}_{i=1,j=1}^{N,~~~n}$ be a finite coordinate charts on $M^n$ and let $\phi^i$ be a partition of unity subordinate to $U_i$. We can simply take $g$ to be $\sum\limits_{i=1}^{N}\phi^i\big((dx^{(i)}_1)^2+\cdots+(dx^{(i)}_n)^2\big)$. Let $w(x)$ be a smooth function on $M^n$ such that $w(x)$ is the distance of $x$ to $\partial M$ w.r.t. the metric $g$ when $x$ is near $\partial M$. Then $\frac{\partial w}{\partial\nu}|_{\partial M}=-1$, where
$\frac{\partial}{\partial\nu}$ is the unit outer normal of $g$ on $\partial M$. Extend the mean curvature $h_g$ to a smooth function defined on $M^n$, still denoted by $h_g$. We can obtain such extension by straightening the boundary and extending any function $\bar{\psi}$ defined on $\partial\Bbb{R}^n_+$ to $\Bbb{R}^n_+$ using $\bar{\psi}(x')(1-x_n)$, where $x=(x',x_n)\in\Bbb{R}^n_+$. However, we want to mention a different way which seems more natural. In fact, we only need to extend $h_g$ smoothly to the interior of $M$ near $\partial M$. Using the partition of unity, we can localize the extension to a small neighborhood of each $x_0\in\partial M$. Notice $h_g$ is the trace of the second fundamental form of $g$ on $\partial M$ whose definition is, at every point $x\in\partial M$,
\[
II(X,Y)=-g(\nabla_X^{\frac{\partial}{\partial \nu}},Y),\quad\forall X,Y\in T_{x}(\partial M).
\]

Let $U$ be a small neighborhood where the tubular neighborhood normal coordinates of $x\in U$ at $x_0$ is smooth and well-defined. Let $\{x_j\}_{j=1}^n$ be such coordinates. Then $-\frac{\partial}{\partial x_n}$ is a smooth extension of $\frac{\partial}{\partial\nu}$ to $U$ and $g(\frac{\partial}{\partial x_k},\frac{\partial}{\partial x_n})=0$ in $U$ for $1\le k\le n-1$, so
\begin{equation}
\frac{1}{n-1}\sum\limits_{i,j=1}^{n-1} g(\nabla_{\frac{\partial}{\partial x_i}}^{\frac{\partial}{\partial x_n}},\frac{\partial}{\partial x_j})g^{ij}
\label{exe1}
\end{equation}
is an extension of $h_g$ to $U$, where $(g^{ij})$ is the inverse of $(g_{ij})=(g(\frac{\partial}{\partial x_i},\frac{\partial}{\partial x_j}))$. From the linear algebra, we know $(g^{ij})=\frac{1}{\det(g_{ij})}adj(g_{ij})$, hence $g^{ij}$ is smooth and (\ref{exe1}) gives a smooth extension of $h_g$ to $U$.

Let $v=h_g w$. Consider $g_1=e^{2v}g$. Then on $\partial M$,

\begin{eqnarray*}
h_{g_1}&=&(\frac{\partial v}{\partial\nu}+h_g)e^{-v}=(h_g\frac{\partial w}{\partial\nu}+w\frac{\partial h_g}{\partial \nu}+h_g)e^{-v}\\
&=&\big(h_g(-1)+(0)\frac{\partial h_g}{\partial \nu}+h_g\big)=0.\\
\end{eqnarray*}

For $g_1$, let $w_1$ be a smooth function such that, near $\partial M$,  $w_1$ is the distance function to $\partial M$ w.r.t. $g_1$. We know that $\frac{\partial w_1}{\partial\nu_1}|_{\partial M}=-1$,
where $\frac{\partial}{\partial\nu_1}$ is the unit outer normal of $g_1$ on $\partial M$. Take $g_2=e^{2A(w_1)^2}g_1$ with $A>0$ being a constant to be chosen later.

Direct computations yield that on $\partial M$
\[
h_{g_2}=(2 A w_1\frac{\partial w_1}{\partial\nu_1}+h_{g_1})e^{A(w_1)^2}=0,
\]
and
\begin{equation}
\begin{array}{lcl}
Ric_{g_2}&=&Ric_{g_1}-(n-2)A\nabla^2_{g_1}(w_1^2)-A(\Delta_{g_1}(w_1^2))g_1+(n-2)A^2 d(w_1^2)\otimes d(w_1^2)\\
&&-(n-2) A^2|\nabla (w_1^2)|_{g_1}^2 g_1\\
&\le& Ric_{g_1}-A(n-2)\nabla^2_{g_1}(w_1^2)-A(\Delta_{g_1}(w_1^2))g_1,\\
\end{array}
\label{ricci}
\end{equation}
where in the last inequality, we used a general fact that $df\otimes df\le |\nabla f|^2_{g_1} g_1$ for any $C^1$ function $f$. The explanation is given as follows. At each $x$, we take a geodesic normal coordinates $\{x_i\}_{i=1}^n$ of $g_1$ at $x$. At $x$,
\[
\begin{array}{lcl}
(df\otimes df)(\frac{\partial}{\partial x_k},\frac{\partial}{\partial x_k})&=&(\frac{\partial f}{\partial x_k})^2\le \sum\limits_{i=1}^n(\frac{\partial f}{\partial x_i})^2= (|\nabla f|^2_{g_1} g_1)(\frac{\partial}{\partial x_k},\frac{\partial}{\partial x_k}),\\
\end{array}
\]
which implies that $df\otimes df\le |\nabla f|^2_g g$ since both $df\otimes df$ and $|\nabla f|^2_g g$ are symmetric $(0,2)$ tensors.

For any $x_0\in\partial M$, we take a tubular neighborhood normal coordinates $\{x_j\}_{j=1}^n$ of $g_1$ at $x_0$. Then $w_1=x_n$ near $x_0$. At $x_0$, by $x_n=0$
\begin{eqnarray*}
\nabla^2_{g_1}[( x_n)^2](\frac{\partial}{\partial x_i},\frac{\partial}{\partial x_j})&=&\nabla_{\frac{\partial}{\partial x_j}}\nabla_{\frac{\partial}{\partial x_i}}[(x_n)^2]-(\nabla_{\frac{\partial}{\partial x_j}}^{\frac{\partial}{\partial x_i}})[(x_n)^2]\\
&=&\nabla_{\frac{\partial}{\partial x_j}}\nabla_{\frac{\partial}{\partial x_i}}[(x_n)^2]-2x_n(\nabla_{\frac{\partial}{\partial x_j}}^{\frac{\partial}{\partial x_i}})[x_n]\\
&=&\nabla_{\frac{\partial}{\partial x_j}}\nabla_{\frac{\partial}{\partial x_i}}[(x_n)^2]
=\nabla_{\frac{\partial}{\partial x_j}}[2 x_n \nabla_{\frac{\partial}{\partial x_i}}(x_n)]\\
&=&2(\nabla_{\frac{\partial}{\partial x_j}}x_n)(\nabla_{\frac{\partial}{\partial x_i}}x_n )+2 x_n \nabla_{\frac{\partial}{\partial x_j}}\nabla_{\frac{\partial}{\partial x_i}}x_n\\
&=&2(\nabla_{\frac{\partial}{\partial x_j}}x_n)(\nabla_{\frac{\partial}{\partial x_i}}x_n )=2\delta_j^n\delta_i^n,\\
\end{eqnarray*}
so at $x_0$,
\[
\nabla_{g_1}^2[w_1^2]=\nabla_{g_1}^2[(x_n)^2]=2 d x_n\otimes d x_n \ge 0,\\
\]
and
\[
\Delta_{g_1}[w_1^2]=\Delta_{g_1}[(x_n)^2]=2.
\]

Substitute the above two into (\ref{ricci}). At $x_0$, we have
\[
\begin{array}{lcl}
Ric_{g_2}&\le& Ric_{g_1}-(n-2)A\nabla^2_{g_1}(w_1^2)-A(\Delta_{g_1}(w_1^2))g_1\\
&\le& Ric_{g_1}-2 A g_1\le C_1 g_1-2 A g_1,\\
\end{array}
\]
where $C_1>0$ is a universal constant depending only on $(M^n,g)$ and independent of $x_0$.

Choose $A\ge \frac{C_1}{2}+\frac{1}{2}$. Then $Ric_{g_2}(x_0)\le -g_1 (x_0)$, which implies that $Ric_{g_2}\le -g_1$ on $\partial M$, hence
\[
Ric_{g_2}<0\quad\mbox{in a tubular neighborhood of}~~\partial M.
\]

By the result in \cite{Loh}, there is a smooth metric $g_3$ on $M$ such that
\[
g_3\equiv g_2\quad\mbox{in a smaller tubular neighborhood of }\partial M,
\]
and
\[
Ric_{g_3}<0\quad\mbox{on}~~M.
\]

Clearly, $h_{g_3}=h_{g_2}=0$ on $\partial M$, and $Ric_{g_3}<0$ on $M$ implies that 
\[
-\lambda_{g_3}(A_{g_3}^t)\in\Gamma_n\subset\Gamma,\quad\forall  t<1.
\]

Thus, by the Theorem~\ref{generalresult}, there exists a unique $C^{4,\alpha_0}$ metric $g_4\in [g_3]$ solving
\[
\{
\begin{array}{lcl}
f(-\lambda_{g_4}(A^t_{g_4}))&=&\phi,\quad -\lambda_{g_4}(A^t_{g_4})\in\Gamma\quad\mbox{on}~~M\\
h_{g_4}&=&\psi\quad\mbox{on}~~\partial M.
\end{array}
\]

In particular, we can take $(f,\Gamma)=(\sigma_n^{\frac 1n},\Gamma_n)$, $t=0$, and $\phi\equiv 1$, $\psi\equiv 0$. Theorem~\ref{existence} has been established. $\clubsuit$
\vskip 5pt

From the arguments in the proof of Theorem~\ref{existence}, it is easy to see that for any smooth compact Riemannian manifold $(M^n, g)$ ($n\ge 3$) with some boundary including those metrics with positive Ricci tensors, there exists some metric $g_3$ which is conformal to $g$ near $\partial M$ satisfying
\[
-\lambda_{g_3}(A_{g_3}^t)\in\Gamma_n\subset\Gamma\quad{on}~~M\quad\mbox{and}\quad  h_{g_3}=0\quad\mbox{on}~~\partial M.
\]
Thus we have the following result
\begin{theorem}
Let $(M^n,g)$ be an $n-$dimensional ($n\ge 3$) compact smooth Riemannian manifold with
$\partial M\neq\emptyset$ and let $f\in C^{2,\alpha_0}(\Gamma)$ ($0<\alpha<1$) satisfy (\ref{cone})-(\ref{f2}). Given
$0<\phi\in C^{2,\alpha_0}(M^n)$, $0\ge \psi\in C^{3,\alpha_0}(\partial M)$ and for
any $t<1$, there exists
a $C^{4,\alpha_0}$ solution $\tilde {g}$ which is conformal to $g$ near $\partial M$ and solves
\[
\{
\begin{array}{lcl}
f(-\lambda_{\tilde{g}}(A^t_{\tilde{g}}))&=&\phi,\quad -\lambda_{\tilde{g}}(A^t_{\tilde{g}})\in\Gamma\quad\mbox{on}~~M\\
h_{\tilde{g}}&=&\psi\quad\mbox{on}~~\partial M.
\end{array}
\]
\end{theorem}


\begin{thebibliography}{99}
\bibitem{Aub} T. Aubin, $\acute{E}$quations diff$\acute{e}$rentielles non lin$\acute{e}$aires et probl$\grave{e}$me de Yamabe concernant la courbure scalaire, J. Math. Pures Appl. (9) 55 (1976), 269-296.
\bibitem{CGY1} S.Y. A. Chang, M. Gursky and P. Yang, An a priori estimate for a fully nonlinear equation on four-manifolds, J. Anal. Math. 87 (2002), 151-186.
\bibitem{Cher} M. Cheerier, Probl$\grave{e}$mes de Neumann non lin$\acute{e}$aires sur les vari$\acute{e}$rt$\acute{e}$s Riemanniennes, J. Funct. Anal. 57 (1984), 154-206.
\bibitem{CNS} L. Caffarelli, L. Nirenberg and J. Spruck, The Dirichlet problem for nonlinear second-order elliptic equations, III. Functions of the eigenvalues of the Hessian, Acta Math. 155 (1985), 261-301.
\bibitem{Eis} L. P. Eisenhart, Riemannian geometry, Princeton Landmarks in Mathematics, Princeton University Press, Princeton, NJ, 1997, Eighth printing, Princeton Paperbacks.
\bibitem{Esc} J. F. Escobar, The Yamabe problem on manifolds with boundary, J. Differential Geom. 35 (1992), 21-84.
\bibitem{Esc1} J. F. Escobar, Conformal deformation of a Riemannian metric to a constant scalar curvature metric with constant mean curvature on the boundary, Indiana Univ. Math. J. 45 (1996), 917-943.
\bibitem{GT} D. Gilbarg and N.S. Trudinger, Elliptic partial differential equations of second order, reprint of the 1998 edition, Classics in Mathematics. Springer, Berlin Heidelberg New York (2001).
\bibitem{GV} M. Gursky and J. Viaclovsky, Fully nonlinear equations on Riemannian manifolds with negative curvature, Indiana Univ. Math. J. 52 (2003), 399-420.
\bibitem{GV1} M. Gursky and J. Viaclovsky, Prescribing symmetric functions of eigenvalues of the Ricci tensor, Anal. of Math. 166 (2007), 475-531.
\bibitem{GW} P. Guan and G. Wang, Local estimates for a class of fully nonlinear equations arising from conformal geometry, Int. Math. Res. Not. 26 (2003), 1413-1432.
\bibitem{HL1} Z.C. Han and Y.Y. Li, The Yamabe problem on manifolds with boundaries: Existence and compactness results, Duke Math. J. 99 (1999), 489-542.
\bibitem{HL2} Z.C. Han and Y.Y. Li, The existence of conformal metrics with constant scalar curvature and constant boundary mean curvature, Comm. Anal. Geom. 8 (2000), 809-869.
\bibitem{JLL} Q. Jin, A. Li and Y.Y. Li, Estimates and existence results for a fully nonlinear Yamabe problem on manifolds with boundary, Calc. Var. PDEs 28 (2007), 509-543.
\bibitem{Loh} J. Lohkamp, Negatively Ricci curved manifolds, Bull. Amer. Math. Soc. (2) 27 (1992), 288-291.
\bibitem{L} Y.Y. Li, Some existence results of fully nonlinear elliptic equations of Monge-Amp$\grave{e}$re type, Comm. Partial Differential Equations, 14 (1989), 1544-1578.
\bibitem{LL1} A. Li and Y.Y. Li, On some conformally invariant fully nonlinear equations, Part II: Liouville, Harnack and Yamabe. Acta Math. 195 (2005), 117-154.
\bibitem{LL2} A. Li and Y.Y. Li, On some conformally invariant fully nonlinear equations, Comm. Pure Appl. Math. (2003), no. 10, 1416-1464.
\bibitem{LT} G.M. Lieberman and N.S. Trudinger, Nonlinear oblique boundary value problems for nonlinear elliptic equations, Trans. Amer. Math. Soc. 295 (2) (1986), 509-546.
\bibitem{STW}, W. Sheng, N.S. Trudinger and X. Wang, The Yamabe problem for higher order curvatures, J. Diff. Geom. 77 (2007), 515-553.
\bibitem{SY} R. Schoen and S.T. Yau, On the proof of the positive mass conjecture in general relativity, Comm. Math. Phys. 65 (1979), 45-76
\bibitem{Tr} N.S. Trudinger, Remarks concerning the conformal deformation of Riemannian structures on compact manifolds, Ann. Scuola Norm. Sup. Pisa (3) 22 (1968), 265-274.
\bibitem{U} J. Urbas, Hessian equations on compact Riemannian manifolds, in Nonlinear Problems in Mathematical Physics and Related Topics, Vol. II, 367-377. Int. Math. Ser. (N.Y.), 2. Kluwer/Plnum, New York, 2002.
\bibitem{Ya} H. Yamabe, On a deformation of Riemannian structures on compact manifolds, Osaka Math. J. 12 (1960), 21-37.
\end{thebibliography}
\end{document}